\documentclass[a4paper,12pt]{article}
\author{Alan Hammond\thanks{Department of Statistics, Oxford University. Research supported in part by NSF grants DMS-0806180 and OISE-0730136 and by EPSRC grant EP/I004378/1. This research was started at New York University.}}
\title{Stable limit laws for randomly biased walks on supercritical trees}

\usepackage{amsfonts,psfrag,epsfig}
\usepackage{appendix}
\usepackage{url}

\def\P{\mathbb{P}}
\def\Z{\mathbb{Z}}
\def\N{\mathbb{N}}
\def\E{\mathbb{E}}
\def\Q{\mathbb{Q}}
\def\R{\mathbb{R}}
\def\binf{q}
\def\bsup{Q}

\def\coup{\Theta}

\def\pf{\mathbb{P}_{f,\nu,\infty} \times \P_{T,\beta}^{\phi}}

\def\pfba{\mathbb{P}_{f,\nu,\infty} \times \P_{B(T),\beta}^{\phi}}
\def\baug{B_{{\rm aug}}}
\def\bext{B_{\rm ext}}
\def\btbase{{v_{\rm base}}}
\def\olddone{d_2}
\def\cpo{d_1}
\def\cpt{\olddone}

\def\barecon{B_0}
\def\conec{C_3}
\def\entrance{{\rm ent}}
\def\trwtbd{v_0}
\def\wgt{weighted }
\def\unidir{end-distinguished }
\def\unidirnosp{end-distinguished}

\def\head{{\rm head}}
\def\hatc{C_1}
\def\base{{\rm base}}

\def\coupl{\Theta}
\def\basedef{{\rm base}}

\def\vbase{{v_{\rm base}}}

\def\connew{c_1}

\def\concaptwo{C_2}
\def\conai{C_4}
\def\falldeep{\mathcal{DE}}
\def\deepfall{\mathcal{FD}}
\def\ero{\mathcal{E}}

\def\omen{\omega_{\entrance}}
\newcommand{\omenn}[1]{\omega_{#1}}
\def\b{\basedef}
\def\gcd{{\rm GCD}}
\def\singent{\overline{T}}
\newcommand{\qstarnoarg}{\Q}
\def\btmeas{\qstarnoarg \times \P_{(B,\singent),\beta}^\entrance}
\def\btmeasbrac{\big( \qstarnoarg \times \P_{(B,\singent),\beta}^\entrance \big)}
\def\btmeasbracbig{\Big( \qstarnoarg \times \P_{(B,\singent),\beta}^\entrance \Big)}
\def\corrfac{2c_f}
\def\halfcorrfac{c_f}
\def\corrfack{2c_{f,k}}

\def\halfcorrfack{c_{f,k}}
\newcommand{\corrfacarg}[1]{2c_{f,#1}}
\newcommand{\halfcorrfacarg}[1]{c_{f,#1}}
\def\pde{p_{\rm de}}
\def\pesc{p_{\rm esc}}
\def\ekm{\pesc(k)}
\def\qkm{\pde(k)}
\def\constr{\eta}
\def\constt{\psi}
\def\chiexp{\gamma}
\def\chiexpdiff{\kappa}
\def\bout{N^{\rm out}}
\def\ik{N^{\rm in}_k}
\def\okun{N^{\rm out}}
\def\ikun{N^{\rm in}}
\def\ok{N^{\rm out}_k}
\def\okmo{N^{\rm out}_{k-1}}
\def\dk{D_k}
\def\ekd{E_k}
\def\pext{q_{{\rm ext}}}
\def\taustar{\tau}
\newcommand{\trar}[1]{{\rm ar}_{#1}}
\newcommand{\tnb}[1]{{\rm tn}_{#1}}
\newcommand{\overh}[1]{\overline{\mathcal{H}}_{#1}}
\newcommand{\hitting}[1]{\mathcal{H}_{#1}}
\newcommand{\queu}{\Q^{\rm out}}
\newcommand{\qstar}[1]{\Q^{#1}}
\newcommand{\qstarl}[1]{\Q_{#1}}
\newcommand{\qstaru}{\Q^u}
\newcommand{\qstarreal}{\Q}
\newcommand{\qsubk}{\Q_k^{\rm out}}
\newcommand{\qsubuksupstar}{\Q^u_k}
\newcommand{\idtrap}{\tilde\Q}
\newcommand{\qold}{\Q^{\rm out}}
\newcommand{\rur}{\mathbb{R}^{u,r}}
\newcommand{\reg}{{\rm r}}

\newcommand{\llone}{M}
\newcommand{\basedefarg}[1]{\basedef(#1)}
\newcommand{\xsub}[1]{X(#1)}
\newcommand{\sigmasub}[1]{\sigma_{(#1)}}
\newcommand{\alphasub}[1]{\lambda_{(#1)}}
\def\offspr{O}
\newcommand{\largev}{{\rm L}}

\newtheorem{corollary}{Corollary}[section]
\newtheorem{lemma}[corollary]{Lemma}
\newtheorem{definition}[corollary]{Definition}

\newtheorem{theorem}{Theorem}
\newtheorem{prop}{Proposition}
\newtheorem{hyp}{Hypothesis}

\def\build#1_#2^#3{\mathrel{ \mathop{\kern 0pt#1}\limits_{#2}^{#3}}}
\def\simk{\build{\simeq}_{}^{K}}
\def\mulaw{\upsilon}
\def\zetax{\zeta^{(x)}}
\def\msupx{\theta^{(x)}}
\def\squone{1}
\def\squtwo{2}

\newcommand{\qed}{\hfill \ensuremath{\Box}}

\def\fff#1{&{{\pageref{#1}}}\cr}
\def\hfff#1{\label{#1}}

\begin{document}
\maketitle
\begin{abstract}
We consider a random walk on a supercritical Galton-Watson tree with leaves, where the transition probabilities of the walk are determined by biases that are randomly assigned to the edges of the tree. The biases are chosen independently on distinct edges, each one according to a given law that satisfies a logarithmic non-lattice condition. We determine the condition under which the walk is sub-ballistic, and, in the sub-ballistic regime, 
we find a formula for the exponent $\chiexp \in (0,1)$ such that the distance $\vert X(n) \vert$ moved by the walk in time $n$ is of the order of $n^\chiexp$. We prove a stable limiting law for walker distance at late time, proving that the rescaled walk $n^{-\chiexp} \vert X(n) \vert$ converges in distribution to an explicitly identified function of the stable law of index $\chiexp$.

This paper is a counterpart to \cite{BAFGH}, in which it is proved that, in the model where the biases on edges are taken to be a given constant, there is a logarithmic periodicity effect that prevents the existence of a stable limit law for scaled walker displacement. It is randomization of edge-biases that is responsible for the emergence of the stable limit in the present article, while also introducing further correlations into the model in comparison with the constant bias case. The derivation requires the development of a detailed understanding of trap geometry and the interplay between traps and backbone. The paper may be considered as a sequel to \cite{paperone}, since it makes use of a result on the regular tail of the total conductance of a randomly biased subcritical Galton-Watson tree. 
\end{abstract}
\footnotesize{
\tableofcontents}
\newpage
\subsection*{Glossary of notation}
Quite a lot of notation will be needed in the article. For the reader's convenience, we begin by listing much of the notation, and provide a summarizing phrase for each item, as well as the page number at which the concept is introduced.

\vspace{-4mm}

\bigskip
\def\qq{&}

\begin{center}
\halign{
#\quad\hfill&#\quad\hfill&\quad\hfill#\cr
$[\binf, \bsup]$ \qq interval containing all edge-biases \fff{suppbias}
$\beta_e$ \qq the bias of the edge $e \in E(T)$ \fff{betae}
$\P_{T,\beta}$   \qq the law of $\beta$-biased walk in $T$ \fff{ptblaw} 
$\P_f,\P_{f,\infty},\P_{f,\nu,\infty}$ \qq notation for Galton-Watson trees (weight law $\nu$) \fff{gwtree}
$\P_{f,\nu,\infty} \times \P_{T,\beta}^\phi$ \qq the joint law of environment and walk \fff{pfnptb}
$\hitting{n}$ \qq time to reach distance $n$ from root \fff{hitdist}
$\pext$ \qq the extinction probability of $\P_f$ \fff{extprob}
$B,\baug,\bext$ \qq terms related to backbone \fff{threebb}
$P_{v,w}$ \qq the path connecting $v,w \in V(T)$ \fff{pvw}
$T_w$ \qq the descendent tree of $w \in V(T)$ \fff{desctree}
$\head,\entrance$ \qq the root and its offspring in a single-entry tree \fff{headent}
$H_w$ \qq hitting time of $w \in V(T)$ \fff{hitw}
$\omega_{x}(y)$ \qq weight of $P_{x,y}$ \fff{omegaxy}
$\omega(T)$ \qq the sum of weights in $V(T)$ \fff{omegat}
$\vbase$ \qq specific vertex attaining depth of $T$ \fff{baset}
$\barecon$-bare \qq notion of ``thin'' tree \fff{baretree}
$\falldeep,\pde$ \qq the event of deep excursion into tree and its probability \fff{evfalldeep}
$\trar{i}$ \qq time of arrival at the $i$-th trap entrance \fff{trapenti}
\textrm{end-distinguished} \qq tree with root infinitely far in the past \fff{edtree}
$(B,\singent)$ \qq backbone-tree pair \fff{btpair}
$\tau_{T_{\entrance}}$ \qq total trap time for walk in $(B,\singent)$ \fff{tautent}
$\pesc$ \qq escape probability for backbone-tree pair \fff{probesc}
$\halfcorrfac$ \qq the correction factor accounting for total trap time \fff{correctionfactor}
$\deepfall$ \qq event of falling deeply into the trap in a backbone-tree pair \fff{evdeepfall}
$\bout_k(u) \in \Psi_k$ \qq $k$-large backbone neighbourhood of $u \in \bext$ \fff{nkoutpsi} 
$\qstarnoarg$ \qq late-time arrival law on backbone-tree pairs \fff{ltabt}
${\rm RG}$ \qq regeneration-time set for walk \fff{regwalk}
$\big\{ C_i: 1 \leq i \leq r(T) \big\}$ \qq components of renewal decomposition of finite tree \fff{renewaldecomp}
$\ik$ \qq the trap interior up to the $k$-th cutpoint $c_k(T_\entrance)$ \fff{bink} 
$\halfcorrfack$ \qq the approximating correction factor \fff{approxcf}
$\Psi_k^+$ \qq state-space of $k$-large neighbourhoods $(\ik,\ok)$  \fff{psikpl}
$\Psi_k^+[K]$ \qq $\big\{ \xi \in \Psi_k^+: \ik(\xi)  \textrm{ contains at most $K$ vertices} \big\}$ \fff{psikplk}
$\qstaru$ \qq  $\qstarnoarg \big( \cdot \big\vert \corrfac \omen(T_{\entrance}) > u \big)$ \fff{qstaru}
$\tau_i$ \qq total time in $i$-th trap under $\pf$ \fff{ttt}
$\overh{n}$ \qq time to reach distance $n$ from root while in backbone \fff{distrootb}
$\tnb{n}$ \qq number of trap entrances encountered before $\overh{n}$ \fff{trenen}
$\idtrap$     \qq idealized sequence of $\qstarnoarg$-distributed backbone-tree pairs \fff{idtrap}
$(\tilde{B}^i,\tilde{T}^i)$,$\tilde{X}_i$,$\tilde{\tau}_i$ \qq environment, walk and total trap time in $i$-th element of  $\idtrap$ \fff{idtrapnot}

\cr}\end{center}
\normalsize
\begin{section}{Introduction}
A thermally agitated particle subject to a constant external force in $\R^3$ may be modelled by a 
random walk with a non-zero drift. Under minor assumptions on the jump law, (such as, that the jumps are bounded),  the walk  will move at a deterministic velocity. If, however, the particle resides in a disordered medium, in which impenetrable objects are present, then the regular progress of the particle is not assured. Indeed, it was argued physically in~\cite{BB} that trapping effects would result in a particle velocity that was not a monotone function of the bias. Later, 
it was further argued (in \cite{Dhar} and \cite{DS}) that strong biases would result in sub-ballistic regimes, in which the particle has asymptotically zero velocity. The latter prediction has been verified for biased random walks on random trees by  Lyons, Pemantle and Peres  \cite{LPP96}.  
In the problem of biased random walk on the supercritical percolation cluster of $\Z^d$, the prediction of zero velocity was verified at sufficiently high values of the bias by \cite{BGP} (for $d=2$) and \cite{Sznitman1} ($d \geq 2$), while the existence of a critical bias delimiting the sub-ballistic and ballistic regimes has been proved in \cite{AlexAlan} (with $d \geq 2$).

Related to this slowdown is the phenomenon of aging, in which a random system becomes trapped in deep wells, spending a time there proportional to the age of the system. 
It has been predicted \cite{Bouchaud} that aging describes the dynamics of low temperature spin-glasses, and it has been proved \cite{BACerny} to arise in effective models, in which random walker jumps are time-changed at a random rate associated to walker location. 
A natural aim is to seek to prove aging in systems where trapping occurs not by being built into the model's definition, but arises naturally in long-time dynamics. 

In this article, we prove aging for biased random walk whose environment is a supercritical Galton-Watson tree having positive extinction probability. This forces us to address in detail the geometry of traps in which the walk falls at advanced time. In the model we study, the bias of the random walk is randomized over edges, so that the jump law at distinct vertices is independent and identically distributed. This random biasing is introduced to permit aging to occur.

We start by defining the general class of biased random walks on trees in which we are interested.
\begin{definition}\label{deftree}
Two constants $\bsup \geq \binf > 1$ are fixed.  \hfff{suppbias}
Let $T$ be a rooted and locally finite tree, with vertex set $V(T)$, edge set $E(T)$, and root $\phi$.
The tree $T$ is said to be weighted if it carries a function $\beta:E(T) \to [\binf,\bsup]$. 
The $\beta$-value associated to an edge $e \in E(T)$ will be called the bias of the edge. It will be denoted by $\beta_e$. \hfff{betae}
\end{definition}
\begin{definition}
We write $\P_{T,\beta}$ \hfff{ptblaw} for the law of 
the $\beta$-biased random walk $\big\{ X(i): i \in \N \big\}$ on the set of vertices of a weighted tree $T$. 
This process is the Markov chain on $V(T)$ with the following transition rules.
If a vertex $v \in V(T)$, $v
\not= \phi$, has offspring $v_1,\ldots,v_k$, then, for each $n \in \N$,  
\begin{eqnarray}
\P_{T,\beta} \Big( X(n + 1) = \overleftarrow{v} \Big\vert \xsub{n} = v \Big) & = &
\frac{1}{1 + \sum_{i=1}^k \beta_{v,v_i}},
        \label{pxyo} \\
\P_{T,\beta} \Big( \xsub{n+1} = v_j \Big\vert \xsub{n} = v \Big)    
& = & \frac{\beta_{v,v_j}}{1 + \sum_{i=1}^k \beta_{v,v_i}},
 \qquad \textrm{for $1 \leq j \leq k$}, \nonumber
\end{eqnarray}
where $\overleftarrow{v}$ denotes the parent of $v$, i.e., the neighbour of $v$ that is the closest to the root. The jump-law from the root is given by
$$
\P_{T,\beta} \Big( \xsub{n+1} = \phi_j \Big\vert \xsub{n} = \phi \Big)    
 =  \frac{\beta_{\phi,\phi_j}}{\sum_{i=1}^k \beta_{\phi,\phi_i}},
 \qquad \textrm{for $1 \leq j \leq k(\phi)$}.
$$  
For $v \in V(T)$, we write $\P_{T,\beta}^v$ for the law of $\P_{T,\beta}$ given that
$X(0) = v$. 
We write $\E_{T,\beta}^v$ for the expectation value under $\P_{T,\beta}^v$.
\end{definition}
The environment of our walk will have the geometry of a Galton-Watson tree.
\begin{definition}\label{deftreelaw}
Let $p = \big\{ p_i: i \in \N \big\}$, $\sum_{i=0}^\infty p_i = 1$, denote an offspring distribution.
Let $f:[0,1] \to [0,1]$, $f(z) = \sum_{i=0}^\infty p_i z^i$, denote the associated moment generating function.
We write $\P_f$ for the Galton-Watson tree with offspring distribution $p$, that is, for the law on rooted trees in which each vertex has an independent and $p$-distributed number of offspring.
We further write $\P_{f,\infty}$ for the law $\P_f$ conditional on the sample being infinite. 
\end{definition}
We will make the following very weak assumption on the offspring distribution.
\begin{hyp}\label{hypf}
The sequence  $\big\{ p_i: i \in \N \big\}$ satisfies $\sum_{k=0}^{\infty}
k p_k > 1$ and $p_0 > 0$. There exists $\alpha > 3$ and $\hatc > 0$ such that 
$\sum_{i \geq k}p_i \leq \hatc k^{-1} (\log k)^{-\alpha}$. 
\end{hyp}
We now describe the law of the random biases.
\begin{definition}
Let $\nu$ denote a probability measure whose support is contained in $\big[ \binf , \bsup \big]$.
We write $\P_{f,\nu}$ \hfff{gwtree} for the following law on weighted trees. A rooted tree is drawn from the distribution $\P_f$. Its weights $\beta_e:E(T) \to \big[ \binf , \bsup  \big]$ are independently chosen according to $\nu$. 
We write $\P_{f,\nu, \infty}$ for the weighted tree that results by so labelling a tree sampled according to the law $\P_{f,\infty}$. 
\end{definition}
The following non-lattice assumption on $\nu$ will play a fundamental role.
\begin{hyp}\label{hypnu}
The measure $\nu$ has  support in $\big[ \binf , \bsup \big]$, where 
$\bsup > \binf   > 1$ are fixed in the definition of a weighted tree from Definition \ref{deftree}.
 The support of $\nu \circ \log^{-1}$ is non-lattice. That is, the $\Z$-linear span of
$\log {\rm supp} (\nu)$ is dense in $\mathbb{R}$.
\end{hyp}
We require notation to describe the joint randomness of tree environment and walk:
\begin{definition}
Let $\P$ denote a law on weighted trees.
We write
 $\P \times \P_{T,\beta}^\phi$ for the joint distribution of a weighted tree sampled according to  
$\P$ and a $\beta$-biased random walk on that tree begun at its root. The weighted tree drawn according to a sample of such a measure will be called the environment of that sample. The law $\P$ on weighted trees given by the environment under $\P \times \P_{T,\beta}^\phi$ will be called the environment marginal of  $\P \times \P_{T,\beta}^\phi$.   
\end{definition}
The aim of this article is to study the long-term behaviour of the walk under the measure 
$\pf$. \hfff{pfnptb}
Randomly biased random walks on infinite trees have previously been studied by \cite{LyonsPemantle}, who provide a criterion for recurrence, or transience, valid for a very broad class of trees. (In our case, the assumption of $\binf > 1$ in Hypothesis \ref{hypnu} trivially implies transience.) In the case of supercritical Galton-Watson trees, randomly biased walk is studied in the recurrent regime by \cite{FHS}, who investigate the high-$n$ asymptotic of the maximal displacement of the walk during $[0,n]$. In \cite{aidekon}, the same model is analysed in the case of zero extinction probability and in the transient regime, the author presenting criteria for zero speed and for positive speed, and providing a formula for the exponent for sublinear growth in walker displacement in the zero-speed case. 

The model $\pf$, which is the object of our attention, has three levels of randomness: tree geometry, edge-biases, and walk. A more basic model arises by setting $\nu = \delta_\beta$, with some $\beta > 1$, so that the second type of randomness disappears, each edge having constant bias $\beta$. Indeed, this model is a natural candidate for an investigation of sub-ballistic walk in random environment, and is the subject of \cite{BAFGH}. It is useful to explain heuristically its behaviour.
Consider for a moment a toy model. Let $\mu$ be the law on
$(0,\infty)$, $\mu = \frac{1 - \alpha}{\alpha} \sum_{n=1}^\infty \alpha^n \delta_{\beta^n}$, where
$\alpha \in (0,1)$ and $\beta > 1$. Attach to each $n \in \N$ an
independent sample $\omega_n$ of $\mu$. A continuous-time walk
$X:[0,\infty) \to \N$, $X(0) = 0$, jumps successively to the right. It
waits at site $n \in \N$ for an exponentially distributed period whose mean
is $\omega_n$. If $\beta$ is sufficiently large, sub-ballistic behaviour
results. However, the time to reach $n$, after normalization, will not
resemble a stable law. This is because the holding times of the walk will
tend to remain within a bounded multiple of powers of $\beta$, so that, if
$\beta$ is high, the sequence of logarithms of these holding times, modulo
$\Z \log \beta$, will contain far more values close to zero than close to
$\frac{\log \beta}{2}$, say. This is inconsistent with a stable limit.
The constant bias walk on  a supercritical Galton-Watson tree resembles the
toy model, at least crudely, because the waiting periods model the delays caused by finite dead-ends for the walker on the tree (and indeed \cite{BAFGH} shows how stable limits do not arise in the tree model).
The distribution $\mu$ has a support whose logarithm lies in a lattice. 

In asking which tree-based random walk models may have stable limits, we are led, then, to the model with edge-biases randomized according to a law satisfying Hypothesis \ref{hypnu}.
This non-lattice condition is similar to that in~\cite{KKS}, which there ensures a stable limit law for the particle trajectory of a random walker in a one-dimensional random environment.  See \cite{ESZ1} and \cite{ESZ2} for explicit representations of constants arising in Kesten's derivation of a stable limit law in this case.

We need a final piece of notation, before stating the principal conclusion of this article.
\begin{definition}\label{hittime}
For $T$ a rooted tree, we set $\vert \cdot \vert: V(T) \to \N$, $\vert v \vert = d \big( \phi,
v\big)$.
For $T$ an infinite weighted tree, we set $\hitting{n} = \inf \big\{ i \in
\N: \vert X_i \vert = n \big\}$ \hfff{hitdist} to be the time that biased random
walk $\P_{T,\beta}^{\phi}$ on $T$ takes to move to distance $n$ from the root.
\end{definition}
\begin{theorem}\label{thm}
Let the offspring distribution $\big\{ p_i: i \in \N \big\}$ 
satisfy Hypothesis \ref{hypf}, and let the measure $\nu$ satisfy Hypothesis \ref{hypnu}.
Define $\chiexp > 0$ according to 
\begin{equation}\label{defchi}
\int_0^\infty y^\chiexp \nu(dy) = \frac{1}{f'(\pext)},
\end{equation}
where $\pext = \P_f\big( \vert V(T) \vert < \infty \big)$ \hfff{extprob} denotes the extinction probability.
Provided that $\chiexp < 1$, under the law $\pf$, the hitting time $\hitting{n}$ satisfies
\begin{equation}\label{statone}
 n^{-1/\chiexp} \hitting{n}   \rightarrow \xi^{1/\chiexp} S_{\chiexp},
\end{equation}
where the limit is in distribution as $n \to \infty$, where $S_{\alpha}$,
$$
\E \exp \big\{ - \lambda S_{\alpha} \big\} = \exp \big\{ - \lambda^{\alpha}
\big\}, \qquad \textrm{for $\lambda > 0$,}
$$
denotes the stable law of index $\alpha$. 
We have further that, under $\pf$,
\begin{equation}\label{stattwo}
 \frac{\big\vert X(n) \big\vert}{n^{\chiexp}} \rightarrow \xi^{-1} \big( S_{\chiexp}
 \big)^{-\chiexp},
\end{equation}
the limit in distribution as $n \to \infty$.
The constant $\xi$ is given by
\begin{equation}\label{eqnxi}
\xi = \constt \constr \cpo \cpt 
\Gamma(1 + \chiexp)
\Gamma(1 - \chiexp), 
\end{equation}
where $\Gamma(z) = \int_0^\infty t^{z-1} \exp \{ -t \} dt$ for $z > 0$. 

Explicit expressions for the constants $\cpo$ and $\cpt$ will be given in
Proposition \ref{theoremone} and Proposition \ref{trprop}. The 
value of $\constr$ is identified after Lemma \ref{lemrexist}, and $\constt$ is defined in 
Lemma \ref{lemf}.

Finally, if $\chiexp > 1$, then motion is ballistic: there exists $v \in (0,\infty)$ such that 
$n^{-1} \hitting{n} \to v$, $\pf$-a.s.
\end{theorem}
We show then that biased random walk with logarithmically non-lattice
biases on a leafy supercritical Galton-Watson tree gives rise in its sub-ballistic regime to a scaled
walk having a stable limiting law. The theorem also exhibits $\gamma = 1$ as the critical point in a phase transition from the ballistic to the sub-ballistic regime. The analogous statement in the case of constant bias appears as Theorem 4.1 in \cite{LPP96} (which also establishes that the speed is zero at the critical value). The final statement of Theorem \ref{thm} will be a simple inference from the techniques that we will develop to prove its assertions regarding the sub-ballistic case.

A fundamental step in our approach is the next result, which describes the law of the total time spent in a trap.
Note that we require no hypothesis on the exponent $\chiexp \in (0,\infty)$, although our interest in the result lies principally in the case that $\chiexp < 1$. Here and throughout, by $A(x) \sim B(x)$ is meant
$\frac{A(x)}{B(x)} \to 1$ as $x \to \infty$.
\begin{theorem}\label{proptraplaw}
Assume Hypotheses \ref{hypf} and \ref{hypnu}. Then, as $x \to \infty$,
 $$
  \btmeasbracbig \Big( \taustar > x \Big)  \sim  \constr \cpo \cpt \Gamma\big( 1 + \chiexp \big)  x^{-\chiexp}.
 $$

\end{theorem}  
We have yet to explain the meaning of the notations in Theorem \ref{proptraplaw}: the measure $\btmeas$ 
and the random variable $\taustar$, as well as the first three constants on the right-hand-side. For now, we invite the reader to keep this picture in mind: that the walk under $\pf$ follows a trajectory under which it moves at linear speed along a backbone in the tree, and that, from time to time, it moves away from this backbone, arriving at a new finite tree, which may then act as a trap, into which the walk may fall and so be delayed. At such a moment of arrival, we regard the walk as being at a new trap entrance. The measure $\btmeas$ has the interpretation of the limiting environment, viewed from the particle, on the occasion of its arrival at a new trap entrance, after many trap entrances have already been encountered. As such, configurations under the law $\btmeas$ have both an ``outside'', that describes the locale of the walker on the backbone, and an ``inside'', that specifies the finite trap that is the descendent tree of the trap entrance. The quantity $\taustar$ denotes the total time that the walk will spend in the trap, after its moment of arrival at the trap entrance. 
As such, Theorem \ref{proptraplaw} describes the law of the total time spent in a trap that is encountered at a very late time by the walk under $\pf$.

At least conceptually, Theorem \ref{thm} is a fairly direct consequence of Theorem~\ref{proptraplaw}.

In the next few subsections, we will define the objects appearing in the statement of Theorem \ref{proptraplaw}. As we do so, we will outline several of the principal ideas involved in the proof of this result, doing so in enough detail that we present the statements of the propositions that 
form the main stepping stones on the route to Theorem \ref{proptraplaw}. 
This first section concludes with accounts of two elements of the broader topic:
how the route to Theorem \ref{proptraplaw} would be much simpler for the case that the bias is constant,
and 
how the two aspects of the clustering/smoothing dichotomy identified in \cite{BAFGH}, and in \cite{paperone} 
and the present article, may both occur in the more physical model of anisotropic walk on a supercritical percolation cluster.

Section \ref{secroutethm} presents the proofs of the propositions in the earlier explanatory sections. It finishes with the proof of Theorem \ref{proptraplaw}.
In Section \ref{seciht},  we present the coupling construction by which Theorem \ref{thm} is derived
from Theorem \ref{proptraplaw}. This section begins with a counterpart to the introductory subsections to which we now turn, in which the strategy of this derivation is presented.\\

\noindent{\bf Acknowledgments.} I thank G{\'e}rard Ben Arous for his central role in initiating this research and for many valuable discussions throughout the period that this project was conducted. I extend thanks to Alex Fribergh and Nina Gantert for interesting conversations.
I thank a referee, for perceptive comments and suggestions to improve the presentation of several proofs, and a second, for detecting an inaccuracy in an earlier version.
\begin{subsection}{The backbone and the traps: the Harris decomposition}\label{secharris}
We now describe the decomposition of the tree
into a backbone, on which the walk moves at linear speed, and traps, by
which it may be waylaid. The assertions made here are proved in \cite{LyonsPeres}, Section 5.7.

Writing $f:[0,1] \to [0,1]$, $f(z) = \sum_{k=0}^{\infty} p_k z^k$, we
set $\pext = \P_f \big( \vert V(T) \vert  < \infty  \big)$ for the extinction
probability, so that $\pext \in (0,1)$ solves $f(\pext) = \pext$.
We introduce $h:[0,1] \to [0,1]$, $h(z) = \sum_{k=0}^{\infty} h_k z^k$, and 
$g:[0,1] \to [0,1]$ according to
\begin{equation}\label{hgform}
 h(t) = \frac{f \big( \pext t \big)}{\pext} \qquad \textrm{and}
\qquad g(t) = \frac{f\big( \pext + (1-\pext)t \big) - \pext}{1-\pext}.
\end{equation}
\begin{figure}
\begin{center}
\includegraphics[width=0.75\textwidth]{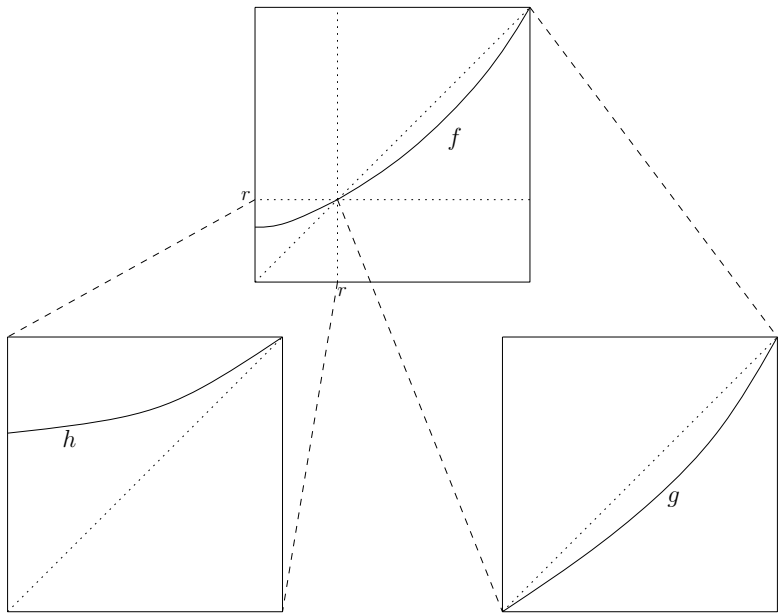} \\
\end{center}
\caption{The moment generating functions associated to the Harris decomposition of a Galton-Watson tree.}\label{figharrisdecomp}
\end{figure}

Then the law $\P_f$ may be sampled as follows. Firstly, an infinite
leafless tree is sampled according to $\P_g$. To each of its vertices is
added a random number of additional edges, whose law is determined by the
number of children of the vertex already constructed. Then, to each of
the endpoints of these new edges, an independent tree sampled 
according to $\P_h$ is attached. Introducing some terminology,
\begin{definition}
Let $T$ be an infinite rooted tree. Throughout, we write $\phi = \phi(T)$
for the root. 
Define the backbone $B = B(T)$ \hfff{threebb}
to be the subtree induced by the set of vertices lying in some infinite
vertex self-avoiding path beginning at $\phi$.
We define the augmented backbone $B_{{\rm aug}}$
to be the subtree induced by 
$\big\{ v \in V(T): d(v,B) \leq 1 \big\}$,
and write $\bext = V\big( B_{\rm aug }\big) \setminus V(B)$
for the exterior vertex boundary of the backbone.
\end{definition}
In these terms, the decomposition states that the backbone $B(T)$ of a sample $T$ of $\P_f$ has law $\P_g$, while the
descendent trees of elements in the exterior vertex boundary $\bext$
are independent, each having law $\P_h$.

For a sample of $\P_{f,\nu,\infty}$, we regard each element of $\bext$ as a trap entrance, with the descendent tree of that element being the corresponding trap. Each such trap is a weighted tree, with distribution $\P_{h,\nu}$. 

Throughout the paper, we use $\big\{ h_i: i \in \N \big\}$ to denote the offspring law arising from the given
 $\big\{ p_i: i \in \N \big\}$ by means of (\ref{hgform}). It is easy to see that, for any supercritical choice of 
 $\big\{ p_i: i \in \N \big\}$, the  $\big\{ h_i: i \in \N \big\}$  are exponentially decaying in $i$. 
\end{subsection}
\begin{subsection}{The return time to the trap entrance}
Under $\pf$, the duration $\big\{ 1,\ldots,\hitting{n} \big\}$ is comprised of moments at which the walk lies in the backbone, and the moments at which it lies in one of the traps. In a regime of sub-ballistic motion, we anticipate that most time is spent among the traps, and that, indeed, it is a small number of deep traps that account for most of the delay in reaching distance $n$ from the root. 

In an effort to describe the sojourn time in a trap, suppose that the walk has just arrived at a new trap entrance, that is, at a previously unvisited element of $\bext$. As a first step to understanding the delay caused by the trap, we present in this subsection results describing the law of the time taken for the walk, from its arrival at the trap entrance,  to return to its parent on the backbone, in the case of a typical trap. We will call this random variable the {\it escape time}.  

This toy problem concerns a walk on a finite tree which we introduce some notation to describe.
\begin{definition}
By a simple path in a graph $G$, we mean a finite list of its vertices, each adjacent to the last, and without repetitions. 
Given a tree $T$ and $v,w \in V(T)$, we write $P_{v,w} = P_{v,w}^T$ \hfff{pvw} for the unique simple path starting at $v$ and ending at $w$. We will also write $P_{v,w}$ for the graph induced from $G$ by the set of elements of this path.
\end{definition}
\begin{definition}
Let $T$ be a rooted tree. We say that $v \in V(T)$ is a descendent of $w \in V(T)$ if $w \in V\big(P_{\phi,v}\big)$. We define the descendent tree $T_w$ \hfff{desctree} of $w$ to be the subtree of $T$ induced by the set of descendents of $w$. The tree $T_w$ is a rooted tree, with $\phi(T_w) = w$.
\end{definition}
\begin{definition}
 A {\it single-entry} tree is a finite rooted tree such that $\phi$ has a unique offspring. For such a tree, we will refer to $\phi$ by $\head$. \hfff{headent} The offspring of $\head$ will be called the entrance and denoted by $\entrance$. 
\end{definition}
\begin{definition}
 Let $T$ be a weighted tree, and let $v,w \in V(T)$. 
Under the law $\P_{T,\beta}^v$, we define the hitting time $H_w = \inf \big\{ n \in \N: \xsub{n} = w \big\}$. \hfff{hitw}
\end{definition}
For a trap in an infinite weighted tree, consider the single-entry weighted tree $T$ induced by the trap entrance, its descendents, and its parent. The parent is identified with $\head$ and the trap entrance with $\entrance$. The escape time for this trap is given by $H_\head$ under $\P_{T,\beta}^\entrance$.

The paper \cite{paperone} is devoted to describing the significant properties of the law of the escape time when the trap is randomly chosen according to the law $\P_{h,\nu}$. We now state the results that will be relevant for our investigation.

\begin{subsubsection}{The mean return time}
The first crucial piece of data is a formula for the mean escape time $\E_{T,\beta}^\entrance(H_\head)$, where
$T$ is a single-entry weighted tree (associated to some trap). 
To present this, we require some definitions.
\begin{definition}\label{defomega}
Given a weighted tree $T$, we define the weight of a simple path 
$P$ in $T$ to be the product of the edge-biases $\beta_e$ over all $e \in E(P)$.
Whenever $x,y \in V(T)$, with $y$ a descendent of $x$, we write $\omega_x(y)$ \hfff{omegaxy}
for the weight of the path $P_{x,y}$. We set $\omega_x(x) = 1$. We also write, for any $x \in V(T)$,
 $\omega_x(T_x) = \sum_{v \in V(T_x)} \omega_x(v)$. We use $\omega(T)$ \hfff{omegat} to denote $\omega_\phi(T)$, the sum of the weights of all simple paths emanating from the root of $T$. 
\end{definition}
The commute-time formula for a reversible network was originally proved in \cite{CRRST}, and is presented as Theorem 3.3 in the overview \cite{Barlow}. The following formula is a special case. Let $T$ be a single-entry weighted tree. Then 
\begin{equation}\label{eqnmean}
\E_{T,\beta}^\entrance \big( H_\head \big) = 2 \omega_\entrance \big( T_\entrance \big) - 1.
\end{equation}
\end{subsubsection}
\begin{subsubsection}{Approximate exponential law for a deep excursion}
Of course, we would like to understand more about the law $\P_{T,\beta}^\entrance(H_\head)$ than the behaviour of its mean. We will say more, for trees $T$ that are typical in the relevant sense, that is, which are likely to occur under the measure $\P_{h,\nu}$ on traps arising in the problem. We now define this notion of a ``typical'' trap.  To do so, we will make use of the lexicographical ordering on the vertices of a finite rooted tree. In our application, we are concerned with Galton-Watson trees. We did not specify an ordering on vertices of such trees in the slightly imprecise Definition \ref{deftreelaw}. Using the formal coding of Galton-Watson trees given in Section 1.1 of \cite{LeGall}, such trees carry a lexicographical ordering. This material is also reviewed in Appendix A of \cite{paperone}.
\begin{definition}\label{defoutgr}
Let $T$ denote a weighted tree. Let $D(T)$, the depth of $T$, denote the maximal distance of $\phi_T$ to an element of $V(T)$. We define $\vbase$ \hfff{baset}
to be the lexicographically minimal vertex at distance $D(T)$ from $\phi$.
Write $P_{\phi,\vbase}$ in the form $\big[ \phi = \psi_0,\psi_1,\ldots,\psi_{D(T)} = \vbase \big]$,
and, for $0 \leq i \leq D(T)$, let $J_i$, the $i$-th outgrowth of $T$,  denote the connected component containing $\psi_i$ of the graph with vertex set $V(T)$ and edge-set 
$E(T) \setminus E \big( P_{\phi,\vbase} \big)$. 
\end{definition}
\begin{definition}\label{defbt}
For $\barecon > 0$, we say that a \wgt tree $T$ is $\barecon$-bare \hfff{baretree} if
$\big\vert V(J_i) \big\vert \leq \barecon \log \log \omega(T)$
for each $i \in \{ 0, \ldots , D(T) - 1 \}$.
\end{definition}
For the present paper, the details of the definition of a $\barecon$-bare tree are inconsequential: it formulates the notion of a tree being long and thin. The only property that we will use in this regard is the following lemma, which appears as Lemma 3 in \cite{paperone}, and which tells us 
that a high-weight trap is probably bare.
\begin{lemma}\label{lemcreg}
There exists a constant $\connew > 0$ such that, fixing
$\barecon > \connew^{-1}$, 
$$
\mathbb{P}_{h,\nu} \Big( T \, \textrm{is $\barecon$-bare}  \Big\vert \omega(T) \geq u  \Big) \geq 1 - (\log u)^{- \connew \barecon},
$$
for all sufficiently high $u$.
\end{lemma}
The next proposition, which is Proposition 1 (in Section 5) of \cite{paperone}, presents an approximate description of the law $H_\head$ under $\P_{T,\beta}^\entrance$, where $T$ is a single-entry weighted tree.
\begin{definition}\label{defde}
Let $T$ denote a single-entry weighted tree. 
Under the law $\P_{T,\beta}^\entrance$, the walk $X$ is said {\it make a deep excursion into} $T$ if $H_\vbase < H_\head$. We write $\falldeep$ \hfff{evfalldeep} for the event that $X$ makes a deep excursion in $T$.
We write $\pde = \pde(T) = \P_{T,\beta}^\entrance \big( \falldeep \big)$. 
\end{definition}
\begin{prop}
\label{proptspold}
Let $\barecon \in (0,\infty)$ be an arbitrary constant. 
There exist $\trwtbd \in (0,\infty)$ and $\concaptwo,\conec > 0$ such that 
the following holds. Let $T$ denote a $\barecon$-bare single-entry weighted tree such that 
$\omen (T_\entrance) > \trwtbd$. 

The distribution of $H_\head$ under $\P_{T,\beta}^{\entrance} \big( \cdot \big\vert \falldeep \big)$ is such that we may construct on this space  an exponential random variable $E$ with $\E_{T,\beta}^\entrance(E) = 2\omen \big( T_\entrance \big)/\pde$ and with 
$$
 \E_{T,\beta}^{\entrance} \Big( \big\vert H_\head - E \big\vert \Big\vert \falldeep \Big)  
 \leq \concaptwo \omen\big(T_{\entrance}\big)^{1/2};
$$
moreover, $E$  under $\P_{T,\beta}^{\entrance} \big( \cdot \big\vert \falldeep \big)$ may be chosen to be independent of the trajectory $X:\{0,\ldots,H_\vbase\} \to V(T)$.

We further have that
\begin{equation}\label{srineq}
\mathbb{E}_{T,\beta}^{\entrance} \Big( H_\head  \Big\vert \falldeep^c  \Big)
 \leq \concaptwo \Big( \log \omen\big(T_{\entrance}\big) \Big)^{\barecon \conec}.
\end{equation}
\end{prop} 
That is, the escape time satisfies the dichotomy: if there is no deep excursion, then the escape time is negligible; if there is, then the escape time is well approximated by an exponential distribution, of mean $2 \omen\big(T_{\entrance}\big)/\pde$.
\end{subsubsection}
\begin{subsubsection}{The tail of the weight of a trap}
Under $\pf$,  each newly encountered trap has law $\P_{h,\nu}$.
We now describe the law of the weight of the trap, or, equivalently by (\ref{eqnmean}), of the mean escape time to the backbone from the trap entrance.  
\begin{hyp}\label{hyph}
Let $\big\{ h_i:  \in \N \big\}$, $\sum_{k=1}^{\infty} k h_k < 1$, 
be a
subcritical offspring distribution for which there exists $c > 0$ such
that $\sum_{l \geq k }h_l \leq \exp \big\{ - ck \big\}$ for each $k \in \N$.  
We write $m_h = \sum_{k=1}^{\infty} k h_k$ for the mean number of offspring.
\end{hyp}
We now state Theorem 2 of \cite{paperone}. We emphasise that Hypothesis \ref{hyph} is not an additional hypothesis for the results in this article; indeed, we have noted that it is necessarily satisfied  for our choice of  $\big\{ h_i:  \in \N \big\}$.
\begin{prop}\label{theoremone}
Let $\big\{ h_i:  \in \N \big\}$ satisfy Hypothesis (\ref{hyph}).
For $\nu$ a distribution on $(1,\infty)$ satisfying (\ref{hypnu}), 
we have that
$$
 \mathbb{P}_{h,\nu} \Big( \omega(T) > u \Big) \sim \cpo u^{-\chiexp}, 
$$ 
as $u \to \infty$, where $\chiexp > 0$ satisfies
\begin{equation}\label{formchi}
\int_0^{\infty} y^\chiexp    \nu (dy) = \frac{1}{m_h},
\end{equation}
and where 
\begin{equation}\label{formcpo}
 \cpo = 
\frac{1}{\chiexp m_h} \Big( \int_1^{\infty} y^\chiexp  \log(y)  \nu (dy) \Big)^{-1}
\lim_{k \to \infty} \E_{h,\nu} \Big( \omega(T)^{\chiexp} 1\!\!1_{D(T) = k}  \Big).
\end{equation}
\end{prop}
The exponent $\chiexp$ in Theorem \ref{thm} has now made its appearance. (This $\chiexp$ is indeed that of Theorem \ref{thm}, since $m_h = f'(\pext)$ by our choice of $h$ given in (\ref{hgform}).) Although the proof of Proposition \ref{theoremone}
appears elsewhere, we mention briefly how this exponent arises: the trees contributing to 
$\mathbb{P}_{h,\nu} \big( \omega(T) > u \big)$ are typically long and thin, so that $\omega(T)$ is approximated up to a constant factor by the $\omega$-value of a vertex $\vbase$ at maximal distance from its root. This distance having a geometric distribution, $\log \omega (\vbase)$ is a geometric sum of independent random variables, each having law $\nu \circ \log^{-1}$. The defective renewal theorem may then be applied to find the asymptotic behaviour of the tail of this random variable. It is at this moment that we invoke Hypothesis~\ref{hypnu} on the law $\nu$. 

This approximation is good enough to explain the form of the exponent $\chiexp$, but not that of the constant $\cpo$. 
Indeed, the limit  in (\ref{formcpo}) may be considered as a constant correction that takes account of the numerous vertices in the environment near the base of the trap under $\P_{h,\nu}$ which contribute to $\omega(T)$, in comparison to the single vertex responsible for the term $\omega(\vbase)$.

We mention also that the scaling in Theorem \ref{thm} may be heuristically identified from Proposition \ref{theoremone}.
Indeed, since we expect the walk under $\pf$ to encounter an order of $n$ traps before time $\hitting{n}$, and for a macroscopic fraction of its time to be spent in the largest trap encountered, we may estimate $\hitting{n}$ by a value of $u$  for which $\mathbb{P}_{h,\nu} \big( \omega(T) > u \big) \approx n^{-1}$. Hence, Proposition \ref{theoremone} points to the scaling $\hitting{n} \approx n^{1/\chiexp}$. 
\end{subsubsection}
\end{subsection}
\begin{subsection}{Modelling a visit to a trap made at late time}
We now formally record the set $\big\{ \trar{i}: i \in \N \big\}$ of moments of {\it arrival} at trap entrances made by the walk.
\begin{definition}\label{deftxb}
Let $T$ be an infinite weighted tree.
Let $X$ have the law $\P_{T,\beta}^{\phi}$.
Let $\trar{1} = \inf \big\{ t \in \N : X(t) \in \bext \big\}$, \hfff{trapenti}
and
$$
\trar{i+1}  = \inf \Big\{ t > \trar{i}: X(t) \in \bext , X(t) \not\in  
     \big\{ X(\trar{1}),\ldots, X(\trar{i})   \big\}  \Big\}, \, \,\, i \geq 1.
$$
\end{definition}
That is, $\big\{ X(\trar{i}): i \in \N \big\}$ and $\big\{ T_{X(\trar{i})}: i \in \N \big\}$,  
are chronological lists of encountered trap entrances, and the corresponding traps.

We are working towards a precise understanding of how traps delay the walk under $\pf$. The results quoted so far tell us about the escape time from trap entrance to backbone. Of course, this is not enough for our purpose,
because the total time spent in a trap may be
interspersed with brief interludes on the backbone, and is, in effect, a
sum of several escape times, their number dependent on the environment
on the backbone neighbouring the trap entrance.
 To discuss this matter precisely, we now introduce a formal description that models a trap entrance encountered at late time. We are modelling a late time limit, at which the root of the global tree will have receded infinitely into the past. Despite the absence of a root, local notions of direction in the tree still make sense. As such, we introduce the following definition. 
\begin{definition}
A locally finite infinite tree $T$ is said to be end-distinguished \hfff{edtree} if its edges are oriented, with each vertex $v$
having exactly one ingoing incident edge $(w,v)$.  
We write $w = \overleftarrow{v}$.
We may extend the definition of biased random walk from the class of rooted
trees to \unidir ones:   
given a labelling $\beta_e: E(T) \to [\binf,\bsup]$ of a \unidir tree by biases, we may define the random
walk by means of the transition probabilties (\ref{pxyo}). 
In an extension of the usage of Definition \ref{deftree}, a \unidir tree carrying such a labelling will be referred to as a \wgt tree (even though such a tree lacks a root).  
\end{definition}
%
\begin{figure}
\begin{center}
\includegraphics[width=0.8\textwidth]{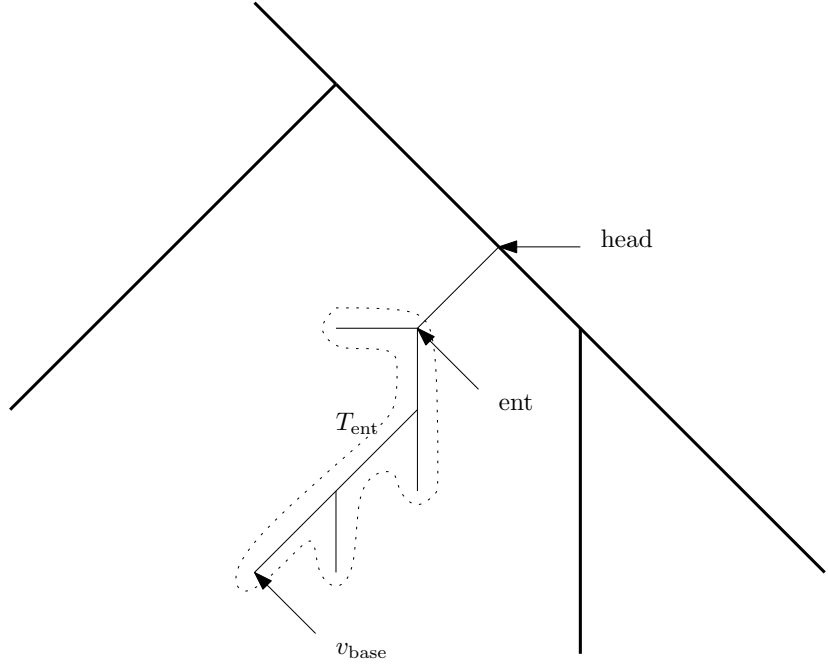} \\
\end{center}
\caption{A backbone-tree pair. The backbone $B$ is drawn with bold lines.}\label{figbackbonetree}
\end{figure}
The following definition is intended to model the backbone and trap, as
seen from the particle at late time, as it visits a trap entrance for the first time.
\begin{definition}\label{defbtp}
Let $B$ be an infinite rooted or \unidir \wgt tree without
leaves, with a
distinguished vertex $u \in V(B)$. 
Let $\singent$
be a single-entry weighted tree. 
We define the {\it backbone-tree pair} $(B,\singent)$ \hfff{btpair}
by contracting the tree
$$
\Big( V(B) \cup V(\singent), E(B) \cup E(\singent) \cup \big\{ \big( u,\head
\big) \big\} \Big)
$$ 
along the edge $\big( u, \head \big)$; that is, by identifying the
vertices
$u \in V(B)$ and $\head \in V(\singent)$. 
The direction of edges and the weights $\beta: E(B,\singent)  \to [\binf,\bsup]$ 
are naturally inherited from the constituent
graphs, as is the root $\phi(B)$, in the case that $B$ is rooted.
As such, $(B,\singent)$ is a \wgt tree, which is either rooted or \unidirnosp,
in accordance with $B$. 
We continue to write $\head$ for the vertex in $V(B,\singent)$ corresponding to both $u \in V(B)$ and
$\head(\singent) \in V(\singent)$. Similarly, the {\it entrance} $\entrance \in V(B,\singent)$ is 
the counterpart of $\entrance(\singent)$. We will occasionally use the notation $\overline{T} = (\head,\entrance) \circ T_\entrance$ to recall that $\overline{T}$ is composed by adjoining to the edge $(\head,\entrance)$ the tree $T_\entrance$.

In the case that $\singent$ is finite, we 
also record by $\btbase$ the lexicographically minimal vertex among those in $V(\singent)$ at maximal distance from $\entrance$.
\end{definition}

For a given backbone-tree pair, we are trying to understand the law of the total time spent in the trap:
\begin{definition}\label{deftauu}
Let $\big( B,\singent \big)$ 
be a backbone-tree pair.  Let  $u
\in V(B,\singent)$.
Under the law $\P_{(B,\singent),\beta}^u$, we define the random variable \hfff{tautent}
$$
\tau_{T_{\entrance}}  = \Big\vert \Big\{  i \in \mathbb{N}: X(i) \in V(T_{\entrance}) \Big\} \Big\vert
$$ 
for the total time spent by $X$ in the vertices of the tree $T_{\entrance}$.
\end{definition}
As we have mentioned, this total time may be spread between several separate visits to the trap. We now present the analogue of Proposition \ref{proptspold} for total visit time of the trap in a backbone-tree pair.
\begin{definition}\label{eqfdef}
Let $\big( B,\singent \big)$ be a backbone-tree pair, with $\singent$ finite. 
Let $\pesc$, \hfff{probesc} the {\it escape probability}, denote the $\P^\head_{(B,\singent),\beta}$-probability that $X$ never visits $\entrance$.
Let $\pde$ denote the $\P^\entrance_{(B,\singent),\beta}$-probability that $X$ 
visits $\btbase$ before its first visit to $\head$, (which, in the notation of Definition \ref{defde}, is the probability that the initial excursion of $X$ into $\singent$ is deep).
We further set $\halfcorrfac = \pesc^{-1} + \pde^{-1} - 1$. \hfff{correctionfactor}

Under  $\P^\entrance_{(B,\singent),\beta}$, we write $\deepfall$ \hfff{evdeepfall} for the event that  $X$ {\it falls deeply} into $T_\entrance$, namely that $X$ visits $\btbase$ at some positive time.
\end{definition}
\begin{prop}\label{proptsp}
Suppose that the edge-law $\nu$ has support in $[\binf,\bsup] \subseteq (1,\infty)$.
Let $\barecon > 0$ be a given constant. 
There exist $\trwtbd \in (0,\infty)$ and $\concaptwo > 0$ such that any backbone-tree pair
$\big( B,\singent \big)$, $\singent = (\head,\entrance) \circ T_{\entrance}$,
for which $\omega(T_{\entrance}) > \trwtbd$, and with $T_{\entrance}$ being a $\barecon$-bare tree,  
has the following property. 
The distribution of $\tau_{T_{\entrance}}$ under $\P_{(B,\singent),\beta}^{\entrance}$ is such that we may write
$$
\tau_{T_{\entrance}} = E 1\!\!1_\deepfall \, + \, \ero,
$$
where $\mathbb{P}_{(B,\singent),\beta}^{\entrance}(\deepfall) = \frac{\pde}{1 - (1-\pesc)(1-\pde)}$, and $E$ is
an exponential random variable of mean $\corrfac \omen (T_{\entrance})$ that is independent of 
$X:\{0,\ldots,H_\vbase\} \to V\big( B,\singent \big)$: that is, $E$ is independent of $\deepfall$; given $\deepfall$, $E$ is independent of the walk until its arrival at $\vbase$; and given $\deepfall^c$, it is entirely independent of the walk. The error $\ero$ satisfies
\begin{equation}\label{properrbd}
 \E_{(B,\singent),\beta}^{\entrance} \Big( \big\vert \ero \big\vert \Big)  
 \leq \concaptwo \omen(T_{\entrance})^{1/2}.
\end{equation}
\end{prop}

In other words, for a backbone-tree pair whose trap has a high weight and is typical, the total time spent in the trap is negligible, if the walk never reaches the base of the trap, and, in the other case, it is well approximated by an exponential distribution. Conditional on falling deeply into the trap, the mean time spent there $\corrfac \omega(T_{\entrance})$  is up to leading order.
The quantity $\halfcorrfac$ is a coefficient modifying the dominant term $2\omen(T_\entrance)$ in the mean time formula (\ref{eqnmean}) that accounts for return visits and the conditioning on falling deeply into $T_\entrance$. 
For this reason, we will call it the {\it correction factor}.
\end{subsection}
\begin{subsection}{Constructing the law for the environment of a late-time arrival at a trap entrance}
Continuing our investigation of the total time spent in among one of the deepest traps encountered before time $\hitting{n}$, we would like to find an analogue of Proposition \ref{theoremone}, for the case when total trap visit time is considered, rather than merely escape time for a single visit. Now, Proposition \ref{theoremone} discusses the law of the weight $\omega(T)$ of the trap, i.e., of the mean time spent before a return to the backbone. 
In searching for an analogue, it is natural to study the law of $\corrfac \omega(T_{\entrance})$, the mean time spent in the trap if the walk falls deeply into it, where we choose a suitable measure on backbone-tree pairs as the environment. Which measure should we choose? Recall that we are trying to study the trap, and its surroundings on the backbone, as witnessed on arrival at a new trap entrance at late time. For the trap $T_\entrance$, we simply use the law $\P_{h,\nu}$, while, for the backbone law, we should use a limiting distribution for the backbone viewed from the particle in a limit of arrival at a high-indexed trap entrance.

To state a result claiming the existence of this limiting law, we require some notation to describe the local backbone environment of a trap entrance.  
\begin{definition}\label{defpsik}
Let $T$ be any infinite \unidir \wgt tree. For
$u \in \bext(T)$, and $k \in \N$, we write $\bout_k(u)$ for the \wgt
tree induced by the vertex $u$ and the set of vertices in the backbone
$B(T)$ of $T$ at distance at most $k+1$ from $u$. We call $\bout_k(u)$ the $k$-large backbone neighbourhood of $u$. 
(Note that the natural choice of root $\phi\big(\bout_k(u)\big)$ for this weighted tree is 
the unique ancestor $v \in V(T)$ of $u$ for which $d(u,v) = k + 1$, since this choice induces the same parent-offspring relations as in the original tree.)
 We write $\Psi_k$ \hfff{nkoutpsi} for the
set of $\bout_k(u)$ arising for some such $T$, $u$ and $k$. For $\xi \in
\Psi_k$, we set $\entrance(\xi)$ to be the leaf corresponding to $u$ in the above
description, and $\head(\xi)$ for its unique neighbour. 
%
\end{definition}
Here, then, is our result asserting the existence of a limiting measure for the neighbourhood on the
backbone of the walk as it reaches a high-indexed trap entrance:
\begin{prop}\label{propqkq}
Let the offspring distribution $\big\{ p_i: i \in \N \big\}$ satisfy
(\ref{hypf}), and let $\nu$ be a measure supported on $(1,\infty)$.
Let $\qsubk(n)$ denote the law on $\Psi_k$ of $\bout_k\big( X(\trar{n}) \big)$ under $\pf$. 

For each $k \in \N$, there exists a measure $\qsubk$ on $\Psi_k$
such that 
$$
{\rm TV} \Big( \qsubk(n),\qsubk \Big)  \to 0,
$$
where ${\rm TV}\big( \cdot,\cdot \big)$ denotes the total variation metric on probability measures on $\Psi_k$. 

The set of measures $\big\{ \qsubk: k \in \N \big\}$ forming a consistent
family, we may define a measure $\queu$, supported on infinite \wgt
\unidir trees having a unique leaf (which leaf we label $\entrance$), such
that the $k+1$-neighbourhood of $\entrance$ under $\queu$ has the law $\qsubk$ for each $k \in \N$.
\end{prop}
Note that $\queu$ is supported on trees that are \unidirnosp, but that lack a root, since, under $\pf$, this root is lost in the receding past from the standpoint of the particle at late time.



We now define the law on backbone-tree pairs that models the trap encountered at late time:
\begin{definition}\label{defqstar}
Let $\qstarnoarg$ \hfff{ltabt}
denote the following law on backbone-tree pairs.
A sample of the law $\queu$ has a unique leaf $\entrance$.
Independently of this sample, we sample a weighted tree according to  $\P_{h,\nu}$, and identify its root with $\entrance$. 
We define $\pesc$, $\pde$  and $\halfcorrfac$
by means of Definition \ref{eqfdef}.
\end{definition}
We are ready to state an analogue of Proposition \ref{theoremone}:
\begin{prop}\label{trprop}
Let the offspring distribution $\big\{ p_i: i \in \N \big\}$ satisfy
Hypothesis \ref{hypf} and the law $\nu$, Hypothesis \ref{hypnu}.
We have that
$$
\qstarnoarg \Big( \corrfac \omen(T_{\entrance}) > u \Big) \sim \cpo \cpt u^{-\chiexp},
$$
as $u \to \infty$, 
where the constants $\cpo \in (0,\infty)$ and  $\chiexp \in (0,\infty)$ were specified in Proposition
\ref{theoremone}. The constant $\cpt$ is given by
$$
\cpt = \big( 1 - p_0/\pext \big)^{-1} \lim_{k \to \infty} \E_{\qstarnoarg}
  \bigg( \Big( \corrfac \omen \big( \basedef(C_k) \big) \Big)^{\chiexp}
  1\!\!1_{r(T_{\entrance}) > k} \bigg),
$$
where the quantities $\basedef(C_k)$ and  $r(T_{\entrance})$ will be introduced in Definition \ref{defndiv}.
\end{prop}
In words, we sample a trap, and its backbone environment, by stopping the walk on its arrival at a late trap entrance. The statistic  $\corrfac \omega(T_{\entrance})$ summarises the mean total time that the walk would spend in the trap, were it to fall deeply into the trap. Proposition \ref{trprop} provides the leading asymptotic for the decay of this random quantity.

The limit in $\cpt$ is an average of a weighting of the correction factor over trap environments. The correction factor is well approximated by data from the environment close to the trap entrance, and the averaging in the limit should be considered as an average over the environment viewed from $\entrance$. A key tool in obtaining this formula for $\cpt$ is a certain decomposition of the trap $T_\entrance$ into renewal blocks, that will be specified in the course of the proof of Proposition \ref{trprop}.
\end{subsection}
\begin{subsection}{Regeneration times for the walk}
We make a further comment regarding Proposition \ref{propqkq}. In the proof, we will use a coupling procedure that depends crucially on the use of regeneration times for the walk. We now define these.
\begin{definition}\label{defpsikt}
Let $T$ denote an infinite rooted \wgt tree.
The set of regeneration times ${\rm RG} \subseteq \N$ \hfff{regwalk} of the walk $X$ under $\P_{T,\beta}^\phi$
consists of those monents $i \in N$ at which $X(i) \in V(B)$
is visiting the backbone $B = B(T)$,  while
$d \big( \phi,
X(j) \big) < d \big( \phi, X(i) \big)$ whenever $X(j) \in V(B)$ and $j <
i$, and  $d \big( \phi,
X(j) \big) > d \big( \phi, X(i) \big)$ whenever $X(j) \in V(B)$ and $j >
i$.
\end{definition}
Moreover, when we come to derive Theorem \ref{thm} from Proposition \ref{theoremone}, we will apply Proposition \ref{propqkq}, in the following guise:
\begin{corollary}\label{corqkq}
 Under the hypotheses of Proposition \ref{propqkq}, let $\qsubk(n)$ now denote the law 
of $\bout_k\big( X(\trar{n}) \big)$ under $\big( \pf \big) \big( \cdot \big\vert 0 \in {\rm RG} \big)$. 
For each $k \in \N$, 
$$
{\rm TV} \Big( \qsubk(n),\qold_k \Big)  \to 0,
$$
as $n \to \infty$, where $\qsubk$ is as stated in Proposition \ref{propqkq}.
\end{corollary}
\end{subsection}  
\begin{subsection}{The first-order asymptotic for the total time spent in a trap encountered at late time}\label{subsectotaltime}
We now define the random variable $\taustar$ appearing in Theorem \ref{proptraplaw}, setting 
$$
 \taustar = \Big\vert \Big\{ j \in \N : X(j) \in V \big( T_\entrance \big) \Big\} \Big\vert
$$
under the law $\btmeas$, 
where $(B,\singent)$ and $X$
denote the backbone-tree pair, and the walk, 
under this law. 

Recall that we are seeking to understand the way in which the walk under $\pf$ spends the duration $\big\{ 1,\ldots,\hitting{n} \big\}$. Most of this time will be spent in a few big traps. 
The measure $\btmeas$
models the environment at the moment of arrival at a trap entrance at late time, and the subsequent behaviour of the walk in that environment. As such, to understand the law of time spent in the big traps before time $\hitting{n}$ under $\pf$,
 we study the behaviour of the tail of the law of $\taustar$.

Proposition \ref{proptsp} provides an approximate representation for total trap visit time,
and Lemma \ref{lemcreg} implies that the geometry of the environment under $\qstarnoarg$ is typically such that this approximation is a good one. As such, the following wil be a simple consequence of these two results:
\begin{lemma}\label{lemsmallerr}
We may construct under $\btmeas$
random variables $E$ and $\ero$ such that
\begin{equation}\label{reptaustar}
 \taustar = E 1\!\!1_{\deepfall} + \ero.
\end{equation}
Given the environment $(B,\singent)$ of a sample of $\btmeas$, $E$ has the conditional distribution of an exponential random variable whose mean is $\corrfac\omen(T_\entrance)$.
Recall that the event $\deepfall$ was specified in Definition \ref{eqfdef}.
The construction of $E$ may and will be carried out so that, given the environment $(B,\singent)$ under $\btmeas$,
$E$ is independent of the trajectory $X:\{0,\ldots,H_\vbase\} \to V(B,\singent)$ (both when $\deepfall = \{ H_\vbase < \infty \}$ occurs and when it does not); hence, $E$ and $\deepfall$ are independent.
The term $\ero$ is an error term whose tail is small, in the sense that there exists $\epsilon > 0$ such that, for all $x > 1$,
\begin{equation}\label{eqsmallerr}
\btmeasbracbig \Big(  \ero > x \Big) \leq \big( \log x \big)^{-\epsilon} 
\btmeasbracbig  \Big( E > x \Big).
\end{equation}
\end{lemma}
Now, the asymptotic decay of $E$ under $\btmeas$ follows directly from Proposition \ref{trprop}:
\begin{lemma}\label{lemvrg}
We have that
$$
 \btmeasbracbig \Big( E > u \Big) \sim  \cpo \cpt \Big( \int_0^{\infty} v^{\chiexp}  
  \exp \{ - v  \} dv \Big)   u^{-\chiexp},
$$
as $u \to \infty$,
where the constants $\cpo$ and $\cpt$ were introduced in Proposition \ref{theoremone} and Proposition \ref{trprop}. 
\end{lemma}

One element remains to derive Theorem \ref{proptraplaw}: it is easy to see that, if we can establish the existence of the limit
\begin{equation}\label{btmeasxlim}
\lim_{x \to \infty} \btmeasbracbig \Big(  \deepfall \Big\vert E = x \Big),
\end{equation}
then Lemmas \ref{lemsmallerr} and \ref{lemvrg} would yield Theorem \ref{proptraplaw}, with the constant $\constr$ being the value of this limit. (It may seem troublesome to condition on the zero-probability event $E = x$ in (\ref{btmeasxlim}). However, $E$ is built by an independent exponential randomization on top of an environment-determined mean, so, as we will see, the effect of this conditioning on the environment is much less singular.)

To prove the existence of this limit, we work as follows. We will argue firstly that $\deepfall$ is an event that is determined by the geometry, edge-bias data and walk local to the trap entrance, in the sense $\deepfall$ may be arbitrarily well approximated in $\pf$-probability by the event that the walk leaves a large neighbourhood of $\entrance$ at a point in the trap $T_\entrance$ (rather than at a point on the backbone). Secondly, we will show that this environment viewed from $\entrance$ has a limiting law as we condition $\btmeas$ by $E = x$ and take $x \to \infty$. We will call this limiting law $\rho$. Since 
$E$ and $\deepfall$ are conditionally independent given the environment $(B,\singent)$,
the upshot of these two statements is that the limiting value of $\btmeasbrac \big( \deepfall \big\vert E = x \big)$ is given by the probability under $\rho$ that a biased walk from its entrance falls indefinitely into $T_\entrance$.

We prove the existence of $\rho$ in two steps. Firstly, we show that there is a limiting distribution for the environment viewed from $\entrance$ for the backbone-tree pair law $\qstarnoarg \big( \cdot \big\vert \corrfac \omen (T_\entrance) > u \big)$ in the limit of high $u$. Secondly, we argue that the analogous limiting law under  $\btmeasbrac \big( \cdot \big\vert E = x \big)$ as $x \to \infty$ (which is $\rho$) exists and coincides with this distribution. We now state precisely the conclusions of these two steps.
\begin{definition}\label{deftemp}
Let $(B,\singent)$ denote a backbone-tree pair. For $k \in \N$,
write $(B,\singent)_k$ for the finite weighted tree
induced by the set of vertices in $V(B,\singent)$ at distance at most $k$ from $\entrance$.
(Similarly to the case of elements of $\Psi_k$ in Definition \ref{defpsik}, the root of  $\phi$ of  $(B,\singent)_k$ is equal to  the unique ancestor $v \in V(B,\singent)$ of $\entrance$ for which $d(\entrance,v) = k$.)

Let $\xi$ denote a law on backbone-tree pairs. For $k \in \N$, let $\xi_{[k]}$ 
 be the law on finite weighted trees which is equal to the distribution of $(B,\singent)_k$, where $(B,\singent)$ has law $\xi$.
\end{definition}
\begin{prop}\label{proprho}
Assume Hypotheses \ref{hypf} and \ref{hypnu}.
There exists a measure $\rho$, supported on backbone-tree pairs for which $T_\entrance$ is infinite, that is the limit as $u \to \infty$ of $\qstarnoarg \big( \cdot \big\vert \corrfac \omen(T_{\entrance}) > u \big)$ in the sense that, for each $k \in \N$,
$\qstarnoarg \big( \cdot \big\vert \corrfac \omen(T_\entrance) \geq u \big)_{[k]} \to \rho_{[k]}$ in total variation norm as $u \to \infty$. 
\end{prop}
(The square bracket notation in $\rho_{[k]}$ is adopted because, when we come to prove Proposition \ref{proprho}, we will make use of a different system of marginal distributions; and to avoid notational complication in the proof, we reserve the notation $\rho_k$ to denote that other use.)
\begin{lemma}\label{lemrhox}
The law $\rho$, as defined by Proposition \ref{proprho}, is the limit as $x \to \infty$ of the environment marginal of $\btmeasbrac \big( \cdot \big\vert E = x \big)$, in the sense of this proposition.
\end{lemma}
\begin{definition}
Let $(B,\singent)$ be a backbone-tree pair where $T_\entrance$ is finite. Recall that $\deepfall$ 
is defined for the walk $\P_{(B,\singent),\beta}^\entrance$ in Definition \ref{eqfdef}. We extend this definition to the case where $T_\entrance$ is infinite as follows.
For such a backbone-tree pair, under the law $\P_{(B,\singent),\beta}^\entrance$, we set 
$$
\deepfall = \Big\{ X(i) \in V\big(T_\entrance\big) \, \, \textrm{for all sufficiently high $i \in \N$} \Big\}.
$$
\end{definition}
\begin{lemma}\label{lemrexist}
The limit (\ref{btmeasxlim})
exists, and is equal to $\big( \rho \times \P_{(B,\singent),\beta}^\entrance \big) \big( \deepfall \big)$.
\end{lemma}
The value of $\constr$  appearing in Proposition \ref{proptraplaw} is equal to the value identified in Lemma \ref{lemrexist}. 
Alternatively, it is given by
\begin{equation}\label{defnr}
 \constr  = \mathbb{E}_{\rho} \Big[  \frac{\pde}{1 - (1-\pesc)(1-\pde)} \Big],
\end{equation}
where the quantites $\pde$  and $\pesc$ are specified in Definition \ref{eqfdef}.
\end{subsection}
\begin{subsection}{The environment viewed from the trap entrance at late time}\label{secenv}
The law $\rho$ is closely related to the environment viewed from the entrance of the trap in which the walker lies in the limit of late time. Here, we give a precise statement in this regard. Our assertion follows in essence from the method by which Theorem \ref{thm} is derived from Theorem \ref{proptraplaw} in Section \ref{seciht}. We sketch the argument at the end of Section \ref{seciht}, although we omit details. 

To obtain this limiting environment $\hat\rho$, the law $\rho$ has to be conditioned  
to take account of deep falling into the trap.
\begin{definition}
Let $\hat\rho$ denote the law on backbone-tree pairs that is the environment marginal of 
$\big( \rho \times \P_{(B,\singent),\beta}^\entrance \big) \big( \cdot \big\vert \deepfall \big)$. 
\end{definition}
\begin{definition}
Under the law $\pf$, define the ``walk on trap entrances" $Y:\N \to \bext$ by setting $Y(n)$ equal to the last element of $\bext$ visited at or before time $n$.
\end{definition}
In fact, $Y$ is defined only after the time of the first visit of $X$ to a trap entrance. This does not matter since we will be concerned only with late time. The walk $X$ spends asymptotically all of its time in traps; hence, $Y$ is in essence a record of the entrance of the trap in which $X$ currently lives.   
\begin{definition}
Under the law $\pf$, for each $n \in \N$, we let
$\mathcal{Y}(n)$ denote the backbone-tree pair given by viewing $T$ from $Y(n)$. That is, $\mathcal{Y}(n) = (B,\singent)$, where the descendent tree $T_\entrance$ in $\mathcal{Y}(n)$ is identified with $T_{Y(n)}$ under $\pf$, and where $\bout_k(\entrance)$ under $\mathcal{Y}(n)$ is identified with  $\bout_k(Y(n))$ under $\pf$ for each $k \in \N$.
\end{definition}
Then, 
assuming Hypotheses \ref{hypf} and \ref{hypnu}, 
$\hat\rho$ has the intepretation of the environment viewed about the trap entrance at late time, in the sense that $\mathcal{Y}(n)$ under $\pf$ converges to $\hat\rho$ as $n \to \infty$, in the sense of Proposition \ref{proprho}. See the final Section \ref{secenvproof} for a sketch of the proof of this assertion.
\end{subsection}
\begin{subsection}{The form of the constant $\xi$ in Theorem \ref{thm}}
We have completed the overview of the proof of Theorem \ref{proptraplaw}.
The various aspects of the proof each leave an imprint in the form of a factor appearing in the constant $\xi$
given in  (\ref{eqnxi}) in the statement of the main result, Theorem \ref{thm}.
The constants $\cpo$ and $\cpt$ are corrections that account respectively for how the trap environment near its base influences the walk's escape from the trap, and how returns of the walk to the trap entrance increase the total sojourn time in the trap.
 The constant $\constr$ is the probability of falling deeply into a random trap with high mean visit time that is encountered at late time.
The remaining constant $\constt$ will be introduced during the derivation of Theorem~\ref{thm} from Theorem~\ref{proptraplaw}: it is the mean number of trap entrances encountered by the walk per unit of advance from the root. 
\end{subsection}
\begin{subsection}{A comparison of the approach with that for the case of constant bias}
The principal aim of this paper is to derive the stable limiting law for walker motion in the form of Theorem \ref{thm}, 
drawing on the regularity of the tail of trap-weight in Proposition \ref{theoremone} that is proved in \cite{paperone}. 
As we have described, Proposition \ref{theoremone} is false for the case of constant bias. Nonetheless, the task undertaken by \cite{BAFGH} for the constant bias case is an analogue of that performed together by \cite{paperone} and the present work: to describe precisely the asymptotic trajectory of the walker in the environment in question. The paper \cite{BAFGH} achieves this understanding in a much shorter space than the combination of the present article and \cite{paperone}. As such, we would like to discuss how the case of random bias  is harder.

We will focus on those difficulties that are specific to the task of the present article, which is in essence to derive Theorem \ref{proptraplaw} from Proposition \ref{theoremone}. Of course, Theorem \ref{proptraplaw} is not true for the constant bias case. However, for the purposes of the comparison, we would like to imagine that it were true, in order to make a direct comparison of what the present paper does, and what it would have to do in the constant bias case.  

The crucial difference between the two models occurs around the trap entrance. Consider a backbone-tree pair 
$(B,\singent)$, where $\singent$ is finite but large. The probability $\pde$ of deep excursion is random in the case of random bias, being approximately expressible in terms of the biases of the edges near $\entrance$ on the path $P_{\entrance,\vbase}$. Consider instead the constant bias case, where the bias is equal to some given $\beta > 1$. For deep enough $\singent$, the quantity $\pde$ is approximable to arbitrary precision by the probability that a nearest-neighbour random walk on $\Z$ with rightward transition probability $\beta/(1 + \beta)$ never reaches $0$ in departing from $1$. Thus $\pde \to 1 - \beta^{-1}$ in the limit of deep $\singent$.

This asymptotic non-randomness permits a dramatically simpler derivation of the law of the total time spent in a trap: as we now explain, the path from Proposition \ref{theoremone} to Theorem \ref{proptraplaw}, which is the subject of Section \ref{secroutethm}, is much easier to travel in the constant bias case.
  
 That $\pde$ is effectively non-random for all large traps has the upshot that the correction factor $\halfcorrfac = \pesc^{-1} + \pde^{-1} - 1$ is a function of $\pesc$ alone, and thus is measurable with respect to the trap exterior. This means that Proposition \ref{trprop} becomes a trivial consequence of Proposition \ref{theoremone}: $\big( 2\halfcorrfac \big) \omen(T_\entrance)$ is a product form whose first term depends only on the trap exterior and whose second depends only on its interior. The first term being uniformly bounded, this factorization into independent terms yields Proposition \ref{trprop}
by means of a simple tool for heavy-tailed random variables (of which we will anyway have need: see Lemma \ref{lemlit}(i)). 
The constant $\olddone$ would be given by $\E_{\qstarnoarg} \big( (2f)^\chiexp \big)$ in this case. (Of course, this formula is valid only under the absurd conceit that Theorem 2 holds for the constant bias case.)

After Proposition \ref{trprop} has been established, the principal step to obtain Theorem \ref{proptraplaw} is the existence of the limit in (\ref{btmeasxlim}), the asymptotic deep-falling probability $\constr$, a step that requires the construction of the limiting measure $\rho$ about the trap entrance. 
The non-randomness of $\pde$ in the constant bias case obviates the need for the non-trivial argument in this step. We may proceed in a slightly different fashion. The adopted approach treats the question of whether the walk falls deeply into the trap as a final step, after the total visit time given deep falling has been understood in Lemma \ref{lemvrg}.
Instead, we may handle deep falling before analysing total visit time. Under the new approach, we consider a trap reached at late time, i.e., a backbone-tree pair sampled under the measure $\qstarnoarg$, and condition straightaway on the event $\deepfall$ of falling deeply into the trap. In the constant bias case, this conditioning has asymptotically no influence on the interior of large traps. It changes only the trap exterior, the backbone nearby the entrance. That is, the law on backbone-tree pairs obtained by this conditioning has the form of a certain measure on the backbone, with the conditional form of the trap being $\P_{h,\nu}$-distributed, just as it was under $\qstarnoarg$. 
We denote this law by $\qstarl{\rm FD}$. We then obtain Lemma \ref{lemvrg} for the measure $\qstarl{\rm FD} \times \P_{(B,\singent),\beta}^\entrance$ as in the proof in this article. To prove Theorem \ref{proptraplaw}, then, we compute 
$\btmeasbrac \big( \taustar > x \big)$ by conditioning firstly on $\deepfall$, and then using the new version of Lemma \ref{lemvrg} to estimate the resulting probability.

 This treatment of deep falling - before conditioning on large time spent in a trap, rather than after it - would raise significant complications for the case of random biases. Traps would no longer be $\P_{h,\nu}$-distributed. Rather, the edges leading into the trap from $\entrance$ would experience a significant reconditioning. This would have the effect that Proposition \ref{theoremone} on the tail of the trap-weight would not be applicable, since it would not accurately describe the measure on traps.

\end{subsection}
\begin{subsection}{Potentially analogous results and the value of the techniques for related lattice models}\label{secaniso}
The present article and its prequel \cite{paperone}, and \cite{BAFGH},
 present a dichotomy for models of random walks on supercritical Galton-Watson trees with leaves. 
On the one hand, there is no scaling limit in the sub-ballistic regime for random walk with a constant bias, 
but, rather, 
there is a persistent discrete inhomogeneity in this motion, in which sojourn times in traps tend to cluster around powers of $\beta$. In contrast, the present paper, by randomizing biases with a law $\nu$, and moving from the previous case of $\nu = \delta_\beta$ to one in which $\nu$ satisfies Hypothesis \ref{hypnu}, perturbs the model enough that this discrete inhomogeneity evaporates, and a stable limit law for the scaled walker may be derived.

It has been suggested in the physics literature~\cite{StaSto} that such a discrete inhomogeneity may occur for an anisotropic walker on a supercritical percolation cluster in $\Z^d$, for $d \geq 2$. 
The slope of the preferred direction of the anisotropic walker in a disordered Euclidean setting is a parameter which is additional to the magnitude, or bias, of the anisotropy. As Alex Fribergh has suggested, it is reasonable to suppose that ``irrational" choices of this slope may interrupt the discrete effects that appear to obtain when the slope is axial, so that a stable limit may arise for such choices of the slope. The irrational slope would replace non-lattice randomization of edge-biases as the mechanism that achieves a stable limit. In this sense, the two behaviours identified here and by \cite{BAFGH} may coexist side-by-side in more physical models of anisotropy in disordered media.

The article \cite{AlexAlan} identifies the exact part of the natural parameter-space under which anisotropic walk on supercritical percolation has zero speed, in dimensions $d \geq 2$. It also specifies the sub-ballistic exponent for walker displacement in this regime, that is to say, the analogue of the exponent $\chiexp$ in Theorem \ref{thm}. 

The perspective of Theorem \ref{proptraplaw}, of constructing and analysing the limiting law of the particle at its arrival at a new trap at late time, may be fruitful for such physical models. As discussed in Section 11 of \cite{AlexAlan}, the traps that are effective at delaying the walker in dimensions $d \geq 3$ have a much narrower cross-section than in the two-dimensional case. This suggests that the point of view of Theorem \ref{proptraplaw} may be of more relevance to dimensions $d \geq 3$, since the notions of trap head and its locale may be well defined here.   
\begin{figure}
\begin{center}
\includegraphics[width=0.5\textwidth]{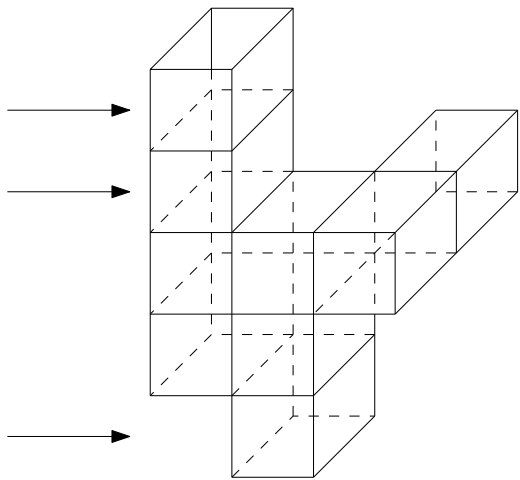} \\
\end{center}
\caption{The outer surface of a trap in $\Z^3$, composed of dual plaquettes. Arrows indicate levels where the cross-section is just a single unit square. These levels might be expected to populate to constant density a typical trap. The centre of the highest box may be regarded as the trap entrance, with the environment viewed from this vertex being analogous to our definition of a backbone-tree pair.}\label{trapinthreed}
\end{figure}
In the preceding section, we discussed how, in the constant bias case, Proposition \ref{trprop} would become a direct consequence of Proposition \ref{theoremone}, by means of a convenient factorization of the interior and exterior about the entrance of a backbone-tree pair. In the random bias case, we are forced to analyse the correction factor without recourse to this independence assumption. With regard to the potential of these techniques to work for more physical models, this technique of proof may be useful, since in such models there will be no independence in the interior and exterior of a trap viewed from its entrance. We mention also that the coupling construction by which we derive Theorem \ref{thm} from Theorem \ref{proptraplaw} in Section \ref{seciht} may provide a useful alternative to the arguments adopted in \cite{BAFGH}, since the notion of a super-regeneration time in Section $2$ of \cite{BAFGH} has no clear analogue for such models as anisotropic random walk in supercritical percolation in $\Z^d$. 

 The present article began in an effort to establish a stable limiting law by perturbing the constant bias walk on a supercritical tree. At this time, however, the potential of the techniques developed here to establish such behaviour in more physical models of anisotropy is at least as powerful a motivation. 
\end{subsection}

\end{section}
\begin{section}{The route to Theorem \ref{proptraplaw}}\label{secroutethm}
In Section \ref{secproptsp}, the brief argument that Proposition \ref{proptspold} implies Proposition \ref{proptsp} is given. Sections \ref{sectrprop}, \ref{secproofs} and \ref{secrho} respectively give the proofs of: Proposition \ref{trprop};
Proposition \ref{propqkq} and Corollary \ref{corqkq}; and Proposition \ref{proprho}. The proof of Theorem \ref{proptraplaw} is completed in Section \ref{secthmfin}, where the proofs of results from Section \ref{subsectotaltime} are given. Note that the values of the large and small positive constants $C$ and 
$c$ may change from line to line.
\begin{subsection}{Proof of Proposition \ref{proptsp}}\label{secproptsp}
Under the law $\P_{(B,\singent),\beta}^{\entrance}$, set $S_1= 0$ and $T_1 = H_\head$. Iteratively, for $j > 1$, set $S_j = \inf \big\{ i > T_{j-1}: X_i = \entrance \big\}$ and $T_j =  \inf \big\{ i > S_j: X_i = \head \big\}$, using the convention that $\inf \emptyset = \infty$. Set $G = \sup\big\{ j \in \N: S_j < \infty \big\}$. For $1 \leq j \leq G$, the time of the $j$-th sojourn in $T_\entrance$, $\hat\tau_j$, is defined by $\hat\tau_j = T_j - S_j$. 

Note that, under  $\P_{(B,\singent),\beta}^{\entrance}$,
 $\tau_{T_{\entrance}} = \sum_{i=1}^G \hat\tau_j$. Note further that $G$ has a geometric distribution, satisfying
 $\P_{(B,\singent),\beta}^{\entrance} \big( G = i \big) = (1-\pesc)^{i-1}\pesc$ for $i \geq 1$. 
 
Recall that the single-entry tree $\singent = (\head,\entrance) \circ T_\entrance$ is formed by prefixing the edge  $(\head,\entrance)$ to $T_\entrance$, so that $\phi(\singent) = \head$.  Note that, for each $i \in \N$, the conditional distribution of 
$\hat\tau_i$ under  $\P_{(B,\singent),\beta}^{\entrance}$, given that $i \leq G$, coincides with $H_{\head}$ under
$\P_{\singent,\beta}^\entrance$. As such, we may apply Proposition \ref{proptspold} to approximate the distributions of each $\hat\tau_i$. To do so, for $i \leq G$, write $\falldeep_i$ for the event that $X_j = \vbase(\singent)$ for some $S_i \leq j \leq T_i$. Note then that the event $\deepfall$ of falling deeply into $T_\entrance$ may be written
$\deepfall =  \cup_{i=1}^G \falldeep_i$.

Applying the construction in the statement of Proposition \ref{proptspold} to the independent sojourns, we obtain
$$
\sum_{i=1}^G \hat\tau_i =  \sum_{i=1}^G  E_i 1\!\!1_{\falldeep_i} \, + \, \mathcal{E},
$$
where the $E_i$ are a sequence of independent exponential random variables of mean $2 \pde^{-1} \omen \big( T_\entrance \big)$ constructed under $\P_{\singent,\beta}^\entrance$. The quantity $\mathcal{E}$ is an error term, whose mean is bounded by Proposition \ref{proptspold}: we have that
\begin{equation}\label{eqerrbd}
 \E \big( \mathcal{E} \big) \leq 
 \concaptwo \E \big( G \big)   \omen \big( T_\entrance \big)^{1/2} \leq \concaptwo \frac{\bsup + \binf + 1}{\binf -1}  \omen \big( T_\entrance \big)^{1/2}.
\end{equation}
where we begin to abbreviate $\P = \P_{\singent,\beta}^\entrance$, and to write $\E$ for the corresponding expectation. The latter inequality relied on $\E(G) = \pesc^{-1}$ and a uniform lower bound on $\pesc$
that we will shortly record in Lemma \ref{lemebfbd}. Note that (\ref{eqerrbd}) gives the error bound (\ref{properrbd}) in the statement of the proposition.

Set $E' = \sum_{i=1}^G  E_i 1\!\!1_{\falldeep_i}$. It is easy to verify that, conditionally on $\deepfall$, 
$E'$ has the memoryless property and thus has an exponential distribution. The events $\big\{ \falldeep_i: 1 \leq i \leq G \big\}$ being conditionally independent given $G$, we have that
\begin{equation}\label{eqprobfd}
\P \big( \deepfall \big) = \sum_{i=1}^\infty \pesc(1-\pesc)^{i-1} \Big( 1 - (1-\pde)^i \Big)
 = \frac{\pde}{1 - (1-\pesc)(1-\pde)}. 
\end{equation}
Noting that $E' = 0$ on $\deepfall^c$, we have that
$\E \big( E' \big\vert \deepfall \big) = \frac{\E(E')}{\P\big( \deepfall \big)}$. Note that
$\E(E') = \E(G) \E(E_1) \P(\falldeep_1) = \pesc^{-1} \cdot 2 \pde^{-1} \omen \big( T_\entrance \big) \cdot \pde = 2 \pesc^{-1}  \omen \big( T_\entrance \big)$. From (\ref{eqprobfd}), we obtain
$\E \big( E' \big\vert \deepfall \big) = 2 \big( \pde^{-1} + \pesc^{-1} - 1\big) \omen \big( T_\entrance \big)$.
Note further that $E'$ is independent of the walk until time $H_{\vbase}$, because each $E_i$ constructed by the use of Proposition \ref{proptspold} has this property. We then define $E = E'$ on $\deepfall$; we choose $E$ on $\deepfall^c$ independently of other details of the walk so that $E$ has the same conditional distributions on $\deepfall$ and on $\deepfall^c$. 
This completes the proof. \qed 
\end{subsection}
\begin{subsection}{The proof of Proposition \ref{trprop}}\label{sectrprop}
The considerations required to prove 
 Proposition \ref{trprop} given Proposition \ref{proptsp}
arise only from the presence of
the factor of $\corrfac$ in the conditioning of $\qstarnoarg$ on $\corrfac \omen(T_\entrance) > u$. Indeed, were this factor to be omitted,
Proposition \ref{trprop} would reduce to Proposition \ref{proptsp},
 since $T_\entrance$ under $\qstarnoarg$
has the law $\P_{h,\nu}$. Our approach to handling the presence of the term $\corrfac$ is to approximate it by a quantity that is defined by the geometry and edge-bias data within a bounded region of $\entrance$, and then condition on $(B,\singent)$ in this region. We do this in such a way that the part of $T_\entrance$ ``below'' the conditioned neighbourhood has the law of $\P_{h,\nu}$. This will enable us to apply Proposition \ref{proptsp}.

To make this plan work, we must split the tree $T_\entrance$ into pieces in a convenient way, which we now describe.
\begin{subsubsection}{The renewal decomposition of a tree}\label{secren}
We recall from   Appendix C of \cite{paperone} the {\it renewal decomposition} of a finite rooted tree. 
\begin{definition}\label{defndiv}
By a root-base tree $T$, we refer to a finite rooted tree, one of whose vertices $\basedef$ at maximal distance from $\phi$ is declared to be the base. \hfff{renewaldecomp}

Given a rooted tree $T$, a vertex $v \in V(T)$, $v \not= \phi$, is called a cutpoint if   
it is not a leaf, and any other vertex in $T$ at the same distance from $\phi$ as $v$ is a leaf.
The set of cutpoints naturally decompose a rooted tree into components in the following manner. We write $r(T)$ for the number of cutpoints of $T$ plus one.
We may then record these cutpoints in the form $c_i$, $1 \leq i \leq r(T) - 1$, in increasing order of distance from the root $\phi$. 
We further set $c_0 = \phi$. 
We set $d_i = d \big( \phi, c_i \big)$ for 
$0 \leq i \leq r(T) - 1$. We also set $d_{r(T)}= D(T)$, where recall that $D(T) = \max \big\{ d \big( \phi,v \big): v \in V(T) \big\}$.  
For $1 \leq i \leq r(T)$, we define the $i$-th component $C_i$ of the tree $T$ to be the subgraph of $T$ induced by the set of vertices in $T$ at a distance from the root of at least $d_{i-1}$ and at most $d_i$.
Then, for $1 \leq i \leq r(T) - 1$, $C_i$ may be regarded as a root-base tree, with
$\phi(C_i) = c_{i-1}$ and $\basedef(C_i) = c_i$.  The final component $C_{r(T)}$, however, is a tree rooted at $c_{r(T) - 1}$ that has no natural choice of base. 

For $1 \leq i \leq r(T)$, we set $S_i$ equal to the subgraph of $T$ induced by the set of vertices in $T$ at a distance from the root of at most $d_i$. Note that $V \big( S_i \big) = \cup_{j=1}^i V \big( C_j \big)$. If $i < r(T)$, then $S_i$ is a root-base tree, with $\phi(S_i) = \phi_T$ and $\basedef(S_i) = c_i$.
\end{definition}
\begin{figure}
\begin{center}
\includegraphics[width=0.4\textwidth]{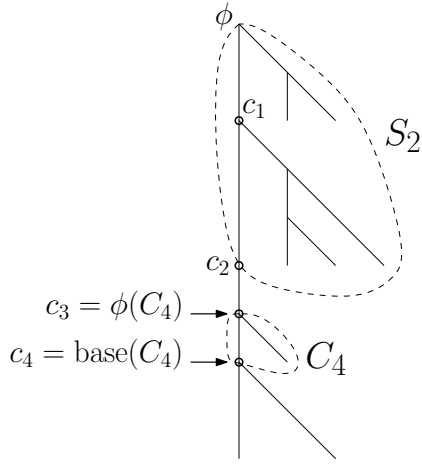} \\
\end{center}
\caption{The renewal decomposition of a tree $T$ satisfying $r(T)=5$.}\label{figrenewaldec}
\end{figure}
In a sense, the renewal decomposition is a counterpart of the ``renewal levels'' for trapping surfaces in $\Z^3$ that are indicated in the sketch in Section \ref{secaniso}.
The key property enjoyed by the decomposition that we will employ is stated in the following lemma, which explains why we may think of the cutpoints under the decomposition as ``regeneration'' points in the tree. The result is Lemma 28 of \cite{paperone}.
\begin{lemma}\label{lemdecom}
Let $k \in \N$.
Consider the law $\P_{h,\nu}$ conditionally on $r(T) \geq k+1$ and on an arbitrary form for the first $k$ weighted components of $T$. Then the conditional distribution of the descendent tree $T_{\basedef(C_k)}$
is given by $\P_{h,\nu}$ conditioned to contain at least one edge.
\end{lemma}

In the proof of Proposition \ref{proprho}, Lemma \ref{lemdecom} will allow us to condition the top $k$ components inside the trap $T_\entrance$
under a backbone-pair sampled from $\qstarnoarg$ so that hanging from the base of this conditioned region is a trap having the original trap law $\P_{h,\nu}$.

We also need to know that the components of the renewal decomposition of a high-weight tree are typically small. The following is 
Proposition 3 
of \cite{paperone}; although the result is plausible, it has a sizable proof.
\begin{lemma}\label{lemdivexp}
For the statement, we take $C_i = \emptyset$ if $i > r(T)$ (for any rooted
tree $T$). There exists $c > 0$ such that, for all $u > 0$ and $i \in \mathbb{N}$,
$$
\mathbb{P}_{h,\nu} \Big( \big\vert V \big( C_i \big) \big\vert \geq k \Big\vert \omega(T) > u \Big) \leq \exp \big\{ - c k  \big\},
$$
for each $k \in \mathbb{N}$.
\end{lemma}
\end{subsubsection}
\begin{subsubsection}{Approximating the correction factor using data in a neighbourhood of the trap entrance}\label{secequil}
We now define our approximation $\halfcorrfack$ to the correction factor $\halfcorrfac$ of Definition~\ref{eqfdef}.
\begin{lemma}\label{lembt}
Let a backbone-tree pair $(B,\singent)$  be given.
Extending Definition \ref{defndiv}, we introduce \hfff{bink}
\begin{equation}
 \ik  =  \left\{ \begin{array}{rl}
  S_k \big(T_\entrance \big) & \textrm{if } r(T_{\entrance}) \geq k,  \\
            T_{\entrance} & \textrm{otherwise.} \label{eqarr}
                  \end{array} \right. 
\end{equation}
for the rooted \wgt tree induced by the union of the first $k$ components
of the tree $T_{\entrance}$. (Note that $\ik$ is a root-base tree precisely
when $r(T_{\entrance}) \geq k+1$.) We also recall the notation $\bout_k(\entrance) \in \Psi_k$
from Definition \ref{defpsik}, and abbreviate $\ok = \bout_k(\entrance)$.

There exists a constant $c > 0$, without dependence on $(B,\singent)$, such that, for each $k \in \mathbb{N}$,
there exists a random variable 
$\halfcorrfack \in \sigma \big\{ \ik, \ok \big\}$ \hfff{approxcf}
satisfying 
\begin{equation}\label{fkrel}
  \halfcorrfac - \exp \big\{ -ck \big\} \leq  \halfcorrfack \leq \halfcorrfac.
\end{equation}
\end{lemma} 
We omit the proof of the following fact.
\begin{lemma}\label{lemebfbd}
For any backbone-tree pair, we have the bounds
$\frac{\binf - 1}{\bsup + \binf + 1}  \leq  \pesc \leq 1$ and
$1 - \binf^{-1} \leq \pde \leq 1$.
Further, $c \leq \halfcorrfac \leq C$, where $c = 1$ and $C = 
 \frac{\bsup + \binf + 2}{\binf - 1}$.
\end{lemma}
\noindent{\bf Proof of Lemma \ref{lembt}.} We will show that there exist
$\ekm \in \sigma \big\{ \ok \big\}$
and 
$\qkm \in \sigma \big\{ \ik \big\}$ 
for which
\begin{equation}\label{edata}
 \pesc + 2 \binf^{1-k/2}  \geq \ekm \geq \pesc
\end{equation}
and
\begin{equation}\label{qdata}
 \pde + 2 \binf^{-k}  \geq \qkm \geq \pde.
\end{equation}
We then set $\halfcorrfack = \ekm^{-1} + \qkm^{-1} - 1$. The bounds in (\ref{fkrel})
then follow from  (\ref{edata}), (\ref{qdata}) and Lemma \ref{lemebfbd}.

We define
$$
\ekm = \P_{(B,\singent),\beta}^\head \Big( H_{V(\ok) \setminus V(\okmo)}  < H_\entrance \Big).
$$
Let $$
\Big\{ p_v :  v \in \big( V(\ok) \setminus V(\okmo)  \big) \cup \{
\entrance \} \Big\}
$$
denote the hitting distribution on $\big( V(\ok) \setminus V(\okmo)
\big) \cup \{ \entrance \}$ under $\P_{(B,\singent),\beta}^\head$. We have then that
\begin{equation}\label{eke}
\ekm - \pesc  = \sum_{v \in V(\ok) \setminus V(\okmo)} p_v \P_{(B,\singent)}^v
 \Big( \xsub{n} = \entrance \, \,  \textrm{for some $n > 0$} \Big).
\end{equation}
\begin{figure}
\begin{center}
\includegraphics[width=0.3\textwidth]{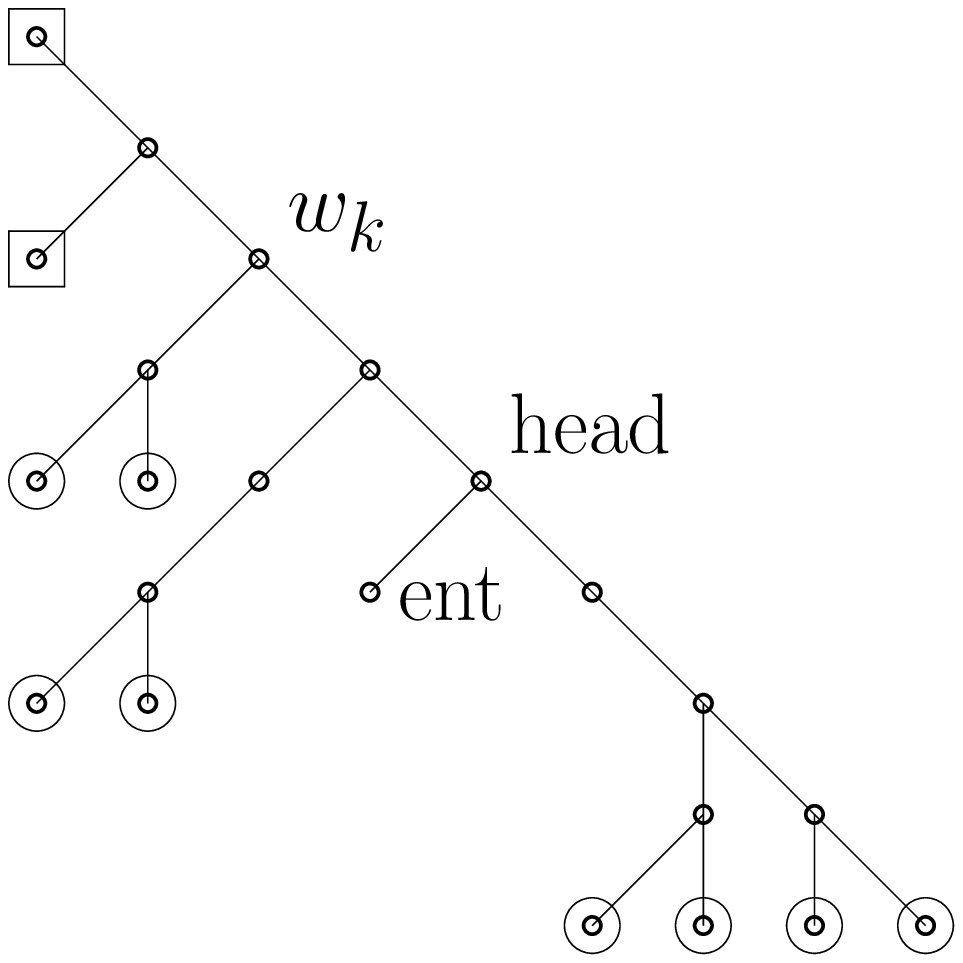} \\
\end{center}
\caption{An instance of $\ok = \bout_k(\entrance)$ with $k=4$. The circles and squares indicate elements of $\dk$ and $\ekd$.}\label{fignkout}
\end{figure}
Figure \ref{fignkout} illustrates the objects which we now use. Let $w_k \in V(\ok)$ denote the element of the backbone $B$ 
whose $\lfloor k/2 \rfloor$-th descendent is $\head$.
We write $V(\ok) \setminus V(\okmo) = \dk \cup \ekd$,
where $\dk$ is the set of $v \in V(\ok) \setminus V(\okmo)$
that are descendents of $w_k$, and $\ekd = \big( V(\ok) \setminus V(\okmo) \big)
\setminus \dk$. If $v \in \dk$, then the unique simple path from $v$ to
$\head$ in $B$, which is necessarily of length $k$, begins by traversing at least $\lfloor k/2
\rfloor$ edges consecutively in the direction from offspring to parent. From this, we see that
\begin{equation}\label{eqak}
\P_{(B,\singent)}^v \Big(  X(n) = \entrance \, \,  \textrm{for some $n > 0$} \Big)
 \leq  \binf^{-\lfloor k/2 \rfloor},
\end{equation}
for $v \in \dk$.  Indeed, the inequality, bounding the probability that biased walk on the backbone ever visits the ancestor at distance $\lfloor k/2 \rfloor$ of its point of departure, is proved by comparison with a biased random walk on $\Z$.

If $v \in \ekd$, on the other hand, the simple path from $\head$ to $v$ in $B$
begins with the $\lfloor k/2 \rfloor$-length 
path from $\head$ to $w_k$, each of whose edges are traversed in the direction from offspring 
to parent.
Hence,
\begin{eqnarray}
 & & \sum_{v \in \ekd} p_v \leq 
\P_{(B,\singent)}^\head \Big(  X(n) \in \ekd \, \,  \textrm{for some $n > 0$} \Big) \nonumber \\
 & \leq &
\P_{(B,\singent)}^\head \Big(  X(n) = w_k \, \,  \textrm{for some $n > 0$} \Big) \leq
\binf^{-k} \label{eqbk} 
\end{eqnarray}
for $v \in \ekd$.

By (\ref{eqak}), (\ref{eqbk}) and (\ref{eke}), then, we see that
$\big\vert \pesc - \ekm \big\vert \leq 2 \binf^{- \lfloor k/2 \rfloor}$. Noting that $\ekm - \pesc  \geq 0$ by (\ref{eke}), we obtain (\ref{edata}).
 
Recall that $\btbase$ is a vertex in $T_\entrance$ at maximal distance from $\entrance$, and that
$\pde = \P_{(B,\singent)}^{\entrance}\big( H_\btbase < H_\head \big)$. 
If $r(T) \leq k$, we simply set $\qkm = \pde$. Otherwise, recalling Definition \ref{defndiv}, abbreviate $C_k = C_k(T_\entrance)$, and
 define 
$$
\qkm = \P_{(B,\singent)}^{\entrance} \Big( H_{\basedef(C_k)} < H_\head \Big). 
$$
By $d\big( \entrance, \basedef(C_k) \big) \geq k$, it follows readily that 
$0 \leq \pde - \qkm \leq \binf^{-k}$ holds. \flushright \qed
\end{subsubsection}
\begin{subsubsection}{The analogue of Proposition \ref{trprop} for the approximate conditioning}
\begin{lemma}\label{qklem}
Let $k \in \N$.
We have that
$$
\qstarnoarg \Big( \corrfack \omen(T_{\entrance}) > u \Big) \sim \cpo \olddone(k) u^{- \chiexp},
$$
as $u \to \infty$,
where 
$$
\olddone(k) = \big( 1 - p_0/\pext \big)^{-1} \E_{\qstarnoarg}
  \bigg( \Big( \corrfack \omen \big( \base(C_k) \big) \Big)^{\chiexp}
  1\!\!1_{r(T_{\entrance}) > k} \bigg),
$$
where $\cpo > 0$ and $\chiexp > 0$ were specified in Proposition \ref{theoremone}.
\end{lemma}
To prove this result, we need a simple lemma.
\begin{lemma}\label{lemlit}
$\empty$
\begin{enumerate}
\item If $U$ and $V$ are independent random variables on a probability space
$(\Omega,\mathbb{P})$, with $V \geq 1$ a.s.,
\begin{equation}\label{uuchi}
 \mathbb{P}\big(U > u\big) \sim c_1 u^{-\chiexpdiff}
\end{equation}
as $u \to \infty$ for some $c_1 > 0$ and $\chiexpdiff > 0$, and, for some $\eta > 0$, 
\begin{equation}\label{vuchi}
 \mathbb{P}\big( V > u \big) \leq u^{-\eta}
 \mathbb{P}\big( U > u \big) 
\end{equation}
for $u$ sufficiently high, 
then, as $u \to \infty$,
\begin{equation}\label{uvchi}
 \mathbb{P}\big( U V > u \big) \sim c_1 \E \big( V^{\chiexpdiff} \big) u^{-\chiexpdiff}. 
\end{equation}
\item Let $U$ and $V$ be non-negative random variables on a probability space $\big( \Omega, \mathbb{P} \big)$.
Suppose that $\mathbb{P} \big( U \geq u \big) \sim c_1 u^{-\chiexpdiff}$ for some
$c_1 > 0$ and $\chiexpdiff > 0$, and also that
\begin{equation}\label{vuratio}
\lim_{u \to \infty} \frac{\mathbb{P} \big( V \geq u \big)}{\mathbb{P} \big( U \geq u \big)} = 0.
\end{equation}
Then 
$$
\mathbb{P} \big( U + V \geq u \big)  \sim c_1 u^{-\chiexpdiff}.
$$
\end{enumerate}
\end{lemma}
\noindent{\bf Proof.} The first statement is Lemma 10 in \cite{paperone}. The simple proof of the second is omitted. 
  \qed \\
\noindent{\bf Proof of Lemma \ref{qklem}.}
We begin by showing that the number of renewal components under the approximate conditioning is asympotically infinite, i.e., for each $k \in \N$,
\begin{equation}\label{qstfku}
\lim_{u \to \infty} \qstarnoarg \Big( r(T_{\entrance}) \leq k \Big\vert \corrfack \omen \big( T_{\entrance} \big) > u \Big) = 0.
\end{equation}
To check this, note that 
\begin{eqnarray}
 & & \qstarnoarg \Big( r(T_{\entrance}) \leq k \Big\vert \corrfack \omen(T_{\entrance}) > u \Big) \label{qstcom} \\
 & = & \frac{\qstarnoarg \Big( \Big\{ r(T_{\entrance}) \leq k \Big\} \cap \Big\{ \corrfack \omen(T_{\entrance}) > u
    \Big\}  \Big)}{\qstarnoarg \Big( \corrfack \omen(T_{\entrance}) > u  \Big)}
  \nonumber \\
 & \leq & \frac{\P_{h,\nu} \Big( \Big\{ r(T) \leq k \Big\} \cap \Big\{ C \omega(T) > u
    \Big\}  \Big)}{\P_{h,\nu} \Big( c \omega(T) > u  \Big)} \nonumber \\
 & \leq &  2 \big( C/c \big)^{\chiexp} \P_{h,\nu} \Big( r(T) \leq k \Big\vert \omega(T) > u/C \Big). \nonumber
\end{eqnarray}
The first inequality here follows from the bounds $c \leq \halfcorrfack \leq C$, which Lemmas \ref{lembt} and \ref{lemebfbd} imply, as well as the law of $T_{\entrance}$ under $\qstarnoarg$ being $\P_{h,\nu}$.
The second inequality is due to Proposition  \ref{theoremone}. It follows directly from Lemma \ref{lemdivexp}
that
\begin{equation}\label{eqlimuom}
\lim_{u \to \infty} \P_{h,\nu}\Big( r(T) \leq k \Big\vert \omega(T) > u \Big) = 0 \quad  \textrm{for each $k \in \N$}.
\end{equation}
We thus find that, indeed, (\ref{qstfku}) holds for each $k \in \N$. 

Note that, in light of  (\ref{qstfku}), the statement of the lemma for some given $k \in \N$ is equivalent to
\begin{eqnarray}
 & & \qstarnoarg \Big(  \corrfack \omen(T_{\entrance}) > u \Big\vert  r(T_{\entrance}) > k
\Big)  \nonumber \\
 & \sim &
\cpo \big( 1 - p_0/\pext \big)^{-1} \E_{\qstarnoarg}
  \bigg( \Big( \corrfack \omen \big( \base(C_k) \big) \Big)^{\chiexp}
 \Big\vert r(T_{\entrance}) > k \bigg)
 u^{-\chiexp}. \label{qstcpr}
\end{eqnarray}
We now prove (\ref{qstcpr}). We partition the trap weight $\omen (T_{\entrance})$ according to contributions by vertices above and below $c_k$:  provided that $r(T_{\entrance}) > k$, note that
\begin{equation}\label{omaluv}
\omen (T_{\entrance}) = \alpha_k(T_{\entrance}) + u_{(k)}(T_{\entrance}) v_{(k)}(T_{\entrance}),
\end{equation}
where here we define 
$\alpha_k\big(T_\entrance\big) = \omen \big( S_k(T_{\entrance})  \big) -  \omen \big( \basedef(C_k) \big)$,
$u_{(k)}\big(T_\entrance\big) =  \omen \big( \basedef(C_k) \big)$
and 
$$
v_{(k)}\big(T_\entrance\big) = \omega_{\basedef(C_k)} \big( T_{\basedef(C_k)} \big).
$$
The notation of Definition \ref{defomega} is used here, for example in the 
last definition, of the weight of $T_{\basedef(C_k)}$ relative to its root $\basedef(C_k)$. In (\ref{omaluv}), we indeed have such a partition of trap weight, in the sense that 
\begin{equation}\label{meaone}
 \alpha_k(T_{\entrance}) \in \sigma \big\{ \ik  \big\}, \, \,   u_{(k)}(T_{\entrance}) \in \sigma \big\{
 \ik   \big\},
\end{equation}
and
\begin{equation}\label{meatwo}
v_{(k)}(T_{\entrance}) \in \sigma \Big\{ T_{\entrance} \setminus \ik \Big\};
\end{equation}
the latter is due to
$v_{(k)}(T_{\entrance})$ being measurable with respect to $T_{\basedef(C_k)}$.

We will derive (\ref{qstcpr}) by multiplying (\ref{omaluv}) by $\corrfack$ and applying both parts of Lemma \ref{lemlit} to analyse the tail of its right-hand side.
To begin this application, note firstly that,
since $T_\entrance$ under $\qstarnoarg$ has the law of $\P_{h,\nu}$, Lemma \ref{lemdecom}
implies the following equality in distribution:
\begin{equation}\label{eqeqdist}
 \qstarnoarg \Big( v_{(k)}\big(T_{\entrance}\big) \geq x \Big\vert r(T_{\entrance}) > k \Big)
 = \P_{h,\nu}\Big( \omega(T) \geq x \Big\vert \vert E(T) \vert \geq 1 \Big)
 \quad \textrm{for all $x > 0$.}
\end{equation}

We thus 
see by Proposition \ref{theoremone} and $h_0 = p_0/\pext$ that 
\begin{equation}\label{eqqsim}
 \qstarnoarg \Big( v_{(k)}(T_{\entrance}) > u \Big\vert r(T_{\entrance}) > k  \Big)
 \sim \cpo \big(1- p_0/\pext \big)^{-1} u^{-\chiexp}.
\end{equation}

We apply 
Lemma \ref{lemlit}(1) with the choices $\P = \qstarnoarg \big( \cdot \big\vert r(T_{\entrance}) >
k   \big)$, $U =  v_{(k)}(T_{\entrance})$ and $V =   \corrfack \omen \big( \basedef(C_k) \big)$; admitting for now that its hypotheses are satisfied, and using (\ref{eqqsim}) after applying the lemma, we learn that 
\begin{eqnarray}
 & & \qstarnoarg \Big( \corrfack u_{(k)}(T_{\entrance}) v_{(k)}(T_{\entrance}) > u \Big\vert r(T_{\entrance}) >
k  \Big) \label{eqcorrfack} \\ 
& \sim & \cpo \big(1- p_0/\pext \big)^{-1}   \E_{\qstarnoarg}
  \bigg( \Big( \corrfack \omen \big( \basedef(C_k) \big) \Big)^{\chiexp}
  \Big\vert r(T_{\entrance}) > k \bigg) u^{-\chiexp}. \nonumber
\end{eqnarray}
We now complete the proof of (\ref{qstcpr}): multiply (\ref{omaluv}) by $\corrfack$, and then apply Lemma \ref{lemlit}(2) for $\P$ as before and with $V$ and $U$ the respective terms on the right-hand side. By (\ref{eqcorrfack}), we obtain (\ref{qstcpr}) and thus reduce the proof of the lemma to verifying that the hypotheses for the two applications of Lemma \ref{lemlit} are satisfied.

For the application of Lemma \ref{lemlit}(1), note that, by  (\ref{meaone}), (\ref{meatwo})
and $\halfcorrfack \in \sigma \big\{ \ik,\ok \big\}$, 
the random variables
$\halfcorrfack u_{(k)}(T_{\entrance})$ and $v_{(k)}(T_{\entrance})$ are indeed independent under $\qstarnoarg \big( \cdot
\big\vert r(T_{\entrance}) > k \big)$. Moreover, (\ref{eqqsim}) validates the hypothesis (\ref{uuchi}). For the remaining hypothesis (\ref{vuchi}), 
it suffices to argue that there exists $\epsilon > 0$ such that, for $u$ sufficiently
high,
\begin{equation}\label{qstu}
\qstarnoarg \Big( \corrfack u_{(k)}(T_{\entrance}) > u \Big\vert r(T_{\entrance}) > k  \Big)
\leq u^{-\epsilon} 
\qstarnoarg \Big( v_{(k)}(T_{\entrance}) > u \Big\vert r(T_{\entrance}) > k  \Big)
\end{equation}
We do so momentarily. Regarding the use of Lemma \ref{lemlit}(2), 
we must confirm (\ref{vuratio}), and for this it suffices to show that 
\begin{eqnarray}
& & \qstarnoarg \Big( \corrfack \alpha_k(T_{\entrance}) > u \Big\vert r(T_{\entrance}) > k  \Big) \nonumber \\
& \leq & u^{-\epsilon} 
\qstarnoarg \Big(  \corrfack u_{(k)}(T_{\entrance}) v_{(k)}(T_{\entrance}) > u \Big\vert r(T_{\entrance}) > k  \Big).\label{qsta}
\end{eqnarray}

We have reduced to showing (\ref{qstu}) and (\ref{qsta}).
For this, note that Lemma~\ref{lemdivexp}
implies that, $k$ being fixed, there exists $\epsilon > 0$
such that, for $u$ sufficiently high,
\begin{equation}\label{eqbom}
\P_{h,\nu} \Big(  \Big\{ \Big\vert \bigcup_{i=1}^k V\big( C_i \big) \Big\vert \geq
\frac{\log u}{\log (2\bsup)} \Big\} \cap \Big\{ r(T) > k \Big\} \Big\vert \omega(T) > u \Big)
 \leq 
u^{-\epsilon}. 
\end{equation}
Noting that
$$
\max \Big\{ \alpha_k(T), u_{(k)}(T) \Big\}
 \leq \Big\vert \bigcup_{i=1}^k V(C_i) \Big\vert \bsup^{\big\vert
   \bigcup_{i=1}^k V(C_i) \big\vert} \leq \big( 2 \bsup \big)^{\big\vert \bigcup_{i=1}^k V(C_i) \big\vert},
$$
and
$\max \big\{ \alpha_k(T), u_{(k)}(T) \big\} \leq \omega(T)$,
we find from (\ref{eqbom}) that
$$
\frac{\P_{h,\nu }\Big( \Big\{ \max \big\{ \alpha_k(T), u_{(k)}(T) \big\} > u \Big\} \cap \Big\{
 r(T) > k \Big\} \Big)}{\P_{h,\nu} \Big( \omega(T) > u \Big)}
\leq 
u^{-\epsilon}, 
$$
or, equivalently, 
\begin{eqnarray}
& & \P_{h,\nu }\Big( \max \big\{ \alpha_k(T), u_{(k)}(T) \big\} > u 
\Big\vert r(T) > k  \Big) \label{eqalphak} \\
& \leq &  \P_{h,\nu} \Big( r(T) > k \Big)^{-1}  u^{-\epsilon} \, \P_{h,\nu} \Big( \omega(T)
> u\Big). \nonumber
\end{eqnarray}
By $\halfcorrfack \leq C$, (\ref{eqeqdist}) and  Proposition \ref{theoremone}, we obtain (\ref{qstu}), by adjusting the value of $\epsilon > 0$.
Regarding (\ref{qsta}), note that, by $\halfcorrfack \geq c$ and $v_{(k)}(T) \geq 1$ $\qstarnoarg$-a.s., the probability on its right-hand side is at least 
$\qstarnoarg \big(   2c  u_{(k)}(T_{\entrance}) > u \big\vert r(T_{\entrance}) > k  \big)$, which, by (\ref{eqeqdist}), is at least
$c  \P_{h,\nu}\big( \omega(T) \geq u/(2c) \big)$. The regularity of the tail of this last term described in  
 Proposition \ref{theoremone} shows that we may replace the term by 
$c  \P_{h,\nu}\big( \omega(T) \geq u \big)$ at the expense of altering $c > 0$. All this means that, to obtain (\ref{qsta}), it is enough to verify that its left-hand side is bounded above by $c u^{-\epsilon}  \P_{h,\nu}\big( \omega(T) \geq u \big)$. Recalling that $T_\entrance$ under $\qstarnoarg$ has the law $\P_{h,\nu}$, and using
 $\halfcorrfack \leq C$, we see that the left-hand side of (\ref{qsta}) is at most
 $\P_{h,\nu} \big( \alpha_k(T) > u/C \big\vert r(T) > k \big)$, so that (\ref{eqalphak}), and an adjustment in the value of $\epsilon > 0$, yields what we need.   \qed
\end{subsubsection}
\begin{subsubsection}{Completing the proof of Proposition \ref{trprop}}
Another general lemma is needed. In this case, the proof is trivial and omitted.
\begin{lemma}\label{lrgfour}
Let $\chiexpdiff > 0$.
Let $s:[0,\infty) \to [0,\infty)$, 
$\big\{ s_{\epsilon}: [0,\infty) \to [0,\infty) , \epsilon > 0 \big\}$
and the collection $\big\{ c_{\epsilon}:
\epsilon > 0 \big\}$ of constants be such that, for all $\epsilon > 0$, 
$$
s_{\epsilon}(u) \sim  c_{\epsilon} u^{- \chiexpdiff}
$$
and
$$
 s_{\epsilon}(u)  \leq  s(u) \leq  s_{\epsilon}\big(u(1 + \epsilon)\big)
$$
for all $u \geq u_0(\epsilon)$ sufficiently high. 
Then $c = \lim_{\epsilon \downarrow 0} c_{\epsilon}$ exists and
$$
 s(u) \sim c u^{-\chiexpdiff}. 
$$
\end{lemma}
\noindent{\bf Proof of Proposition \ref{trprop}.} Note that
\begin{eqnarray}
 & & \qstarnoarg \Big( \corrfack \omen(T_{\entrance}) > u \Big) \leq
 \qstarnoarg \Big( \corrfac \omen (T_{\entrance}) > u \Big) \nonumber \\
 & \leq &
 \qstarnoarg \Big( \corrfack \big( 1 + c^{-1} \exp \big\{ - c k \big\} \big) 
 \omen(T_{\entrance}) > u \Big), \nonumber 
\end{eqnarray}
the first inequality due to $\halfcorrfack \leq \halfcorrfac$, and the second to
$\halfcorrfack \geq c$ and $\vert \halfcorrfac - \halfcorrfack \vert \leq \exp \{ - ck \}$ (which bounds are shown   for high choices of $k$ by Lemma \ref{lembt} in tandem with Lemma \ref{lemebfbd}).
The result then follows from Lemmas \ref{qklem} and  \ref{lrgfour}, with 
$$
\olddone =  \big( 1 - p_0/\pext \big)^{-1} \lim_{k \to \infty} \E_{\qstarnoarg}
  \bigg( \Big( \corrfack \omen \big( \basedef(C_k) \big) \Big)^{\chiexp}
  1\!\!1_{r(T_{\entrance}) > k} \bigg).
$$
We may replace $\halfcorrfack$ by $\halfcorrfac$ in this formula, due to $\halfcorrfac \geq c$ and
(\ref{fkrel}). \qed
\end{subsubsection}
\end{subsection}
\begin{subsection}{Deriving Proposition \ref{propqkq} and Corollary \ref{corqkq}}\label{secproofs}
The proof is undertaken in several stages:
\begin{itemize}
 \item in Subsection \ref{subsecone}, we reformulate  Proposition \ref{propqkq} as the assertion that a certain coupling exists;
 \item in \ref{subsectwo} and \ref{subsecthree}, we introduce notation and construct this coupling; 
 \item in \ref{subsecfour}, we prove that the coupling has the required properties; 
 \item in \ref{subsecfive}, we prove two lemmas that were needed in \ref{subsecfour}; and, in \ref{subsecsix}, we explain the change needed to obtain Corollary \ref{corqkq}. 
\end{itemize}
\subsubsection{Rephrasing Proposition \ref{propqkq} in terms of a coupling}\label{subsecone}
Let $\big( T^{[1]},X_{[1]} \big)$ and $\big( T^{[2]},X_{[2]} \big)$
each have the distribution  $\pf$.
Recalling Definitions \ref{deftxb} and \ref{defpsik}, 
we write 
$\big\{ E_k^{[i]}(n): n \in \N \big\}$,
$i \in \{ 1,2 \}$, 
for the random variables $\big\{ \bout_k\big(X_{[i]}(\trar{n})\big): n \in \N \big\}$
arising from the sample $\big( T^{[i]} , X_{[i]} \big)$.
\begin{prop}\label{propcoup}
For each  $k \in \N$ and $\epsilon >
 0$, 
there exists
 $n_0 \in \N$ with the following property. Whenever $n_1,n_2 \geq n_0$, there exists a coupling $\coupl_{n_1,n_2}$ of  $\big( T^{[1]} , X_{[1]} \big)$ and
 $\big( T^{[2]} , X_{[2]} \big)$ such that
$$
 \Theta_{n_1,n_2} \Big( E_k^{[1]}(n_1) \not= E_k^{[2]}(n_2) \Big)
 \leq \epsilon.
$$ 
\end{prop}
\noindent{\bf Proof of Proposition \ref{propqkq}.}
Proposition \ref{propcoup} implies that,
for each $k \in \N$ and for any $\epsilon > 0$, there exists $n_0 \in \N$
such that, for $n_1,n_2 \geq n_0$,
\begin{equation}\label{eqqkpair}
 {\rm TV} \Big(  \qsubk(n_1) , \qsubk(n_2) \Big) \leq  \epsilon.
\end{equation}
It is a simple matter to verify that 
the space of probability measures on an arbitrary measure space, endowed with the total variation metric, is complete. 
Thus,~(\ref{eqqkpair}) implies Proposition~\ref{propqkq}. \qed

\subsubsection{Preliminaries and an overview for constructing $\Theta_{n_1,n_2}$}\label{subsectwo}

In seeking to construct the coupling $\Theta_{n_1,n_2}$ in Proposition \ref{propcoup}, we may,  without loss of generality, take $n_1 > n_2$.

We will construct $\Theta_{n_1,n_2}$ by beginning with an independent coupling of the two marginal processes and making an appropriate modification. Let $\big(\Omega,\Theta'\big)$ denote a probability space where the measure $\Theta'$ is the product coupling of two copies of the law $\pf$. We will refer to the marginals of $\Theta'$ as $(T^{\squone},X_{\squone})$ and $(T^{\squtwo},X_{\squtwo})$; we also write ${\rm RG}_{i}$
 and $\big\{ \trar{j}^{i}: j \in \N^+ \big\}$,  $i \in \big\{ 1,2 \big\}$, for the set of regeneration times and the sequence of trap entrance arrivals of the two processes.

Some further notation is needed.
\begin{definition}\label{defprop}
Under $\pf$, enumerate ${\rm RG} = \big\{ r_i: i \in \N^+ \big\}$.
For $\ell \in \N^+$, define the ``tree until time $r_\ell$'', $\mathcal{T}_\ell$, to be the weighted tree given by the connected component of $X(r_\ell) \in B(T)$ in the graph
 $(V(T),E')$ whose edge-set $E'$ is formed by removing from $E(T)$ the set of elements of $V(T)$ that connect $X(r_\ell)$ to the offspring of this vertex; informally,  $\mathcal{T}_\ell = T \setminus T_{X(r_\ell)}$.
  We define the  ``walk until time $r_\ell$''
to be $X:\big\{ 0, \ldots r_\ell \big\} \to V(\mathcal{T}_\ell)$. 
We then define the ``history until time $r_\ell$'' 
 to be the ordered pair given by the tree and the walk until this time.
\end{definition}
\noindent{\bf Remark.} Each regeneration time separates the walk's future from its past. In constructing~$\Theta_{n_1,n_2}$, we will exploit the following easily verified but essential property of ${\rm RG}$ under $\pf$. For each $\ell \in \N^+$, given the history until time $r_\ell$, the conditional distribution of the future environment $T_{X(r_\ell)}$ and the future walk $X(r_\ell + \cdot):\N \to V\big(T_{X(r_\ell)} \big)$ is the law of $(T,X)$ under 
\begin{equation}\label{pfzerorg}
 \Big( \pf \Big) \Big(  \cdot \Big\vert 0 \in {\rm RG} \Big).
\end{equation}

In studying $\Theta'$,  we will use the term first history (and walk) until a given element in  ${\rm RG}_{\squone}$ in the sense of Definition \ref{defprop} applied to the first marginal of $\Theta'$; and similarly, of course, for the second marginal.

For $i \in \{ 1,2 \}$ and $n \in \N$, we record the number of trap entrances found by $X_{i}$ (under $\Theta'$) before time $n$,
$$
h_i(n) = \sup \big\{  j \in \N : \trar{j}^{i} < n \big\}.
$$
A pair $(t_1,t_2) \in \big( {\rm RG}_{\squone}, {\rm RG}_{\squtwo} \big)$ is said to be good if 
$$
 h_1(t_1) = h_2(t_2) + \big( n_1 - n_2 \big).
$$
A good pair $(t_1,t_2)$ is a useful concept for the following reason. Write $s = n_1 - h_1(t_1)$, which is also $n_2 - h_2(t_2)$. Suppose that $s > 0$, and consider the conditional distribution of $\Theta'$ 
given the first history until time $t_1$ and the second until time $t_2$. Then the $n_1$-st trap entrance for the first walk will be the $s$-th to be encountered subsequently, while, for the second walk, it is the $n_2$-nd trap entrance that will be the $s$-th to be so encountered. Our procedure for constructing $\Theta_{n_1,n_2}$ will rely on this fact. The ``earliest'' good pair under $\Theta'$ will be located; then, under $\Theta'$ given the two histories until the respective times in the given good pair, the two conditional futures will be the same, in view of the remark made a few moments ago. We may then modify $\Theta'$ to define $\Theta_{n_1,n_2}$, by insisting that these walk futures after the respective times be equal under the new measure. In so doing, we expect that the $k$-neighbourhood environments of the walks in the two marginals to be equal $s$ trap entrances later; and this is what we need to prove Proposition~\ref{propcoup}, because these environments have the laws of $E_k^{[1]}(n_1)$
and $E_k^{[2]}(n_2)$. For this plan to work, we must have $s > 0$, in fact $s \gg k$; but, as we will see, this is reasonable if $n_2$ is large.

\subsubsection{Constructing $\Theta_{n_1,n_2}$ by playing {\rm \bf LEAPFROG}}\label{subsecthree}

We now define a stopping-time procedure that selects a good pair under $\Theta'$. Under the procedure, called {\rm LEAPFROG}, the two histories under $\Theta'$ are alternately sampled until a good pair is located. 

The first stage of {\rm LEAPFROG} begins by sampling the first history until time $\reg_1^{\squone}$, where  $\reg_1^{\squone}$ denotes the first element of ${\rm RG}_{\squone}$ that exceeds 
 $\trar{n_1 - n_2}^{\squone}$. In this history, $h_1\big(\reg_1^{\squone}\big)$ trap entrances have been witnessed; there is a non-negative ``overshoot'' of 
$h_1\big(\reg_1^{\squone}\big) - (n_1 - n_2)$. The second history is then sampled until  time $\reg_1^{\squtwo}$, where  $\reg_1^{\squtwo}$ denotes the smallest element of ${\rm RG}_{\squtwo}$ that exceeds 
 $\trar{h_1\big(\reg_1^{\squone}\big) - (n_1 - n_2)}^{\squtwo}$.  
The second history has leapfrogged the first, and its trap entrance counter has a non-negative overshoot given by $h_2 \big( \reg_1^{\squtwo} \big)$ minus the first history overshoot. If the overshoot is zero, in other words, if
\begin{equation}\label{htwoeq}
 h_2 \big( \reg_1^2 \big) = h_1 \big( \reg_1^1 \big) - \big( n_1 - n_2 \big), 
\end{equation} 
then {\rm LEAPFROG} terminates here; otherwise, it continues to the second stage.

We now describe {\rm LEAPFROG}'s generic step.
For $\ell \geq 2$, at the start of the $\ell$-th stage, the two histories have been sampled until times
$\reg_{\ell-1}^{\squone} \in {\rm RG}_{\squone}$ and $\reg_{\ell-1}^{\squtwo} \in {\rm RG}_{\squtwo}$. 
If the $\ell$-th step is to take place at all,
then, necessarily, 
$$
 h_2 \big( \reg_{\ell-1}^{\squtwo} \big) \not=  h_1 \big( \reg_{\ell-1}^{\squone} \big) - \big( n_1 -
 n_2 \big). 
$$
Without loss of generality, the left-hand-side is the greater. The first process will now leapfrog the second, and then, if need be, the second the first.  We extend the sampling of the first history  until time $\reg_\ell^{\squone}$, which we define to be the smallest element of ${\rm RG}_{\squone}$ that exceeds
$\trar{h_2(\reg_{\ell-1}^{\squtwo}) + n_1 - n_2}^{\squone}$.
We then extend the sampling of the second history  until time $\reg_\ell^{\squtwo}$, which is the smallest element 
of ${\rm RG}_{\squtwo}$ that exceeds
$\trar{h_1(r_\ell^1) - ( n_1 - n_2 )}^{\squtwo}$. (Note that, in the case that $\reg_\ell^{\squtwo} = \reg_{\ell -1}^{\squtwo}$, no non-trivial extension is made here.)
If the condition
\begin{equation}\label{htwort}
 h_2 \big( \reg_\ell^{\squtwo} \big) =  h_1 \big( \reg_\ell^{\squone} \big) - \big( n_1 -
 n_2 \big)
\end{equation}  
is met, then {\rm LEAGFROG} terminates. Otherwise, it continues to 
its $(\ell+1)$-st stage.

If {\rm LEAPFROG} terminates at some finite stage, we write $\ell^* \in \N^+$ for the index of the terminating stage.

We define under $\Theta'$ the random sets
$H_i = \big\{ h_i(t): t \in {\rm RG}_{i} \big\}$
for $i \in \{ 1,2 \}$. 
In these terms, note that there is a good pair of the form 
$(t_1^*,t_2^*)$ where
\begin{equation}\label{eqhone}
h_1\big( t_1^* \big) = \inf \bigg( H_1 \cap \Big( H_2 + (n_1 - n_2) \Big) \bigg)
\end{equation}
and
\begin{equation}\label{eqhtwo}
h_2\big( t_2^* \big) = \inf \bigg( \Big( H_1 - (n_1 - n_2 ) \Big)  \cap H_2  \bigg),
\end{equation}
(provided that the infimum is in fact over a non-empty set).

We claim that $H_1 \cap \big( H_2 + (n_1 - n_2) \big)  \not=\emptyset$ if and only if {\rm LEAPFROG} terminates. Moreover, if this is so, then {\rm LEAPFROG} locates this particular good pair, in the sense that $t_1^* = \reg_{\ell^*}^{\squone}$ and  $t_2^* = \reg_{\ell^*}^{\squtwo}$. See Figure \ref{figleapfrog} for an explanation of why this claim holds.

\begin{figure}
\begin{center}
\includegraphics[width=0.75\textwidth]{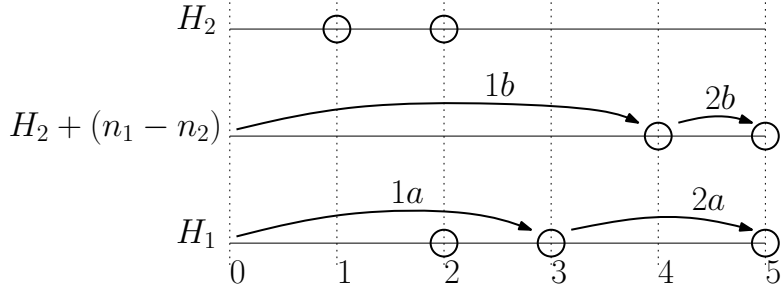} \\
\end{center}
\caption{A little game of {\rm LEAPFROG}: in this example, $n_1 = 7$ and $n_2 = 4$. In stage $1$, the first history is constructed so that the first counter reaches $h_1(\reg_1^{\squone}) = 3$ (arrow $1a$); then, in arrow $1b$, the second history leapfrogs so that the second counter reaches $h_2(\reg_1^{\squtwo}) + (n_1 - n_2) = 1 + (7 - 4) = 4$. In stage $2$: (2a) the first history leapfrogs so that $h_2(\reg_2^{\squone}) = 5$; then (2b) the second leapfrogs to the same position, 
$h_2(\reg_2^{\squtwo}) + (n_1 - n_2) = 5$. The terminal stage is thus $\ell^*=2$, with both counters then at the first meeting point of the sets in the second and third rows, i.e., at $h : = \inf \big\{ H_1 \cap \big( H_2 + (n_1 - n_2) \big) \big\}$. That $h < \infty$ implies that the counters meet at $h$ may be established by induction on $h$; the simple argument is omitted.}\label{figleapfrog}
\end{figure}

Let $\sigma^*$ denote the $\sigma$-algebra of the probability space $\big( \Omega, \Theta'\big)$ generated by the first history until time $\reg_{\ell^*}^{\squone}$ and the second history until time $\reg_{\ell^*}^{\squtwo}$. The random variable $\ell^*$ being a stopping time for the procedure, it follows from the remark after Definition~\ref{defprop} that, 
$\Theta' \big( \cdot \big\vert \sigma^* \big)$-almost surely, the conditional joint distribution of
$$
\Big(T^{i}_{X(\reg_{\ell^*}^{i})}, X \big( \reg_{\ell^*}^{i} + \ldots  \big): \N \to T^{i}_{X(\reg_{\ell^*}^{i})} \Big)
$$
has the distribution (\ref{pfzerorg}), for both $i=1$ and $2$. We augment the probability space $(\Omega,\Theta')$ to include a (tree,walk)-valued random variable $(\tilde{T},\tilde{X})$ with law (\ref{pfzerorg}) which is independent of all existing data. We then define $(T^{[i]},X_{[i]})$ for $i \in \{1,2\}$ as follows.
For the process $(T^{i},X_{i})$, recall the tree until time~$\reg_{\ell^*}^{i}$  from Definition \ref{defprop}. Rather than call this tree the unwieldy $\mathcal{T}_{\reg_{\ell^*}^{i}}^{i}$,
we use the shorthand $\mathcal{T}_{\ell^*}^{i}$.
Let $T^{[i]}$ be formed from the trees  $\mathcal{T}_{\ell^*}^{i}$  and $\tilde{T}$
by identifying their respective vertices $X\big(\reg_{\ell^*}^{i}\big)$  and 
$\phi$. In a slight abuse of notation, we may thus regard $\mathcal{T}_{\ell^*}^{i}$  and $\tilde{T}$
as induced sub-trees in $T^{[i]}$.  We then set $X_{[i]}: \N \to T^{[i]}$ according to 
\begin{equation}
X_{[i]}(j)  = \left\{ \begin{array}{rl} X_{i}(j) &  \textrm{for $0 \leq j \leq \reg_{\ell^*}^{i}$,} \\
          \tilde{X}_{i}(j) &  \textrm{for $j > \reg_{\ell^*}^{i}$.} \end{array} \right. 
\end{equation}

The coupling $\Theta_{n_1,n_2}$ is then defined to be the joint distribution of $(T^{[i]},X_{[i]})$ for $i \in \{ 1,2 \}$.

\subsubsection{Deriving Proposition \ref{propcoup}}\label{subsecfour}

We now prove Proposition \ref{propcoup} by verfying that $\Theta_{n_1,n_2}$ has the property asserted by this proposition. \\
\noindent{\bf Proof of Proposition~\ref{propcoup}.}
We begin by describing a set of 
circumstances under which $E_k^{[1]}(n_1)$ and $E_k^{[2]}(n_2)$ are equal under $\Theta_{n_1,n_2}$. For the process $(\tilde{T},\tilde{X})$ used to define  $\Theta_{n_1,n_2}$, we
write $\big\{ \tilde\trar{j}: j \in \N^+ \big\}$ for the sequence of trap entrance arrival times.

Note that, since we are assuming that $n_1 > n_2$,
we have that
$E_k^{[1]}(n_1) = E_k^{[2]}(n_2)$
under $\Theta$, provided that $h_1\big( \reg_{\ell^*}^{\squone} \big) < n_1$, and that
\begin{equation}\label{dphi}
d \Big( \phi(\tilde{T}), \tilde{X} \big( \tilde\trar{n_1 - h_1(r_{\ell^*}^{\squone})} \big)   \Big) > k;
\end{equation}
indeed, in this case, $E_k^{[1]}(n_1)$ and $E_k^{[1]}(n_2)$ each coincide with $\bout_k\big(\tilde{X}_{\tilde\trar{n_1 - h_1(r_{\ell^*}^{\squone})}}\big)$, because this $k$-large exterior neighbourhood lies entirely inside the copy of~$\tilde{T}$.

The process $(\tilde{T},\tilde{X})$ being independent of $\sigma^*$ and having the law of $\pf\big( \cdot \big\vert 0 \in {\rm RG} \big)$, we thus find that
\begin{eqnarray}
 & & \Theta_{n_1,n_2} \Big( E_k^{[1]}(n_1) = E_k^{[2]}(n_2) \Big) \label{cop} \\
 & \geq & \coupl' \Big( h_1\big( \reg_{\ell^*}^{\squone} \big) < n_1 - \frac{n_2}{2} \Big) \nonumber \\
 & & \quad \times \, \,
 \Big( \pf \Big) \Big( d \big( \phi, X(\trar{m}) \big) > k \, \textrm{for all $m \geq
   \frac{n_2}{2}$}  \Big\vert 0 \in {\rm RG} \Big)  . \nonumber
\end{eqnarray}

In light of (\ref{htwort}) with $\ell = \ell^*$, 
\begin{equation}\label{eqhonehtwo}
\textrm{$h_1\big( \reg_{\ell^*}^{\squone} \big) < n_1 - \frac{n_2}{2} \;$  is equivalent to   $\; h_2(\reg_{\ell^*}^{\squtwo}) < \frac{n_2}{2}$.} 
\end{equation}
The quantity
 $h_2(\reg_{\ell^*}^{\squtwo})$, which is non-negative but  may take values as small as zero, may be considered to be a measurement of how long {\rm LEAPFROG} takes to terminate. 
The next lemma bounds the law of this termination time. We defer the proof to the next subsection.
\begin{lemma}\label{lemleapfrog}
There exists a sequence $\big\{ c_\ell: \ell \in \N \big\}$
such that $c_\ell \to 0$ as $\ell \to \infty$ for which
\begin{equation}\label{hoto}
 \Theta' \Big(   h_2\big(r_{\ell^*}^{\squtwo}\big)  \geq \ell \Big)
 \leq c_\ell,
\end{equation}
for all choices $n_1 > n_2$.
\end{lemma}
We must also bound below the second term on the right-hand side of~(\ref{cop}).
Recalling that $\big\{ r_i: i \in \N \big\}$ enumerates the set of regeneration times
of  the walk under $\pf$, note that $\big\{ d \big( \phi, \xsub{r_i}\big): i \in \N \big\}$ is a
strictly increasing sequence, and that $X(j) \geq X(r)$ whenever $j \geq r$
with $r \in {\rm RG}$. Hence, if $\ell \in \N$ is such that $r_k \leq
\trar{\ell}$, then $d \big( \phi,
X(\trar{j}) \big) \geq k + 1$ for all $j \geq \ell$. 
Thus, 
\begin{eqnarray}
  & & \Big( \pf \Big) \Big( d \big( \phi, X(\trar{m}) \big) > k \, \textrm{for all $m \geq
   \frac{n_2}{2}$}  \Big\vert 0 \in {\rm RG} \Big) \nonumber \\  
 & \geq & \Big( \pf \Big) \Big( r_k \leq \trar{\frac{n_2}{2}}  \Big\vert 0 \in {\rm RG}
 \Big). 
\label{folpk}
\end{eqnarray}

To bound below the right-hand side requires some understanding of how trap entrance arrivals are punctuated by regeneration times. In this regard, we state a definition and a lemma. The lemma also contains information to be used in the proof of Lemma~\ref{lemleapfrog}; its proof is deferred to the next subsection.
\begin{definition}\label{defai}
Let $(T,X)$ be sampled from
$\pf$. 
Let $\big\{ r_i: i \geq 1  \big\}$
enumerate the set of regeneration times for $X$, and set  $r_0 = 0$. For $i \geq 0$, set
$$
A_i = \Big\vert  \Big\{ j \in \N : r_i < \trar{j} < r_{i+1}  \Big\} \Big\vert,
$$
where the set $\big\{ \trar{i}: i \in \N \big\}$ is specified in Definition \ref{deftxb}.
\end{definition}
\begin{lemma}\label{lemai}
Assume that $\big\{ p_j: j \in \N \big\}$ satisfies
Hypothesis \ref{hypf}, and that the measure $\nu$ has support in $(1,\infty)$.
The sequence $\big\{ A_j: j \geq 1 \big\}$
is independent and identically distributed. We have that there exists $\conai > 0$ such that, for any $j,\ell \in \N$,
\begin{equation}\label{aikbd}
  \Big( \pf \Big) \Big( A_j \geq \ell \Big) \leq \conai \ell^{-1} \big( \log \ell \big)^{2 - \alpha},
\end{equation}
where the constant $\alpha > 3$ is specified in Hypothesis \ref{hypf}.
Thus,
$$
\E_{\pf} \big( A_1  \big) < \infty.
$$
These statements hold also for $A_0$. Furthermore, $\gcd\big( {\rm
  supp}(A_1) \cap \N^+ \big) = 1$.
\end{lemma}

We make the claim that there exists a constant $C > 0$ such that, for $n,k \in \N$ such that $n > 2k$,
\begin{equation}\label{eqrgclaim}
\Big( \pf \Big) \Big( r_k > \trar{\frac{n}{2}}  \Big\vert 0 \in {\rm RG}
 \Big) \leq C k^2 n^{-1}
 \Big( \log \big( n/(2k) \big) \Big)^{2 - \alpha},
\end{equation}
where the constant $\alpha > 3$ is specified in Hypothesis \ref{hypf}.
To check the claim, note that the condition $r_k > \trar{n/2}$
is equivalent to $\sum_{j=0}^{k-1}A_j \geq n/2$. For the latter condition to be satisfied, one of the summands must exceed $n/(2k)$. A union bound and Lemma \ref{lemai} then yield (\ref{eqrgclaim}), with the choice $C = 2 \conai$.

We learn from (\ref{cop}), (\ref{eqhonehtwo}), 
Lemma \ref{lemleapfrog}, (\ref{folpk}), and (\ref{eqrgclaim}) with applied with $n = n_2$, 
that, for $n_1 \geq n_2 > 2k$,
$$
 \coupl_{n_1,n_2} \Big( E_k^{[1]}(n_1) = E_k^{[2]}(n_2) \Big) \geq  \Big( 1 -  c_{n_2/2}   \Big)
 \bigg( 1 -  C k^2 n_2^{-1}
 \Big( \log \big( n_2/(2k) \big) \Big)^{2 - \alpha}    \bigg).
$$
This completes the proof of Proposition \ref{propqkq}, subject to proving Lemmas~\ref{lemleapfrog} and~\ref{lemai}. \qed

\subsubsection{Proofs of the technical results}\label{subsecfive}

For the proof of Lemma \ref{lemleapfrog}, we need a simple general result concerning  the intersection of 
sets of partial sums of independent sequences.
\begin{lemma}\label{lemaonetwo}
Let $\big\{ A_j^{(1)}: j \geq 1 \big\}$ and $\big\{ A_j^{(2)}: j \geq 1
\big\}$
be two sequences of independent random variables, each term having the same law, satisfying 
$\E A_1^1 < \infty$, ${\rm supp} A_1^{(1)} \subseteq \N^+$, 
and with $\gcd \big( {\rm supp} A_1^{(1)}  \big) = 1$.
Further, let $A_0^{(1)}$ and $A_0^{(2)}$ be independent random
variables of finite mean.
Let
$$
 R_i = \Big\{ \sum_{j=1}^\ell A_j^{(i)}: \ell \geq 0 \Big\}, \qquad i \in \{ 1,2 \}.
$$
Then there exists a sequence $\big\{ c_j: j \in \N \big\}$, $c_j \to 0$,  such that, for
all $\llone,j \in \N$,
$$
\P \Big(  \big( A_0^{(1)} + R_1 \big) \cap \big( \llone + A_0^{(2)} + R_2 \big)
 \cap \big\{ \llone, \ldots, \llone + j \big\} \not= \emptyset \Big)
 \geq 1 - c_j.
$$ 
\end{lemma}
\noindent{\bf Proof of Lemma \ref{lemleapfrog}.}
For each $i \in \{ 1,2 \}$, we let $H_i$ be enumerated as
$\big\{ \sum_{k=0}^\ell A_k^{(i)}  : \ell \geq 0 \big\}$.
Under $\Theta'$, 
the two sequences $\big\{ A_j^{(i)}: j \in \N \big\}$ each have the distribution of $\big\{ A_j: j \in \N \big\}$ under $\pf$.

We apply Lemma \ref{lemaonetwo} 
with the choice of the two independent sequences
 $\big\{ A_i^{(1)} 1\!\!1_{A_i^{(1)} \not= 0} : i \in \N \big\}$ and $\big\{ A_i^{(2)} 1\!\!1_{A_i^{(2)} \not= 0 } : i \in \N
 \big\}$. (Indicators are included because Lemma~\ref{lemaonetwo} is applicable to sequences of strictly positive random variables.) It is Lemma~\ref{lemai} which establishes that the hypotheses of Lemma~\ref{lemaonetwo} are indeed satisfied.  In this way, we learn that
$$
 \Theta' \bigg(  \inf \Big( \big( H_1 - (n_1 - n_2 ) \big)  \cap H_2  \Big) \geq \ell  \bigg) \leq c_\ell  
$$
for some $c_\ell \to 0$. By (\ref{eqhtwo}) and $t_2^* = \reg_{\ell^*}^{\squtwo}$, we obtain Lemma~\ref{lemleapfrog}. \qed

We now provide the proofs of Lemmas \ref{lemai} and \ref{lemaonetwo}. For the former, we will use the next result.
\begin{lemma}\label{lemrgexp}
Let $(B,X)$ have the law $\P_{f,\nu,\infty} \times \P_{B(T),\beta}^{\phi}$
of the backbone and the associated biased random walk. 
Let $\big\{ r_i: i \in \N \big\}$ enumerate the set of regeneration times specified in Definition \ref{defpsikt}  for
this walk. Then the increments $\big\{ r_{i+1} - r_i: i \geq 1 \big\}$ are
independent and identically distributed. There exists $c > 0$ such that, for
each $\ell \in \N$, 
$$
 \Big( \P_{f,\nu,\infty} \times \P_{B(T),\beta}^{\phi} \Big) \Big( r_2 - r_1 \geq \ell  \Big)
 \leq \exp \big\{ - c \ell \big\},
$$ 
and the same statement holds for $r_1$.
\end{lemma}
\noindent{\bf Proof.}
Set 
$Z(n) = d \big( \phi, X(n) \big)$ for $n \in \N$.  Write   
$$
\mathcal{R} = \Big\{ n \in \mathbb{N}: \textrm{$Z(m) < Z(n)$ for $m < n$ and $Z(m) > Z(n)$ for $m > n$} \Big\},
$$
and note that $\mathcal{R}$ coincides with the set of regeneration times. 
It suffices, then, for the required bounds on the laws of $r_1$ and $r_2 - r_1$
to prove that, for some $c > 0$, 
and for each $\ell \in \N$,
\begin{equation}\label{inftp}
  \Big( \P_{f,\nu,\infty} \times \P_{B(T),\beta}^{\phi} \Big) 
 \Big( \inf \mathcal{R} > \ell \Big) \leq \exp \{ - c\ell \}, 
\end{equation}
because $\big\{ 0 \in {\rm RG} \big\}$ has positive $\P_{f,\nu,\infty} \times \P_{B(T),\beta}^{\phi}$-probability.

To this end, let $y_i = \inf \big\{ j \in \N : Z(j) = i \big\}$. As a shorthand, we write $\P$
for the probability measure under which the process $Z$ is defined. Note that $y_{i+1} \geq y_i + \ell$ implies that $Z(y_i + \ell) \leq i = Z(y_i)$.  
Thus, for some $c > 0$, for each $i \in \N$
and $\ell \in \N$,
\begin{eqnarray}
 & & \P \Big( y_{i+1} - y_i \geq \ell \Big\vert Z(1),\ldots, Z(y_i) \Big) \label{eqyi} \\
 & \leq &  \P \Big( Z(y_i + \ell) - Z(y_i) \leq 0 \Big\vert Z(1),\ldots, Z(y_i) \Big) \leq \exp \{ - c\ell \}, \nonumber
\end{eqnarray}
$\sigma \big\{ Z(1),\ldots, Z(y_i) \big\}$-a.s. The second inequality follows directly after 
noting that $Z$ is a walk on $\N$ whose
increments stochastically dominate those of an independent and identically
distributed sequence of $\{ -1,1 \}$-valued random variables, with
probability exceeding $1/2$ of $1$.

We also have that 
\begin{equation}\label{eqzji}
 \P \Big( Z(j) > i \, \textrm{for all $j > y_i$} \Big\vert Z(1),\ldots,
 Z(y_i) \Big) \geq c,
\end{equation}
$\sigma \big\{ Z(1),\ldots, Z(y_i) \big\}$-a.s.
In light of (\ref{eqyi}) and (\ref{eqzji}), 
$\inf \mathcal{R}$ is stochastically dominated by a sum of independent and
identically distributed random variables, each having an exponentially decaying
tail, the number of summands being an independent geometric random
variable. We have obtained (\ref{inftp}), as required. \qed \\
The next remark will be useful on several occasions.\\
\noindent{\bf Remark.} Let $T$ denote an infinite rooted \wgt 
tree. Let $S$ denote the subtree of $T$ induced by a connected set of vertices that contains $\phi(T)$.
If $S$ has the property that the removal of this set of vertices from $T$ results in a graph all of whose components are finite, 
then note that the biased walk on $S$ may be obtained from that on $T$ by restricting the domain of the latter. 
That is, the measure
$\P_{S,\beta}^{\phi}$ is obtained from $\P_{T,\beta}^{\phi}$
by mapping $X: \N \to V(T)$ to $Y:\N \to V(S)$ via 
$Y(n) = X(s_n)$, where  $s_0 = 0$ and
$$
s_n = \inf \Big\{ t > s_{n-1}: X(t) \in V(S), X(t) \not=
X(s_{n-1}) \Big\}, \quad \textrm{for each $n \geq 1$}.
$$ 
\noindent{\bf Proof of Lemma \ref{lemai}.}
Let $\big( T, X \big)$ denote a sample of the law $\pf$.
Using the notation introduced in Definition~\ref{defprop}, it is easy to confirm that, for each $j \geq 0$, the set $\big\{ \ell \in \N : r_j < \trar{\ell} < r_{j+1} \big\}$ is 
measurable with respect to the history until time $r_{j+1}$ but is independent of the history until time $r_j$.
Hence, this sequence of sets is independent. For
any $j \geq 1$, the distribution of $A_j$ is that of $A_0$ under $\big(
\pf \big) \big( \cdot \big\vert 0 \in {\rm RG} \big)$.

We couple the laws $\pf$ and $\pfba$ by associating to $X$ the walk $Y:\N \to B(T)$ given by the preceding remark in the case that $S = B(T)$.
 
Write $\big\{ r_j: j \geq 1 \big\}$ for the set of regeneration times
for $X$, and indicate this set for the process $Y$ with a superscript $Y$; further set $r_0^Y = r_0 = 0$. For $j \in \N$, let $B_j$ denote the set of $\trar{\ell}$ satisfying $r_j
< \trar{\ell} < r_{j+1}$ (so that $A_j = \vert B_j \vert$).
For $v \in V(T)$, write $\offspr_v$ for the set of offspring of $v$ in $T$. Note
that $B_0 \subseteq \bigcup_{j=0}^{r_1^Y}\offspr_{Y(j)}$. Thus, for any $C > 0$,
\begin{eqnarray}
  & &  \Big( \pf \Big) \Big( A_0 \geq \ell \Big) \label{pfby} \\
 & \leq & \Big( \pf \Big) \Big(  \Big\vert \bigcup_{j=0}^{\lfloor C \log \ell
   \rfloor} \offspr_{Y(j)} \Big\vert \geq \ell \Big)  +   \Big( \pf \Big) \Big( r_1^Y
  \geq C \log \ell \Big) \nonumber \\
 & \leq & \sum_{j=0}^{\lfloor C \log \ell \rfloor} \Big( \pf \Big) \Big(
 \big\vert \offspr_{Y(j)} \big\vert \geq \frac{\ell}{C \log \ell}, Y(j) \not\in \big\{
 Y(0),\ldots, Y(j-1)\big\} \Big) \,    +\,  \ell^{-cC}, \nonumber
\end{eqnarray}
where, in the second inequality, we used Lemma \ref{lemrgexp} in regard to the regeneration times of the backbone walk $Y$. 
We claim that, for each $j \in \N$ and $\ell \in \N$, 
\begin{equation}\label{eqclaim}
 \Big( \pf \Big) \Big(
 \big\vert \offspr_{Y(j)} \big\vert \geq \ell  \Big\vert Y(j) \not\in \big\{
 Y(0),\ldots, Y(j-1)\big\} \Big) \leq (1-\pext)^{-1} \sum_{m \geq \ell} p_m,
\end{equation}
where recall that $\pext = \P_f \big( \vert V(T) \vert < \infty \big)$.
To confirm this, note that it is easy to see that, for each $j \in \N$, under the law $\pf$ conditioned by $Y$ visiting a new backbone vertex at time $j$, the conditional distribution of 
the descendent tree $T_{Y(j)}$ is $\P_{f,\infty}$, so that 
$\vert \offspr_{Y(j)} \vert$ has the conditional law of
 $\vert \offspr_\phi \vert$ under $\P_{f,\infty}$; then (\ref{eqclaim}) follows because its right-hand side is an upper bound on $\P_{f,\infty} \big( \vert \offspr_\phi \vert \geq \ell \big)$.
 
Returning then to (\ref{pfby}) and applying Hypothesis \ref{hypf}, we
find that
$$
  \Big( \pf \Big) \Big( A_0 \geq \ell \Big)
  \leq 2 (1-\pext)^{-1}  C^2 \hatc \ell^{-1} \big( \log \ell \big)^{2 - \alpha} + \ell^{-cC},
$$
whence (\ref{aikbd}) for $j = 0$. For $j \geq 1$, (\ref{aikbd})
follows from $\big( \pf \big) \big( 0 \in {\rm RG} \big) > 0$.

The final statement of the lemma follows from $1 \in {\rm supp}(A_1)$, an almost trivial fact whose proof we omit.
 \qed\\
\noindent{\bf Proof of Lemma \ref{lemaonetwo}.}
We will find a non-increasing sequence $\big\{ d_j: j \in \N \big\}$
satisfying $d_j \to 0$ such that, for all $L \geq 0$,
\begin{equation}\label{prl}
\P \Big( R_1 \cap \big( L + R_2 \big) \cap \big\{ L,\ldots, L + j \big\} \not=
\emptyset \Big) \geq 1 - d_j.
\end{equation}
Firstly, we will see why this suffices.

For given $c_1,c_2 \in \N$ satisfying $0 \leq c_1 \leq c_2 \leq j/2$, 
\begin{eqnarray}
 & & \P \Big(  \big( A_0^{(1)} + R_1 \big) \cap \big( \llone + A_0^{(2)} + R_2 \big)
 \cap \big\{ \llone, \ldots, \llone + j \big\} \not= \emptyset \Big\vert
  A_0^{(1)} = c_1, A_0^{(2)} = c_2 \Big) \nonumber \\
 & = &   \P \Big(   R_1  \cap \big( \llone + c_2 - c_1 + R_2 \big)
 \cap \big\{ \llone - c_1, \ldots, \llone - c_1 + j \big\} \not= \emptyset  \Big)
 \nonumber \\
 & \geq & 1 - d_{j - c_2} \geq 1 - d_{j/2}, \nonumber
\end{eqnarray}
the first inequality by (\ref{prl}) with the choice $L = \llone + c_2 - c_1$.
For given $c_1,c_2 \in \N$ satisfying $0 \leq c_2 < c_1 \leq j/2$,
\begin{eqnarray}
 & & \P \Big(  \big( A_0^{(1)} + R_1 \big) \cap \big( \llone + A_0^{(2)} + R_2 \big)
 \cap \big\{ \llone, \ldots, \llone + j \big\} \not= \emptyset \Big\vert
  A_0^{(1)} = c_1, A_0^{(2)} = c_2 \Big) \nonumber \\
 & \geq & \P   \Big(   R_1  \cap \big( \llone + c_2 - c_1 + R_2 \big)
 \cap \big\{ \llone  + c_2 -  c_1, \ldots, \llone -c_1 + j \big\} \not=
 \emptyset  \Big), \label{ronel}
\end{eqnarray}
which, in the case that $\llone + c_2 - c_1 \geq 0$, is, by (\ref{prl}), at least $1 - d_{j -
  c_2}$, which in turn is at least $1 - d_{j/2}$. In the case that  $\llone + c_2 - c_1 < 0$, the term
on the right-hand-side of (\ref{ronel}) is at least
$$
 \P   \Big(  \big(  R_1 - (\llone + c_2 - c_1) \big)  \cap R_2 
 \cap \big\{ - \big( \llone + c_2 - c_1 \big) , \ldots, j - c_2 \big\} \not=
 \emptyset  \Big),
$$
which, by means of (\ref{prl}) with the roles of $R_1$ and $R_2$ interchanged, exceeds $1 - d_{j/2}$, because $(j - c_2) + \big( \llone + c_2 - c_1
\big) \geq j/2$, due to $c_1 \leq j/2$ and $\llone \geq 0$.

Hence,
\begin{eqnarray}
 & & \P \Big(  \big( A_0^{(1)} + R_1 \big) \cap \big( \llone + A_0^{(2)} + R_2 \big)
 \cap \big\{ \llone, \ldots, \llone + j \big\} \not= \emptyset \Big) \nonumber \\
 & \geq & 1 - d_{j/2} - \P \Big( \Big\{ A_0^{(1)} > j/2 \Big\} \cup  \Big\{
 A_0^{(2)} > j/2 \Big\} \Big). \nonumber
\end{eqnarray}
We see that the lemma holds with $c_j = d_{j/2} +  \P \big( \big\{ A_0^{(1)} > j/2 \big\} \cup  \big\{
 A_0^{(2)} > j/2 \big\} \big)$.

As for the proof of (\ref{prl}), it suffices to find a sequence
$\big\{ d_j: j \in \N \big\}$
such that for all $j \in \N$, some $K \in \N$ and all $\llone \in \N$,
\begin{equation}\label{krort}
\P \Big(  \big( K + R_1 \big) \cap \big( K + \llone  + R_2 \big)
 \cap \big\{ K + \llone, \ldots, K + \llone + j \big\} \not= \emptyset \Big)
 \geq 1 - d_j.
\end{equation}

Consider two further independent random subsets $R_1^*$ and $R_2^*$ of $\N$
each having the law of $R_1$. Note that, conditionally on $K \in R_1^*$ and
$K + \llone \in R_2^*$, we may define a coupling 
$K + R_1 = R_1^* \big\vert_{\{ K, \ldots \}}$
and 
$K + \llone + R_2 = R_2^* \big\vert_{\{ K + \llone, \ldots \}}$.

That is, the probability in (\ref{krort}) may be written
\begin{equation}\label{rsrs}
\P \Big(   R_1^*  \cap R_2^*
 \cap \big\{ K + \llone, \ldots, K + \llone + j \big\} \not= \emptyset \Big\vert
 K \in R_1^* , K + \llone \in R_2^* \Big).
\end{equation}
Note that
$$
\P \Big( m \in R_1^* \Big) \to \frac{1}{\E A_1^1},
$$
as $m \to \infty$, by the renewal theorem (\cite{feller}, page 360).
By choosing $K \in \N$ sufficiently high, and invoking independence,
\begin{equation}\label{pkrone}
\P \Big(   
 K \in R_1^* , K + \llone \in R_2^* \Big) \geq  \frac{1}{2 \big( \E A_1^1 \big)^2}.
\end{equation}
The set $R_1^* \cap R_2^*$ is the range of the set of partial sums of a
sequence $\big\{ Y_i: i \in \N \big\} $
of independent and identically distributed random variables, 
with $\E Y_1 = \big( \E A_1^1 \big)^2$, this following from
$$
\P \Big( m \in R_1^* \cap R_2^* \Big) \to \frac{1}{\big( \E A_1^1 \big)^2}
$$ 
as $m \to \infty$, and the renewal theorem, once more.

By (\cite{feller}, (4.10), page 370), then, $\inf R_1^* \cap R_2^* \cap \{ m,\ldots  \} - m$
converges in distribution as $m \to \infty$ (to a law with finite mean).
Hence, 
for $\epsilon > 0$, for any $\llone \geq 0$, $K \geq 0$ and for $j$ sufficiently
high,
\begin{equation}\label{pkrtwo}
\P \Big(   
 R_1^*  \cap R_2^* \cap \big\{ K + \llone,\ldots , K + \llone + j \big\} =
 \emptyset  \Big) \leq  \epsilon.
\end{equation}
By (\ref{pkrone}) and (\ref{pkrtwo}), the probability (\ref{rsrs}) is at least 
$1 - 2 \big( \E A_1^1 \big)^2 \epsilon$,
for $K \in \N$ chosen so that (\ref{pkrone}) holds for each $\llone \geq 0$,
and 
 for sufficiently high $j$. In this way, we obtain (\ref{krort}),  and complete the proof. \qed
\subsubsection{Proof of Corollary \ref{corqkq}}\label{subsecsix}
In the proof of Lemma~\ref{lemleapfrog}, take $A_0^1 = A_0^2 = 0$.  \qed 
\end{subsection}
\begin{subsection}{Deriving Proposition \ref{proprho}}\label{secrho}
We begin by restating Proposition \ref{proprho} in terms of convergence of marginal distributions defined by dividing backbone-tree pairs using the renewal decomposition of the trap that was specified in Definition \ref{defndiv}.
To do this, we extend the notation $\Psi_k$ of Definition \ref{defpsik} to include the interior of a trap. The new notation is illustrated in the upcoming Figure \ref{figxiinpsiplus}.
\begin{definition}\label{defpsikplus}
Given a backbone-tree pair $(B,\singent)$,  $\singent = (\head,\entrance) \circ T_{\entrance}$, that satisfies $r(T_{\entrance})
\geq k+1$, we write $N_k$, the $k$-large neighbourhood of $\entrance \in V(B,\singent)$, for the \wgt tree induced by the set of
vertices of $(\ik \cup \ok)(B,\singent)$.

Write $\Psi_k^+$ \hfff{psikpl} for the set of values of $N_k$ as  $(B,\singent)$ ranges over all such back-bone tree pairs. 
Note that if $\xi \in \Psi_k^+$, then $\ik(\xi)$ has a renewal decomposition containing $k$ components, all of which, including the last, are root-base trees. Recording these components in the form $C_1(\xi),\ldots,C_k(\xi)$, note that $\ik(\xi)$ is itself a root-base tree; its base $\basedef(C_k(\xi))$ we will refer to by $\basedefarg{\xi}$. 
Note also that each $\xi \in \Psi_k^+$ has a vertex labelled $\entrance$, which is $\phi(C_1(\xi))$.

For $K \in \N$, set $\Psi_k^+[K] = \big\{ \xi \in \Psi_k^+: \big\vert V\big(\ik(\xi)\big) \big\vert \leq K \big\}$. \hfff{psikplk}

For any measure $\zeta$ on
backbone-tree pairs supported on those $(B,\singent)$ for which
$r(T_{\entrance}) \geq k + 1$,
we define $\zeta_k$ to be the marginal of $\zeta$ on $\Psi_k^+$.
\end{definition}
We will adopt the shorthand that, for $u > 0$, $\qstaru$  \hfff{qstaru} denotes 
the law of $\qstarnoarg$ conditioned on $\corrfac \omen(T_{\entrance}) > u$.

Using Definition \ref{defpsikplus}, we write $\qsubuksupstar = \big( \qstaru \big)_k$, for each $k \in \N$. 

Proposition \ref{proprho} was stated for the marginal distributions $\rho_{[k]}$ to permit its statement to be as simple as possible; to prove it, we will use the system $\rho_{k}$ instead.
\begin{prop}\label{propqtlim}
Let $k \in \N$.
There exists a measure $\rho_k$ on $\Psi_k^+$ such that 
\begin{equation}
 {\rm TV} \Big( \qsubuksupstar , \rho_k \Big) \to 0
\end{equation}
as $u \to \infty$.
\end{prop}
\noindent{\bf Remark.} 
We remark that the notation $\qsubuksupstar$ is
technically illegitimate, because a sample $(B,\singent)$ of $\qstaru$ may
have $r(T_{\entrance}) \leq  k$, so that $\qsubuksupstar$ is not supported on
$\Psi_k^+$.
However, this notational abuse is permissible for our
purpose, since  $\lim_{u \to \infty} \qstaru \big( r(T_{\entrance}) \leq k \big)
= 0$. This assertion is merely (\ref{qstfku}) after $\halfcorrfack$ has been replaced by $\halfcorrfac$; the modified form has a verbatim proof to the original.   \\
\noindent{\bf Proof of Proposition \ref{proprho}.}
The measure $\rho$ may be constructed from 
the consistent collection $\big\{
\rho_k: k \in \N \big\}$ in Proposition \ref{propqtlim}.
This establishes the statement of Proposition \ref{proprho} with  $\rho_{k}$
replacing $\rho_{[k]}$. Clearly this is equivalent to the actual statement.  \qed

The strategy for proving Proposition \ref{propqtlim} is as follows. Recall the locally defined approximation 
$\halfcorrfack$ to the correction factor $\halfcorrfac$ that was introduced in Lemma \ref{lembt}.
Firstly, we prove an analogue of  Proposition \ref{propqtlim}
where the measure is $\qsubuksupstar$ is replaced by  $\qstarnoarg \big( \cdot \big\vert \corrfack \omen (T_{\entrance}) > u \big)_k$, that is, where the appearance of the correction factor $\halfcorrfac$ in the conditioning is replaced by its approximation $\halfcorrfack$. Secondly, we derive Proposition \ref{propqtlim} by comparison with the approximating version. The next two subsections perform these two steps.
\subsubsection{Understanding 
$\qstarnoarg \big( \cdot \big\vert \corrfack \omen (T_{\entrance}) > u \big)_k$
 as $u \to \infty$}
The next lemma accomplishes the first step. For its statement, note that $\halfcorrfack$ is determined by the data $N_k$, so that $\halfcorrfack(\xi)$ is well defined for each $\xi \in \Psi_k^+$.  Recall from Section \ref{secharris} the subcritical Galton-Watson law $\big\{ h_i:i \in \N \big\}$, and note that Definition \ref{defomega} is used to define the relative weight $\omen \big( \basedefarg{\xi} \big)$.  
\begin{lemma}\label{fllim}
$\empty$ 
\begin{enumerate}
 \item 
For each $k \in \N$, and $\xi \in \Psi_k^+$, we have that
\begin{eqnarray}
 & & \lim_{u \to \infty}
\frac{d \qstarnoarg \big( \cdot \big\vert \corrfack \omen(T_{\entrance}) > u
  \big)_k}{d \qstarnoarg_k}
\big( \xi \big) \label{fora} \\
& = & \frac{ \Big( \corrfack(\xi) \omen \big( \basedefarg{\xi} \big) \Big)^{\chiexp}
  \P_{h,\nu}\Big( r(T) > k \Big)^{-1}}{\int_{\Psi_k^+} \Big( \corrfack(\psi) \omen \big( \basedefarg{\psi} \big) \Big)^{\chiexp} d \alphasub{k}(\psi)}, \label{forb}
\end{eqnarray}
where $\alphasub{k}$ is the measure on $\Psi_k^+$ given by
$$
\alphasub{k} = \qstarnoarg \Big( \cdot \Big\vert r(T_{\entrance}) > k  \Big)_k.
$$
 For each $K \in \N$, the convergence is uniform in $\xi \in \Psi_k^+[K]$. 
\item The limit
$$
\sigmasub{k} : = 
\lim_{u \to \infty}  \qstarnoarg \Big( \cdot \Big\vert \corrfack \omen (T_{\entrance}) > u \Big)_k
$$
exists in total variation, with the measure so defined having support in $\Psi_k^+$,
and with $\frac{d \sigmasub{k}}{d \qstarnoarg_k} \big( \xi \big)$
being given by (\ref{forb}) for each $\xi \in \Psi_k^+$.
\end{enumerate}
\end{lemma}
The ample notation in the statement of Lemma \ref{fllim} may lend the result a foreboding aspect. 
The proof is, however, split into three further lemmas with slightly simpler statements.
Central to the argument is a formula for the Radon-Nikodym derivative in question:
\begin{lemma}\label{lemform}
For $\xi \in \Psi_k^+$,
\begin{eqnarray}
& & \frac{d \qstarnoarg \big( \cdot \big\vert \corrfack \omen(T_{\entrance}) > u
  \big)_k}{d \qstarnoarg_k}
\big( \xi \big) \nonumber \\
& = & \frac{\qstarnoarg \Big( \corrfack \omen\big( T_{\entrance} \big) > u \Big\vert
   \Big\{ N_k = \xi
   \Big\} \cap \Big\{ r\big(T_{\entrance} \big) > k \Big\} \Big)}{\qstarnoarg \Big( \corrfack \omen\big( T_{\entrance} \big) > u \Big)}.  \label{quone}
\end{eqnarray}
\end{lemma}

Recall that $\nu$ is the compactly supported measure from which are drawn the random biases on edges under the law $\P_{f,\nu,\infty}$.
In the case that $\nu$ is discrete, Lemma~\ref{fllim} follows directly from the definition of conditional probability. We defer the general proof of Lemma~\ref{lemform} until the other elements needed 
for the proof of Lemma~\ref{fllim} have been gathered together.

\begin{figure}
\begin{center}
\includegraphics[width=0.6\textwidth]{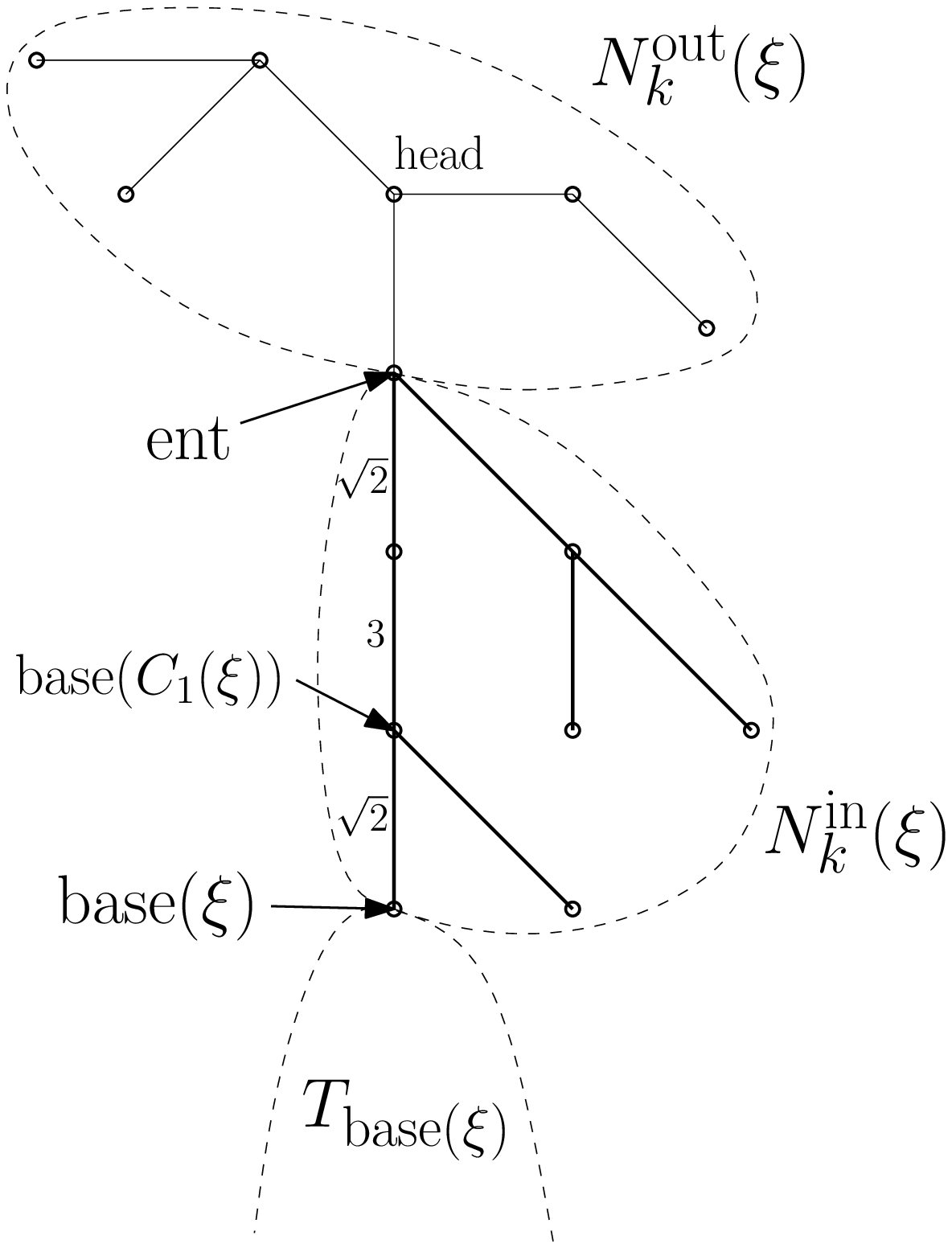} \\
\end{center}
\caption{An instance of $\xi \in \Psi_k^+$; $k=2$ and the edge-bias law $\nu$ has support $\big\{ \sqrt{2},3 \big\}$.
The bias of edges is shown only for those on the path $P_{\entrance,\basedef(\xi)}$. Note that $\omen\big( \basedef(\xi) \big) = \sqrt{2} \cdot 3 \cdot \sqrt{2}$. The figure illustrates the motivation for working with the conditioning $r(T_{\entrance}) > k$ in Lemma \ref{fllim}: conditioning also on the specific form $N_k = \xi$, 
the trap below $\ik(\xi)$, which could (with other definitions) have several components, 
in fact takes the form of the single non-trivial descendent tree $T_{{\rm base}(\xi)}$.}\label{figxiinpsiplus}
\end{figure}
The next two lemmas analyse the numerator and the denominator on the right-hand side of (\ref{quone}).
\begin{lemma}\label{lemintstep}
For each $k \in \N$ and $\xi \in \Psi_k^+$,
\begin{eqnarray}
 & & \qstarnoarg \Big( \corrfack \omen\big( T_{\entrance} \big) > u \Big\vert  \Big\{ N_k = \xi
   \Big\} \cap \Big\{ r\big(T_{\entrance} \big) > k \Big\}  \Big) \nonumber \\
 & \sim &  u^{-\chiexp}
 \big( 1 - h_0 \big)^{-1}  \cpo  \Big(  \corrfack(\xi) 
     \omen \big(  \basedefarg{\xi} \big)   \Big)^\chiexp \label{qutwonew}
\end{eqnarray}
as $u \to \infty$. 
For each $k \in \N$, the implied convergence is uniform in $\xi \in \Psi_k^+[K]$ for any $K \in \N$.
\end{lemma}
\noindent{\bf Proof.}
For a backbone-tree pair $(B,\singent)$ for which $r(T_\entrance) > k$, and with $N_k = \xi$, note that we have the representation 
$$
  \omen(T_{\entrance})  =  \omen \big( \ik(\xi) \big) \, + \,
   \omen \big(  \basedefarg{\xi}  \big) \Big( \omega_{\basedefarg{\xi}} \big( T_{\basedefarg{\xi}} \big) - 1 \Big) 
$$
where recall that 
$\omen \big( \ik(\xi) \big) 
= \sum_{v \in \cup_{i=1}^k V \big( C_i(\xi) \big)} \omen(v)$.
For such a backbone-tree pair, 
we set $\varepsilon(\xi) =  \omen \big( \ik(\xi) \big) - \omen \big( \b( C_k(\xi) ) \big)$, and note that the condition
$\corrfack \omen(T_{\entrance}) > u$  
may be rewritten in the form
\begin{equation}\label{ineqcorr}
\omega_{\basedefarg{\xi}}\big( T_{\basedefarg{\xi}}  \big) > \frac{\frac{u}{\corrfack(\xi)} - \varepsilon(\xi) }{  \omen \big(  \basedefarg{\xi}  \big)  }.
\end{equation}
This rewriting is useful because it expresses the requirement 
$\corrfack \omen(T_{\entrance}) > u$  purely in terms of a condition on the descendent tree 
$T_{\basedefarg{\xi}}$ which hangs off the conditioned region $N_k = \xi$. 

Regarding the conditional law of this descendent tree, 
we claim that
\begin{equation}\label{eqclrd}
T_{\basedefarg{\xi}} \, \, \textrm{under $\qstarnoarg \big( \cdot \big\vert  \big\{ N_k = \xi
   \big\} \cap \big\{ r\big(T_{\entrance} \big) > k \big\}  \big)$}
 \, \, \textrm{has the law $\P_{h,\nu}\big( \cdot \big\vert  \vert E(T) \vert \geq 1 \big)$}.
\end{equation}
This claim depends on the renewal decomposition in a fundamental way: it is here that we exploit this decomposition to find that the portion of $T_\entrance$ below $\basedefarg{\xi}$ is unencumbered by the conditioning we impose higher up this tree (and in the backbone). To derive it, recall that $T_\entrance$
under $\qstarnoarg$ is $\P_{h,\nu}$-distributed. Hence, (\ref{eqclrd}) is implied by Lemma \ref{lemdecom}. 

That 
$\corrfack \omen(T_{\entrance}) > u$   is equivalent to (\ref{ineqcorr}), and (\ref{eqclrd}), yield the equality in
\begin{eqnarray}
 & & \qstarnoarg \Big( \corrfack \omen\big( T_{\entrance} \big) > u \Big\vert  \Big\{ N_k = \xi
   \Big\} \cap \Big\{ r\big(T_{\entrance} \big) > k \Big\}  \Big) \nonumber \\
& = & \P_{h,\nu}
  \Bigg(  \omega(T) >   \frac{\frac{u}{\corrfack(\xi)} -  \varepsilon(\xi) }{
     \omen \big(  \basedefarg{\xi}  \big) } \Bigg\vert   \vert E(T) \vert \geq 1  \Bigg). \label{eqbdamin} \\
  & \sim & 
 \big( 1 - h_0 \big)^{-1}\P_{h,\nu} \Bigg(  \omega(T) >   \frac{u}{\corrfack(\xi) 
     \omen \big(  \basedefarg{\xi}  \big)  }   \Bigg). \label{eqbda}
\end{eqnarray}
The asymptotic equality (asserted as $u \to \infty$) 
is due to the tail regularity stated in Proposition \ref{theoremone}; the implied convergence is uniform over $\xi \in \Psi_k^+[K]$, because, for such $\xi$, $\varepsilon(\xi) \leq  \omen \big( \ik(\xi) \big) \leq K \bsup^K$. Applying Proposition \ref{theoremone} to (\ref{eqbda}) completes the proof, since the convergence is again uniform, due to $\corrfack(\xi)  \omen \big(  \basedefarg{\xi}  \big) \leq C   \bsup^K$ for  $\xi \in \Psi_k^+[K]$. \qed
\begin{lemma}\label{lemintsteptwo}
Let $k \in \N$.
As $u \to \infty$,
\begin{eqnarray}
 & & \qstarnoarg \Big( \corrfack \omen\big( T_{\entrance} \big) > u \Big) \nonumber \\
 & \sim &
 u^{-\chiexp}  \P_{h,\nu} \big( r(T) > k \big) (1 - h_0)^{-1} \cpo \int_{\Psi_k^+}
  \Big(
\corrfack(\psi) \omen \big(  \basedefarg{\psi}  \big) \Big)^{\chiexp} d \alphasub{k}(\psi). \nonumber
\end{eqnarray}
\end{lemma}
\noindent{\bf Proof.}
In essence, the proof works by expressing  
$\qstarnoarg \big( \corrfack \omen\big( T_{\entrance} \big) > u \big)$ as an integral over possible values of $\ik$, and then applying Lemma \ref{lemintstep} to find a high-$u$ asymptotic expression for the conditioned probability. Some care is needed because Lemma \ref{lemintstep} must be applied simultaneously over these events; we may only assert this over choices of $\ik$ belonging to $\Psi_k^+[K]$ for some large $K$. The non-trivial Lemma \ref{lemdivexp}
is the tool which tells us that the contribution of $\Psi_k^+ \setminus \Psi_k^+[K]$ to 
$\qstarnoarg \big( \corrfack \omen\big( T_{\entrance} \big) > u \big)$ is indeed small.

Some notation is convenient as we turn to formulate this precisely. 
Suppose given two sequences of functions $A_K,B_K:(0,\infty) \to (0,\infty)$ parameterized by $K \in \N$; the dependence on $K$ will sometimes be trivial.  
We write $A_K \simk B_K$ as $u \to \infty$ to mean that, for each $\epsilon > 0$, there exists $K_0 \in \N$ such that, for $K \geq K_0$ and for $u \geq u_0(K)$, $\big\vert A_K(u)/B_K(u) - 1 \big\vert \leq \epsilon$. 

We now argue that we have the following concentration on $\Psi_k^+[K]$: for each $k \in \N$,
\begin{eqnarray}
& & \Q \Big(    \corrfack \omen (T_{\entrance}) > u,  r(T_\entrance) > k  \Big) \label{eqlocal} \\
& \simk &
\Q \Big(   \big\vert \ik \big\vert \leq K ,  \corrfack \omen (T_{\entrance}) > u , r(T_\entrance) > k \Big)
\nonumber
\end{eqnarray}
as $u \to \infty$. To derive this, note the following.
In light of (\ref{eqarr}), which defines $\ik$, the derivation (\ref{qstcom}) with $r(T_\entrance)$ replaced by $\ik$
and with $r(T)$ replaced by $\cup_{i=1}^k V(C_i)$ shows that
\begin{eqnarray}
& & \Q \Big(  \big\vert \ik \big\vert > K \Big\vert   \corrfack \omen (T_{\entrance}) > u  \Big) \label{eqik} \\
& \sim & 2 \big(C/c\big)^\chiexp \P_{h,\nu} \Big( \big\vert \cup_{i=1}^k V(C_i) \big\vert > K 
  \Big\vert \omega(T) > u/C \Big). \nonumber
\end{eqnarray}
as $u \to \infty$. However, for some $c \in (0,1)$,
\begin{eqnarray}
& & \P_{h,\nu} \Big( \big\vert \cup_{i=1}^k V(C_i) \big\vert > K 
  \Big\vert \omega(T) > u/C \Big) \label{eqvci} \\
& \leq & \sum_{i=1}^k \P_{h,\nu} \Big( \big\vert  V(C_i) \big\vert > K/k 
  \Big\vert \omega(T) > u/C \Big) \leq k c^{K/k}. \nonumber
\end{eqnarray}
Lemma \ref{lemdivexp} is used for the latter inequality in the last display.
Note then that (\ref{qstfku}), (\ref{eqik}) and (\ref{eqvci}) imply (\ref{eqlocal}).

Dividing (\ref{eqlocal}) by $\Q\big( r(T_\entrance) > k \big)$ and expressing the right-hand side as an integral over values of $\ik$, we find that, for any given $k \in \N$,
\begin{eqnarray}
& & \qstarnoarg \Big(  \corrfack \omen (T_{\entrance}) > u  \Big\vert r(T_{\entrance}) > k
 \Big)  \label{eqqsco} \\
 & \simk &  \int_{\Psi_k^+[K]} \qstarnoarg \Big( \corrfack \omen (T_{\entrance}) > u
 \Big\vert \Big\{ N_k = \psi \Big\} \cap \Big\{ r(T_{\entrance}) > k \Big\} \Big) d
\alphasub{k}(\psi) \nonumber
\end{eqnarray}
as $u \to \infty$.
 By Lemma \ref{lemintstep}, including the asserted uniformity of convergence, we have that,  for any given $k \in \N$,
\begin{eqnarray}
& &   \int_{\Psi_k^+[K]} \qstarnoarg \Big( \corrfack \omen (T_{\entrance}) > u
 \Big\vert \Big\{ N_k = \psi \Big\} \cap \Big\{ r(T_{\entrance}) > k \Big\} \Big) d
\alphasub{k}(\psi)  \nonumber \\
 & \sim & u^{-\chiexp}  \cpo \big( 1 - h_0 \big)^{-1} \int_{\Psi_k^+[K]} 
 \Big( \corrfack(\psi) 
     \omen \big(  \basedefarg{\psi}  \big)     \Big)^\chiexp  d
\alphasub{k}(\psi) \label{eqcork}
\end{eqnarray}
as $u \to \infty$. 
Note that (\ref{eqqsco}) and (\ref{eqcork}) imply that, for any given $k \in \N$,
\begin{eqnarray}
& & \qstarnoarg \Big(  \corrfack \omen (T_{\entrance}) > u  \Big\vert r(T_{\entrance}) > k
 \Big)  \label{eqcofkearly} \\
 & \simk & u^{-\chiexp}  \big( 1 - h_0 \big)^{-1} \cpo \int_{\Psi_k^+[K]} 
 \Big( \corrfack(\psi) 
     \omen \big(  \basedefarg{\psi}  \big)     \Big)^\chiexp  d
\alphasub{k}(\psi)  \label{eqcofk}
\end{eqnarray}
as $u \to \infty$. Set $a_{K}^{(k)}$ equal to the integral expression in (\ref{eqcofk}); write $a_\infty^{(k)}$ for the same integral over $\Psi_k^+$. That the expression in (\ref{eqcofkearly}) is independent of $K \in \N$ implies that, for each $k \in \N$, and for all $\epsilon > 0$, there exists $K_0 \in \N$ such that, if $K_1,K_2 \geq K_0$, then
 $\Big\vert \frac{a_{K_1}^{(k)}}{a_{K_2}^{(k)}} - 1 \Big\vert < \epsilon$. Of course, for each $k \in \N$,
$a_{K}^{(k)} \uparrow a_\infty^{(k)}$. Hence, for any $k \in \N$, 
\begin{equation}\label{eqakk}
u^{-\chiexp}  a_{K}^{(k)} \simk  u^{-\chiexp}  a_\infty^{(k)}
\end{equation}
as $u \to \infty$. We learn that,
for any given $k \in \N$,
\begin{eqnarray}
& & \qstarnoarg \Big(  \corrfack \omen (T_{\entrance}) > u  \Big\vert r(T_{\entrance}) > k
 \Big)  \label{eqclat} \\
 & \sim & 
 u^{-\chiexp}  \big( 1 - h_0 \big)^{-1} \cpo \int_{\Psi_k^+} 
 \Big( \corrfack(\psi) 
     \omen \big(  \basedefarg{\psi}  \big)     \Big)^\chiexp  d
\alphasub{k}(\psi) \nonumber
\end{eqnarray}
as $u \to \infty$;
indeed,  the relation $\simk$ being transitive, (\ref{eqcofk}) and (\ref{eqakk}) imply (\ref{eqclat}) with $\sim$ replaced by $\simk$; but since the functions on each side do not depend on $K$, we obtain (\ref{eqclat}).
The statement of the lemma follows by multiplying (\ref{eqclat}) by $\Q \big(  r(T_{\entrance}) > k \big)$ and then applying (\ref{qstfku}). \qed \\
\noindent{\bf Proof of Lemma \ref{lemform}.}
The numerator on the right-hand side of (\ref{quone}) may involve conditioning on an event of zero $\qstarnoarg$-probability. To define this object formally, and to derive Lemma \ref{lemform}, we invoke the concept of regular conditional probability; see Section 4.1(c) of \cite{durrett}.
Let $\big(\Omega,\mathcal{F},P\big)$ be a probability space, $X:\big(\Omega,\mathcal{F}\big) \to \big( S,\mathcal{S}\big)$ be a measurable map, and $\mathcal{G}$ a sub $\sigma$-field of $\mathcal{F}$. A map $\mu:\Omega \times \mathcal{S} \to [0,1]$ is said to be a regular conditional distribution for $X$ given $\mathcal{G}$ if:
\begin{itemize}
 \item for each $A$, $\omega \to \mu\big( \omega, A \big)$ is a version of $P\big( X \in A \big\vert \mathcal{G} \big)$;
 \item for almost every $\omega$, $A \to \mu \big( \omega , A \big)$ is a probability measure on $\big(S,\mathcal{S} \big)$.
\end{itemize}

In our case, 
set $\Omega = \big\{ \big( B,\singent \big): r(T_\entrance) > k \big\}$. (We choose $\mathcal{F}$ to be a $\sigma$-field rich enough to specify the geometry and edge-bias data of such backbone-tree pairs.)
We take $P = \qstarnoarg \big( \cdot \big\vert r(T_\entrance > k) \big)$, 
and $\mathcal{G} = \sigma \big\{ N_k(B,\singent): (B,\singent) \in \Omega \big\}$; set $S = [0,\infty)$ and let $\mathcal{S}$ be the Borel $\sigma$-field on $S$.  For $(B,\singent) \in \Omega$, set $X(B,\singent) = 2 \halfcorrfack \omen \big( T_\entrance \big)$. We may now explicitly construct the regular conditional distribution by setting
\begin{equation}\label{eqmuform}
\mu \big( (B,\singent), A \big) =  \P_{h,\nu}
  \Bigg(  \omega(T) \in   \frac{\frac{A}{\corrfack(\xi)} -  \varepsilon(\xi) }{
     \omen \big(  \basedefarg{\xi}  \big) } \Bigg\vert   \vert E(T) \vert \geq 1  \Bigg)
\end{equation}
for $(B,\singent) \in \Omega$ satisfying $N_k(B,\singent) = \xi$, and for $A \in \mathcal{S}$; here, we use the notation $b A + c = \big\{ ba + c: a \in A \big\}$ for $A \subseteq [0,\infty)$. This definition  codifies the meaning of $\qstarnoarg \big( 2 \halfcorrfack \omen(T_\entrance) \in \cdot \, \big\vert  \big\{ N_k = \xi
   \big\} \cap \big\{ r\big(T_{\entrance} \big) > k \big\}  \big)$ in terms of $\P_{h,\nu}$ positive probability events, and is inspired by the equality (\ref{eqbdamin}). (Indeed, the present discussion gives precise meaning to this equality.) The second of the defining properties of a regular conditional distribution is trivially satisfied by $\mu$. For the first, we must argue that, for any $G \in \mathcal{G}$ and $A \in \mathcal{S}$,
\begin{equation}\label{eqcint}
\int_G  \mu \Big( \big(B,\singent\big), A  \Big) d  P (B,\singent)
 = P \Big(  N_k \in G, 2 \halfcorrfack \omen(T_\entrance) \in A \Big).
\end{equation}
We derive (\ref{eqcint}) by explaining its left-hand side.  The law $P$ may be formed by firstly sampling $N_k$ according to the appropriate marginal distribution and then appending a $\P_{h,\nu}$-distributed tree (conditioned to be non-trivial) to $\basedef(C_k)$. If $N_k = \xi$, then the quantity $\mu \big( (B,\singent), A  \big)$ is the conditional probability that the sample $(B,\singent)$ of  $P$
satisfies $2 \halfcorrfack \omen(T_\entrance) \in A$. This is due to (\ref{ineqcorr}) and (\ref{eqclrd}). Hence, the left-hand side of (\ref{eqcint}) is equal to the $P$-probability that $N_k \in G$ and that  $2 \halfcorrfack \omen(T_\entrance) \in A$, so that (\ref{eqcint}) indeed holds.

Having formally defined the right-hand side of (\ref{quone}), it remains to justify this formula. For this, it is enough to integrate each side over an arbitrary $G \in \mathcal{G}$ with respect to the law $\qstarnoarg_k$ and check that the  the same result is obtained in each case. For the left-hand side, the integral is given by 
$\qstarnoarg \big( G \big\vert   2 \halfcorrfack \omen \big( T_\entrance \big) > u \big)$ by definition. For the right-hand side, the outcome is $\frac{\qstarnoarg \big( 2 \halfcorrfack \omen \big( T_\entrance \big) > u, N_k \in G  \big)}{\qstarnoarg \big(  2 \halfcorrfack \omen \big( T_\entrance \big)}$, because the integral in question may be interpreted by the decomposition in the preceding paragraph. The two values are equal by the law of conditional probability, completing the proof of Lemma \ref{lemform}. {\flushright \qed} \\  
\noindent{\bf Proof of Lemma \ref{fllim}.} \\
\noindent{\bf (1).} The statement is obtained by applying Lemma \ref{lemintstep} to the denominator and Lemma \ref{lemintsteptwo} to the numerator of (\ref{quone}).  \\
\noindent{\bf (2).} Note that (\ref{eqik}) and (\ref{eqvci}) imply that for some $c \in (0,1)$, for each $k,K \in \N$ and for all $u > 0$,
$$
 \Q \Big(  \big\vert \ik \big\vert > K \Big\vert   \corrfack \omen (T_{\entrance}) > u  \Big)
 \leq 2 \big(C/c \big)^\chiexp k c^{K/k}.
$$
Hence, for any $\epsilon > 0$, there exists $K \in \N$ such that $\sigmasub{k} \big( \big\vert \ik \big\vert > K \big) < \epsilon$. The uniform convergence of the Radon-Nikodym derivative over $\xi \in \Psi_k^+[K]$ then yields Lemma \ref{fllim}(2). \qed
\subsubsection{Handling the error from the approximate correction factor}
We now prove Proposition \ref{propqtlim} by means of Lemma \ref{fllim}.
\begin{lemma}\label{lemlong}
There exlsts $c > 0$ such that, for each $\ell \in \N$, and for $u$ sufficiently high, 
$$
 1 \leq  \frac{\qstarnoarg \Big( \corrfac \omen(T_{\entrance}) > u \Big)}{\qstarnoarg \Big(
   \corrfacarg{\ell} \omen(T) > u \Big)} \leq
 1 + \exp \big\{ - c\ell  \big\}. 
$$
\end{lemma}
\noindent{\bf Proof.}
The first inequality follows from $\halfcorrfacarg{\ell} \leq \halfcorrfac$. Regarding the second, note that
(\ref{fkrel}) and $\halfcorrfac \leq c$ imply that there exists $c > 0$ such that, for each $\ell \in \N$,
$\halfcorrfac \leq \halfcorrfacarg{\ell} \big( 1 + \exp \big\{ - c \ell \big\} \big)$. 
Hence,
$$
 \qstarnoarg
 \Big(  \corrfac \omen(T_{\entrance}) > u  \Big) 
\leq
  \qstarnoarg
 \Big(  \corrfacarg{\ell} \omen(T_{\entrance}) > u\big(1 + \exp\{ - c \ell \}\big)^{-1}  \Big). 
$$
Lemma \ref{lemintsteptwo} shows that the right-hand side has a pure power-law decay in $u$. Thus, the second inequality in the statement is obtained with an adjustment to the value of $c > 0$. \qed \\
\noindent{\bf Proof of Proposition \ref{propqtlim}.} 
Let $A$ be an event defined on the space of backbone-tree pairs that satisfies $A \in
\sigma \big\{ N_k \big\}$ for some given $k \in \N$.
We claim that there exists $c > 0$ such that, for each $\ell \in \N$ and for all high enough $u$,
\begin{eqnarray}
 & &   \qstarnoarg \Big( A \cap \Big\{  \corrfacarg{\ell} \omen(T_{\entrance}) > u \Big\} \Big) 
  \leq  \qstarnoarg \Big( A \cap \Big\{  \corrfac \omen(T_{\entrance}) > u \Big\} \Big) \label{longeqn} \\
 & \leq & \qstarnoarg \Big( A \cap \Big\{  \corrfacarg{\ell} \omen(T_{\entrance}) > u \Big\} \Big) 
 \, + \, \qstarnoarg \Big(  \corrfacarg{\ell} \omen(T_{\entrance}) > u \Big) \exp \big\{ - c\ell \big\}. 
  \nonumber
\end{eqnarray}
The first inequality is due to $\halfcorrfacarg{\ell} \leq \halfcorrfac$; the second follows from
\begin{eqnarray}
  & & \qstarnoarg \Big( \Big\{ \corrfac \omen(T_{\entrance}) > u \Big\} \Delta  \Big\{  \corrfacarg{\ell} \omen(T_{\entrance}) > u \Big\} \Big) \nonumber \\
 & \leq & \exp \big\{ - c \ell \big\} \qstarnoarg \Big(    \corrfacarg{\ell} \omen(T_{\entrance}) > u  \Big), \nonumber
\end{eqnarray}
which is implied by Lemma \ref{lemlong}.

Using 
\begin{equation}\label{cieqn}
  \qstaru \big( A  \big) 
 = \frac{ \qstarnoarg \Big( A \cap \Big\{  \corrfac \omen(T_{\entrance}) > u \Big\}
   \Big)}{ \qstarnoarg \Big(  \corrfac \omen(T_{\entrance}) > u \Big)},
\end{equation}
we find that, for $u$ sufficiently high,
\begin{eqnarray}
 & & \frac{\qstarnoarg \Big( A \cap \Big\{  \corrfacarg{\ell} \omen(T_{\entrance}) > u \Big\}
   \Big)}{\qstarnoarg \Big( \corrfacarg{\ell} \omen(T_{\entrance}) > u  \Big) \Big( 1 + \exp \big\{ - c\ell \big\} \Big)}  \leq \qstaru \big( A \big) \nonumber \\
  & \leq & \frac{\qstarnoarg \Big( A \cap \Big\{  \corrfacarg{\ell} \omen(T_{\entrance}) > u
    \Big\} \Big) \, + \,\qstarnoarg \Big( \corrfacarg{\ell} \omen(T_{\entrance}) > u  \Big) \exp
    \big\{ - c\ell \big\} }{\qstarnoarg \Big( \corrfacarg{\ell} \omen(T_{\entrance}) > u  \Big)}, \nonumber 
\end{eqnarray}
where the numerator and the denominator in (\ref{cieqn}) were bounded, both above and below, 
respectively by (\ref{longeqn}) and Lemma \ref{lemlong}. 

That is,
\begin{eqnarray}
 & & \frac{\qstarnoarg \Big( A \Big\vert \corrfacarg{\ell} \omen(T_{\entrance}) > u  \Big)}{1 + \exp \big\{ - c\ell \big\}} \leq  \qstaru \big( A \big) \nonumber \\
 & \leq & \qstarnoarg \Big( A \Big\vert  \corrfacarg{\ell} \omen(T_{\entrance}) > u  \Big) \, + \, \exp \big\{ - c\ell \big\} . \label{twoineq}
\end{eqnarray}
Recall from Lemma \ref{fllim}(2) the total variation limit
$\sigmasub{\ell}$. This result implies that, by making a small decrease in the value of $c> 0$, we obtain the following uniform control.
For  $\ell \in \N$ sufficiently high,
there exists $u_{\ell} \in (0,\infty)$
such that, for any $A \in \sigma \big\{ \ikun_{k},\okun_{k} \big\}$, and for all $u \geq u_\ell$,
$$
\frac{\sigmasub{\ell}(A) - \exp \big\{ - c\ell \big\}}{1 + \exp \big\{ - c\ell \big\}}
\leq \qstaru \big( A \big)
\leq  \sigmasub{\ell}(A) + \exp \big\{ - c\ell \big\}.
$$
It follows directly from these bounds that, for all $u,u' \geq u_{\ell}$ and for any $A \in \sigma \big\{ \ikun_{k},\okun_{k} \big\}$,
$\big\vert \qstaru \big( A \big) - \qstarreal^{u'} \big( A \big)  \big\vert \leq 3 \exp \big\{ -c\ell \big\}$. The integer $\ell$ being arbitrary, we see that $\qstaru \big( A \big)$  converges as $u \to \infty$ uniformly over  $A \in \sigma \big\{ \ikun_{k},\okun_{k} \big\}$. 
 Noting that, by definition, 
$\qsubuksupstar(A) = \qstaru(A)$, we confirm the existence of the total variation limit asserted in 
 Proposition~\ref{propqtlim}.
 This completes the proof of Proposition~\ref{propqtlim}. \qed  
\end{subsection}
\begin{subsection}{Proofs for Section \ref{subsectotaltime}}\label{secthmfin}
Here, we present the proofs of Lemmas \ref{lemsmallerr}, \ref{lemvrg},  \ref{lemrhox} and \ref{lemrexist} and Theorem~\ref{proptraplaw}.

We begin by proving the first two of these results.
 We do so in three steps. Firstly, we establish Lemma \ref{lemsmallerr} except for the error bound (\ref{eqsmallerr}); then we prove Lemma \ref{lemvrg}; and then we prove (\ref{eqsmallerr}).
\begin{subsubsection}{Proofs of Lemmas \ref{lemsmallerr} and \ref{lemvrg}}
\noindent{\bf Construction of $E$ and $\ero$.}
Let $\trwtbd \in (0,\infty)$ be given by Proposition \ref{proptsp}, and fix a choice of $\barecon > 0$ (subject to a condition to be specified later).
Let $(B,\singent)$ denote a given backbone-tree pair in the support of $\qstarnoarg$.
If $\omen\big( T_\entrance \big) > \trwtbd$ and $T_\entrance$ is $\barecon$-bare, then we may construct under $\P_{(B,\singent),\beta}^\entrance$
random variables $E$ and $\ero$ by means of Proposition \ref{proptsp}. If $(B,\singent)$ violates one or other of these conditions, we simply choose under $\P_{(B,\singent),\beta}^\entrance$ the random variable $E$ to be an independent exponential random variable of mean $\corrfac \omen\big( T_\entrance \big)$ and then set 
$\ero = \tau_{T_\entrance} - E 1\!\!1_\deepfall$.
This construction evidently satisfies the properties stated in Lemma \ref{lemsmallerr}, subject to verifying 
(\ref{eqsmallerr}).
 
\vspace{3mm}

\noindent{\bf Proof of Lemma \ref{lemvrg}.}
Recall that, under $\btmeas$ conditional on its environment $(B,\singent)$, $E$ has the exponential distribution of mean 
$\corrfac \omen\big( T_\entrance \big)$.  Writing $\mu$ for the law of $\corrfac
\omen(T_{\entrance})$ under $\qstarnoarg$ , we have then that
\begin{eqnarray}
   \btmeasbrac \big( E > x \big)
 & = & \int_0^{\infty} \mu \big( x/v, \infty \big) \exp \{ - v \} dv
 \nonumber \\
 & \sim & \cpo \cpt x^{-\chiexp}  \int_0^{\infty} v^{\chiexp}  
  \exp \{ - v  \} dv, 
\nonumber 
\end{eqnarray}
where the asymptotic equality follows directly from Proposition \ref{trprop}. 
\qed

\vspace{3mm}

\noindent{\bf Proof of Lemma \ref{lemsmallerr}: derivation of (\ref{eqsmallerr}).}
The construction of $E$ and $\ero$ just made entails that, if a backbone-tree pair $(B,\singent)$ is such that $T_\entrance$ is of high enough weight and is $B_0$-bare, then $\E_{(B,\singent),\beta} \big( \big\vert \ero \big\vert \big) \leq C_2 \omega \big( T_\entrance \big)^{1/2}$. In any case, we have the weaker bound $\E_{(B,\singent),\beta} \big( \big\vert \ero \big\vert \big) \leq C \omen \big( T_\entrance \big)$, by $\halfcorrfac \leq C$ (see Lemma \ref{lemebfbd}) and Proposition \ref{proptsp}. By Markov's inequality, we obtain
\begin{equation}\label{eqmarkovone}
  \Big( \btmeas \Big) \Big( \ero > x \Big\vert \omen\big(T_\entrance \big) = y, \, \textrm{$T_\entrance$ is $B_0$-bare} \Big) \leq C \frac{y^{1/2}}{x}
\end{equation}
and
\begin{equation}\label{eqmarkovtwo}
  \Big( \btmeas \Big) \Big( \ero > x \Big\vert \omen\big(T_\entrance \big) = y, \, \textrm{$T_\entrance$ is not $B_0$-bare} \Big) \leq C \frac{y}{x}.
\end{equation}
We now bound the tail of the error $\ero$:
$$
 \qstarnoarg \big( \ero > x \big) \leq A_1 + A_2 +   A_3 + A_4, 
$$
where
$$
A_1 = \qstarnoarg \Big( \omen\big(T_\entrance \big) > \frac{x}{\log x}, \, \, \textrm{$T_\entrance$ is not $B_0$-bare}   \Big),
$$
$$
A_2 = \btmeasbracbig \Big( \omen\big(T_\entrance \big) \leq \frac{x}{\log x}, \, \, \textrm{$T_\entrance$ is not $B_0$-bare}, \ero > x   \Big),
$$
$$
A_3 = \qstarnoarg \Big( \omen\big(T_\entrance \big) > x \log x  \Big),
$$
and
$$
A_4 = \btmeasbracbig \Big( \omen\big(T_\entrance \big) \leq  x \log x, \, \, \textrm{$T_\entrance$ is $B_0$-bare}, \ero > x    \Big).
$$
By fixing $\barecon > 0$  so that $\barecon \geq 2 \connew^{-1} \chiexp$, where the constant $\connew > 0$ appears in Lemma \ref{lemcreg}, we find from this lemma and  Proposition \ref{theoremone} that
\begin{eqnarray}
A_1 & \leq & 2 \big( \log x \big)^{-\connew \barecon} \qstarnoarg \Big(  \omen\big(T_\entrance \big) > \frac{x}{\log x} \Big)
 \nonumber \\
 & \leq & C x^{-\chiexp} \big( \log x \big)^{\chiexp - \connew \barecon} \leq  C x^{-\chiexp} \big( \log x \big)^{- \connew \barecon/2}.
 \nonumber
\end{eqnarray}
We write $\mulaw$ for the law of $\omega(T)$ under $\P_{h,\nu}$, which coincides with the law of $\omen\big(T_\entrance\big)$ under $\qstarnoarg$.  We have that
\begin{eqnarray}
A_2 & \leq & \int_0^{\frac{x}{\log x}} \big( \qstarnoarg \times  \P_{(B,\singent),\beta}^\entrance \big) \Big(  \ero > x \Big\vert \omen\big(T_\entrance \big) = y, \, \textrm{$T_\entrance$ is not $B_0$-bare} \Big) d \mulaw(y) \nonumber \\
 & \leq & C x^{-1} \int_0^{\frac{x}{\log x}} y d \mulaw(y) \leq C x^{-\chiexp} \big( \log x \big)^{\chiexp - 1},
\nonumber
\end{eqnarray}
the first inequality by (\ref{eqmarkovtwo}) and the second a straightforward inference from Proposition \ref{theoremone}.  
We also find that
$$
A_3 \leq C x^{-\chiexp} \big( \log x \big)^{-\chiexp},
$$
by Proposition \ref{theoremone}. Finally, by (\ref{eqmarkovone}), we have that
\begin{eqnarray}
A_4 & \leq & \int_0^{x \log x} \big( \qstarnoarg \times \P_{(B,\singent),\beta}^\entrance \big) \Big(  \ero > x \Big\vert \omen\big(T_\entrance \big) = y, \, \textrm{$T_\entrance$ is $B_0$-bare} \Big) d \mulaw(y) \nonumber \\
 & \leq & C x^{-1} \int_0^{x \log x} y^{1/2} d \mulaw(y) \leq C x^{-1/2 - \chiexp} \big( \log x \big)^{1/2 - \chiexp}.
 \nonumber
\end{eqnarray}
Recalling that $\chiexp < 1$ for the case of $A_2$, we compare these bounds with the decay given in Lemma \ref{lemvrg} and thereby complete the proof of (\ref{eqsmallerr}). \qed 
\end{subsubsection}
\begin{subsubsection}{Proof of Lemma \ref{lemrhox}}
Let $\zetax$ denote the environment marginal of  
$\btmeasbrac \big( \cdot \big\vert E = x \big)$.
It suffices to show that,
for each $k \in \N$, 
\begin{equation}\label{sxlim}
{\rm TV} \Big( \big( \zetax \big)_k, \rho_k  \Big) \to 0
\end{equation}
as $x \to \infty$, where the subscript $k$ indicates the marginal distributions introduced in Definition \ref{defpsikplus}.  To prove (\ref{sxlim}), we require:
\begin{lemma}\label{lemb}
Let $r > 0$. Write $\rur$ for the 
measure $\qstarnoarg$, 
conditional on $\corrfac \omen \big( T_\entrance \big) \in \big[ u,u(1+r) \big]$.
Then, for each $r > 0$,
$$
{\rm TV} \Big( \big( \rur \big)_k , \rho_k \Big) \to 0
$$
as $u \to \infty$. 
\end{lemma}
\noindent{\bf Proof.} 
As in the proof of Proposition \ref{proprho},
for $u > 0$, $\qstaru$ denotes 
the law of $\qstarnoarg$ conditioned on $\corrfac \omen(T_{\entrance}) > u$.
Note that
\begin{equation}\label{abav}
\qstaru 
 = a_{u,r} \rur + b_{u,r} 
\qstar{u(1+r)},
\end{equation}
where
\begin{displaymath}
b_{u,r} = \qstaru \Big( \corrfac \omen \big( T_\entrance \big) > u \big( 1 + r \big) \Big),
\end{displaymath}
$a_{u,r} = 1 - b_{u,r}$.
The result then follows from Propositions \ref{trprop}
and \ref{proprho}.  \qed 

We now derive (\ref{sxlim}). Write again $\mu$ for the law of $\corrfac \omen \big( T_{\entrance} \big)$
under $\qstarnoarg$, and recall that $E$ under $\btmeas$ 
is the product of a $\mu$-distributed random variable and an independent mean one exponential random variable.
Let $\mu^x$ denote the law of 
$\corrfac \omen \big( T_{\entrance} \big)$, under $\btmeas$
conditionally on $E=x$. 

Set $h_x: [0,\infty) \to [0,\infty)$, $h_x(y) = I_x^{-1} y^{-1} \exp \big\{ - x/y \big\}$,
with 
$$
I_x = \int_0^{\infty} u^{-1} \exp \big\{ - x/u \big\} d \mu(u). 
$$ 
Then we have that, for $y \geq 0$,
$$
 \frac{d \mu^x(y)}{d \mu(y)} = h_x(y).
$$
Let $\epsilon > 0$.
By Proposition \ref{trprop}, there exist constants $0 < c < C < \infty$ such that, for
all $x$ sufficiently high,
\begin{equation}\label{xyexp}
\int_0^{cx} h_x(y) d \mu(y) +
\int_{Cx}^{\infty} h_x(y) d \mu(y) < \epsilon. 
\end{equation}
For $i \in \N$, set $P_{i,n}^x = \big[ xi/n, x(i+1)/n \big]$
and $r_{i,n}^x = h_x(xi/n)$. Further set $q_{\epsilon,x}^{(n)}:[0,\infty)
\to [0,\infty)$,
$$
q_{\epsilon,x}^{(n)}(y) = A_{n,x}^{-1} \sum_{i = \lfloor nc \rfloor}^{ \lfloor nC
  \rfloor + 1} r_{i,n}^x 1\!\!1_{P_{i,n}^x}(y),
$$
where $A_{n,x}$ is a normalization chosen to ensure that 
$\int_0^{\infty} q_{\epsilon,x}^{(n)}(y) d \mu(y) = 1$.
It is readily verified that, for $\epsilon > 0$, fixing $n$ high enough,
and for all sufficiently large $x$, 
$\big\vert h_x(y) - q_{\epsilon,x}^{(n)}(y)  
\big\vert  < \epsilon$ for all $y \in (cx,Cx)$. By this and (\ref{xyexp}), and
by choosing $n$ high enough, we obtain, for all sufficiently high $x$,
\begin{equation}\label{yix}
\Big\vert \Big\vert h_x - q_{\epsilon,x}^{(n)}  
\Big\vert \Big\vert_{L^1(d \mu)}  <  \epsilon.
\end{equation}
Set
$$
\msupx : = A_{n,x}^{-1} \sum_{i = \lfloor nc \rfloor}^{ \lfloor nC \rfloor + 1} r_{i,n}^x q_{i,n}^x
 \qstarnoarg \Big( \cdot \Big\vert \corrfac \omen \big( T_\entrance \big) \in P_{i,n}^x \Big),
$$
with $q_{i,n}^x = \qstarnoarg \big( \corrfac \omen ( T_\entrance ) \in P_{i,n}^x \big)$.
We learn from (\ref{yix}) that 
$$
 {\rm TV} \Big(  \zetax, \msupx \Big) \leq \epsilon.
$$
Noting that 
$P_{i,n}^x = xin^{-1} \big[ 1, 1 + r \big]$, with $r = i^{-1}$, 
we see that
Lemma \ref{lemb} yields (\ref{sxlim}). This completes the proof of Lemma \ref{lemrhox}. \qed
\end{subsubsection}
\begin{subsubsection}{Deriving Lemma \ref{lemrexist}}
\begin{lemma}\label{lemaak}
Let $k \in \N$, and let $(B,\singent)$ be a backbone-tree pair
such that the number of components $r\big(T_\entrance \big)$ in the renewal decomposition of $T_\entrance$ is at least $k+1$. (The decomposition was introduced in Subsection \ref{secren}.)
Recalling the notation of Lemma \ref{lembt}, 
define under $\P_{(B,\singent),\beta}^\entrance$ the event
$$
\deepfall_k = \Big\{ H_{\basedef(C_k)} < H_{V(\ok) \setminus V(\okmo)} \Big\}.
$$
In the case that  $r\big(T_\entrance \big) \leq k$, we set $\deepfall_k = \big\{ H_{\vbase} < H_{V(\ok) \setminus V(\okmo)} \big\}$.

There exists $c \in (0,1)$
such that, for any backbone-tree pair $(B,\singent)$ and for each $k \in \N$,
\begin{equation}\label{eqaak}
 \P_{(B,\singent),\beta}^\entrance \big( \deepfall \Delta \deepfall_k \big)  \leq c^k.
\end{equation}
\end{lemma}
\noindent{\bf Proof.} 
Note that $\deepfall \Delta \deepfall_k$ entails that the walk under $\P_{(B,\singent),\beta}^\entrance$ either fails to reach $\vbase$ after arriving at  $\b\big( C_k \big)$, or that it returns to $\entrance$ after reaching  a backbone vertex at distance $k+1$ from $\entrance$. Now, the probability of the first of these alternatives is at most $\binf^{-k}$,
since $d \big( \entrance  , \b( C_k ) \big) \geq k$; and the probability of the second is given by $\ekm - \pesc $, which is at most $2 \binf^{1 - k/2}$, by (\ref{edata}).\qed \\
Recall that $\zetax$ denotes the environment marginal of $\btmeasbrac \big( \cdot \big\vert E = x \big)$.
\begin{lemma}\label{lemenvmar}
Under each of the laws 
$\btmeasbrac \big( \cdot \big\vert E = x \big)$ and
$\zetax \times \P_{(B,\singent),\beta}^\entrance$, the walk trajectory $X:\{0,\ldots,H_\vbase \} \to V(B,\singent)$ has the same distribution.
\end{lemma}
\noindent{\bf Proof.} By the definition of $\zetax$, the two laws have the same environment marginal.
Lemma \ref{lemsmallerr} states that, under $\btmeasbrac \big( \cdot \big\vert E = x \big)$ given $(B,\singent)$, the law of the walk until hitting $\vbase$ (possibly at infinite time) is unaffected by dispensing with the conditioning $E=x$. The latter walk law is simply  $\P_{(B,\singent),\beta}^\entrance$ stopped at time $H_\vbase$. \qed \\
\noindent{\bf Proof of Lemma \ref{lemrexist}.} 
We begin by arguing that
\begin{equation}\label{eqqak}
\lim_{x \to \infty} \btmeasbracbig \Big( \deepfall_k  \Big\vert E = x\Big) = \big( \rho \times \P_{(B,\overline{T}),\beta}^\entrance \big)  \big( \deepfall_k \big).
\end{equation}

Note that, given any $(B,\singent)$, under $\P_{(B,\singent),\beta}^\entrance$, 
\begin{equation}\label{eqsenend}
\deepfall_k \in \sigma \big\{ X(0),\ldots, X(H_\vbase) \big\}.
\end{equation}

Hence, Lemma \ref{lemenvmar} implies that
\begin{equation}\label{eqdpex}
\btmeasbrac \big( \deepfall_k  \big\vert E = x \big) = \big( \zetax \times \P_{(B,\singent),\beta}^\entrance \big) \big( \deepfall_k \big).
\end{equation}
Adopting the usage of Definition \ref{defpsikplus} and
abbreviating $\zeta^x_k = \big( \zetax \big)_k$, note that the right-hand probability in (\ref{eqdpex}) is given by
$\big( \zeta^x_k \times \P_{(B,\singent),\beta}^\entrance \big) \big( \deepfall_k \big)$; this 
is  because $\zeta^x_k$ has support in $\Psi_k^+$ and $\deepfall_k \in \sigma \big\{ \ik,\ok \big\}$.

We see from
 Lemma \ref{lemrhox} that  the law $\rho_k$ (on $\{\ik,\ok\}$) is the total variation limit as $x \to \infty$ of 
 $\zeta^x_k$ (for convergence of all $\rho_{k}$ to $\zeta^x_k$ and of all $\rho_{[k]}$ to $(\zeta^{(x)})_{[k]}$ are obviously equivalent).
Hence,
$$
\lim_{x \to \infty}
\big( \zeta^x_k \times \P_{(B,\singent),\beta}^\entrance \big) \big( \deepfall_k \big)
= \big( \rho_k \times \P_{(B,\overline{T}),\beta}^\entrance \big) \big(\deepfall_k\big),
$$
whose right-hand side is equal to $\big( \rho \times \P_{(B,\overline{T}),\beta}^\entrance \big) \big(\deepfall_k\big)$. We thus obtain (\ref{eqqak}).

Recall that the measure $\rho$ is supported on backbone-tree pairs $(B,\singent)$, where $T_\entrance$ is
an infinite tree. By averaging Lemma \ref{lemaak} by $\rho$, we obtain

\begin{equation}\label{eqrhoaka}
  \Big\vert \big( \rho \times \P_{(B,\overline{T}),\beta}^\entrance \big) \big( \deepfall_k \big) -  \big( \rho \times \P_{(B,\overline{T}),\beta}^\entrance \big)  \big( \deepfall \big)  \Big\vert \leq c^k.
\end{equation}

Note that the sentence including (\ref{eqsenend}) holds equally for the event  $\deepfall \Delta \deepfall_k$.
Hence, Lemma \ref{lemenvmar} implies that 
\begin{equation}\label{eqqstmin}
   \btmeasbracbig \Big( \deepfall \Delta \deepfall_k  \Big\vert E = x \Big)   = 
   \big( \zetax \times \P_{(B,\singent),\beta}^\entrance  \big) \big( \deepfall \Delta \deepfall_k  \big)
\end{equation}
Averaging Lemma \ref{lemaak} over $\zetax$ shows that the right-hand side of (\ref{eqqstmin}) is at most $c^k$. Hence,
\begin{equation}\label{eqqst}
   \btmeasbracbig \Big( \deepfall \Delta \deepfall_k  \Big\vert E = x \Big)    \leq c^k.
\end{equation}
By (\ref{eqqak}), (\ref{eqrhoaka}) and (\ref{eqqst}),
we find that the limit $\lim_{x \to \infty} \big( \btmeas \big) \big( \deepfall  \big\vert E = x \big)$
exists and equals $\big( \rho \times \P_{(B,\overline{T}),\beta}^\entrance \big)(\deepfall)$. \qed
\end{subsubsection}
\begin{subsubsection}{The proof of Theorem \ref{proptraplaw}}
Lemmas \ref{lemvrg} and \ref{lemrexist} imply that 
$\big( \btmeas \big) \big( E 1\!\!1_{\deepfall} > x \big) \sim  \constr \cpo \cpt \Gamma\big( 1 + \chiexp \big)  x^{-\chiexp}$. We now apply Lemma \ref{lemlit}(2) with the choice $U = E 1\!\!1_{\deepfall}$ and $V = \ero$.
Lemma \ref{lemsmallerr} tells us that the hypotheses of Lemma \ref{lemlit}(2) are satisfied, and also that the conclusion that we reach is the statement of Theorem \ref{proptraplaw}. \qed    
\end{subsubsection}
\end{subsection}
\end{section}
\begin{section}{Comparison with independent heavy tail}\label{seciht}
The task remaining is to derive Theorem \ref{thm} from Theorem \ref{proptraplaw}.
To obtain Theorem \ref{thm}, we must describe the asymptotic decay of the law of the time $\hitting{n}$ for a walk under $\pf$ to reach distance $n$
from the root. This time may be considered as a sum of holding times in traps encountered along the way, a sum which, we anticipate, will be dominated by a few large terms, corresponding to the deepest of the encountered traps. These deep traps are well separated, so that, after the walk departs from one of them, 
the environment from the next one should resemble the late-time trap environment measure $\qstarnoarg$, conditionally on the trap $T_\entrance$ being deep. This suggests a means of deriving Theorem \ref{thm}: we will seek to couple the law $\pf$ to a sequence of independent copies of $\btmeas$ in such a way that the hitting time $\hitting{n}$ under $\pf$
approximates the total time spent in a certain initial number of traps in the idealized trap sequence.
This coupling made, we may invoke Theorem \ref{proptraplaw} to derive Theorem \ref{thm}, since this result describes the total time spent in one of the traps in the idealized sequence.

We now summarize the proof strategy in more detail, and introduce some more of the required notation.

\vspace{3mm}

\noindent{\bf Step 1: the hitting time $\hitting{n}$ is close to a sum of $\constt n$ trap holding times.}
To implement the plan that we have just sketched, we firstly have to understand how long the initial sequence of traps should be in the comparison made there. That is, how many trap entrances does the walk under $\pf$ typically visit before reaching a distance $n$ from the root?

Recall from Definition \ref{deftxb} that, under $\pf$, the sequence $\big\{ \trar{i}: i \in \N \big\}$
enumerates the successive trap entrances encountered by the walk. 
The descendent tree $T_{X(\trar{i})}$ of the $i$-th trap entrance is then the $i$-th trap so encountered. We may also record the total time $\tau_i$ spent by the walk in that trap: \hfff{ttt}
$$
\tau_i = \Big\vert \Big\{  j \in \N: X(j) \in V \big( T_{X(\trar{i})} \big)  \Big\} \Big\vert.
$$
It is convenient to use in place of $\hitting{n}$ its close cousin $\overh{n}$, here defined:
\begin{definition}
For $X$ having law 
$\P_{T,\beta}^{\phi}$, let \hfff{distrootb}
$\overh{n} =  \inf \big\{ i \in \N: d \big( \phi, X(i) \big) = n, X(i) \in V(B) \big\}$.
\end{definition}
In this first step, we establish two basic properties: firstly that the hitting time $\overh{n}$ under $\pf$
may be approximated by the sum of the holding times $\tau_i$ over those trap entrances $\trar{i}$ reached by the walk before time $\overh{n}$, and, secondly, that the number of these trap entrances typically grows linearly in $n$ at a deterministic rate. Regarding the latter, we will prove:  
\begin{lemma}\label{lemf}
Define the ``trap-number'' function  \hfff{trenen} $\tnb{\empty}:\N \to \N$,
$$
\tnb{n} = \sup \big\{ j \in \N: \trar{j} < \overh{n} \big\}.
$$
Then there exists a constant $\constt  \in (0,\infty)$ such that
$$
\lim_{n \to \infty} \frac{\tnb{n}}{n}   = \constt,  
$$
$\P_{f,\nu,\infty} \times \P_{T,\beta}^{\phi}$-a.s.
\end{lemma}
The conclusion of this first step is then:
\begin{lemma}\label{lemt}
With $\constt \in (0,\infty)$ 
denoting the constant provided by Lemma \ref{lemf}, and $C \in (0,\infty)$
a sufficiently high constant, for any $\epsilon > 0$ and for $n$ sufficiently high,
$$
\Big( \P_{f,\nu,\infty} \times \P_{T,\beta}^{\phi} \Big)
 \Big( \sum_{i=1}^{ \constt n(1-\epsilon)} \tau_i \leq \overh{n} \leq 
  \sum_{i=1}^{ \constt n(1+\epsilon)} \tau_i \, + \, Cn   \Big) \geq 1 - \epsilon.
$$
\end{lemma}
We are discussing a walk in a sub-ballistic regime, so that the linear quantity $Cn$
in Lemma \ref{lemt} should be regarded as a small error term. In the ballistic case, of course, the durations spent in backbone and traps are both linear: we mention that, having carried out Step 1, we will give a short argument (in Lemma \ref{lemdnconv}) that shows that motion is ballistic in the case that $\chiexp > 1$, and that   proves, in essence, the final assertion of Theorem \ref{thm}.

\vspace{3mm}

\noindent{\bf Step 2: coupling the walk with an idealized model of holding times in traps.}
The holding times that form the summands in the approximation $\sum_{i=1}^{ \constt n} \tau_i$, are, we anticipate, close to independent, except if their indices are 
are very close to one another. As such, it is natural to model this sequence of holding times by an independent and identically distributed sequence, with each term being distributed as 
the holding time in a trap encountered at late time by the walk. The law of this idealized trap, viewed from its entrance, 
is the limiting law $\qstarnoarg$ that was specified in Definition \ref{defqstar}. For this reason, we specify the idealized sequence of holding times in the following way.
\begin{definition}\label{defidseq}
Let $\idtrap$ \hfff{idtrap}
denote a law under which
is defined an independent sequence 
$\big( \tilde{B}^i, \tilde{T}^i \big)$, 
$\tilde{T}^i = \big( \head,\entrance \big) \circ \tilde{T}^i_\entrance$, $i \in \N$,
of backbone-tree pairs, each having law $\qstarnoarg$.
Recall that $\tilde{T}^i_\entrance$ refers to the descendent tree of the
vertex $\entrance$ in the backbone-tree pair $\big( \tilde{B}^i
,\tilde{T}^i \big)$. 
For $i \in \N$, we introduce $\tilde{\omega}_i$ under $\idtrap$ by $\tilde\omega_i = \omen \big( \tilde{T}^i_\entrance
\big)$. 
Further, let 
$$
\tilde{X}_i: \N
 \to \big( \tilde{B}^i , \tilde{T}^i \big), \qquad i \in \N,
$$
be independent under $\idtrap$, $\tilde{X}_i$ 
having law $\P_{(\tilde{B}^i,\tilde{T}^i),\beta}^{\entrance}$. \hfff{idtrapnot}
Set the holding time $\tilde{\tau}_i$ according to
$$
\tilde{\tau}_i = \Big\vert \Big\{ j \in \N : \tilde{X}_i(j) \in V \big( \tilde{T}^i_\entrance \big) \Big\} \Big\vert.
$$
\end{definition}
Our aim is to couple $\pf$ and $\idtrap$ so that 
the sum of a long initial sequence of holding times $\tau_i$ 
is close to the counterpart sum under the idealized sequence of traps under $\idtrap$. 
This second step of the proof leads then to the construction of the following coupling. 
\begin{lemma}\label{lemp}
Assume that $\chiexp$, specified in (\ref{defchi}),
satisfies $\chiexp < 1$.
 For any $\epsilon > 0$, there exist $n_0 \in \N$ and, for each $n \geq
n_0$, a coupling $\coup_n$ of the distributions 
$\idtrap$ and $\pf$
such that
$$
\coup_n \Big( \Big\vert \sum_{i=1}^n \tau_i - \sum_{i=1}^n \tilde{\tau}_i
\Big\vert > \epsilon \sum_{i=1}^n \omen \big(  T_{X(\trar{i})} \big)  \Big) < \epsilon,
$$
and such that $\omega \big( T_{X(\trar{i})}\big) = \tilde\omega_i$ for all $i \in \N$ under $\Theta_n$.
\end{lemma}
With these elements in place, the proof of Theorem \ref{thm} is essentially complete. 
Indeed, Lemmas \ref{lemt} and \ref{lemp} provide a coupling 
according to which $\hitting{n}$ under $\pf$ is well approximated by the independent and identically distributed sum 
$\sum_{i=1}^{\constt n} \tilde{\tau}_i$ under $\idtrap$, while Theorem \ref{proptraplaw} describes the regularly decaying tail of $\tilde{\tau}_1$ under  $\idtrap$.

The next two subsections present the proofs of the results stated in the respective steps. A brief subsection then follows in which Theorem \ref{thm} is derived from Proposition \ref{theoremone} by means of Lemmas \ref{lemt} and \ref{lemp}.
\begin{subsection}{Proofs for step 1}
Here, we furnish the proofs of Lemmas \ref{lemf} and \ref{lemt}. \\
\noindent{\bf Proof of Lemma \ref{lemf}.}
Under $\pf$, enumerate ${\rm RG} = \{ r_i: i \in \N \}$ and  set $\kappa = \E \big( d(\phi,X(r_2)) - d(\phi,X(r_1))  \big)$, the mean being with respect to this law.
That $\kappa \in (0,\infty)$ follows directly from Lemma \ref{lemrgexp}.
Renewal blocks being independent, we find
that
$$
 \frac{d \big( \phi  , X(r_n) \big)}{n} \to \kappa,
$$ 
$\P_{f,\nu,\infty} \times \P_{T,\beta}^{\phi}$-a.s., by the
strong law of large numbers. 
 Setting $g_n = \sup \big\{ j \in \N: r_j < \overh{n}\big\}$, note that 
$d \big( \phi, X(r_{g_n}) \big) \leq  n <
 d \big( \phi, X(r_{g_n + 1}) \big)$, whence  
\begin{equation}\label{kapinv}
\frac{g_n}{n} \to \kappa^{-1},
\end{equation}
$\P_{f,\nu,\infty} \times \P_{T,\beta}^\phi$-a.s.

Recalling Definition \ref{defai}, we find that
$$
   A_0 + \sum_{i=1}^{g_n - 1} A_i  \leq \tnb{n} \leq 
 A_0 + \sum_{i=1}^{g_n} A_i.  
$$
Invoking Lemma \ref{lemai}, and (\ref{kapinv}),
$$
\frac{\tnb{n}}{n} \to
\frac{\E_{\P_{f,\nu,\infty} \times \P_{T,\beta}}^{\phi}
  A_1}{\kappa},
$$
$\P_{f,\nu,\infty} \times \P_{T,\beta}^{\phi}$-a.s, by the
strong law of large numbers applied to the independent and identically distributed sequence $\big\{ A_i: i \geq 1 \big\}$. 
Thus, the statement of the lemma holds with 
$\constt  = \kappa^{-1} \E_{\pf} A_1$.
\qed \\
\noindent{\bf Proof of Lemma \ref{lemt}.}
For $(T,X)$ sampled from 
$\P_{f,\nu,\infty} \times \P_{T,\beta}^{\phi}$,
let $J(n)$ denote the event that, if $i,j \in \{ 0,\ldots, n \}$
satisfy $j \geq i + \delta n$, then there exists $r \in {\rm RG}$ such that
$\trar{i} < r < \trar{j}$. (We set $\trar{0} = 0$ here.) Adopting once more the notation introduced in 
Definition \ref{defai}, note that
$$
 \bigcap_{i=0}^{Cn}  \Big\{ A_i \leq \delta n \Big\}
 \, \cap \, \Big\{ r_{Cn} > \trar{n} \Big\}  \subseteq J(n),
$$
for $C$ sufficiently high.
Note that, by (\ref{aikbd})
in the statement of Lemma \ref{lemai}, 
$$
\Big( \pf \Big) \Big( A_i > \delta n \, \, \textrm{for some $i \leq Cn$} \Big) \leq 2 C \conai  \delta^{-1}  \big( \log n \big)^{2 - \alpha}
$$
for $C$ sufficiently high.
Each regeneration epoch independently witnesses a new trap entrance with positive probability, so that
$r_{Cn} \leq \trar{n}$ has a probability under $\pf$ that  decays exponentially in $n$,  for $C$ sufficiently high. Hence, 
\begin{equation}\label{jnineq}
\Big( \P_{f,\nu,\infty} \times \P_{T,\beta}^{\phi} \Big) 
 \Big( J(n) \Big) \geq 1 - C  \big( \log n \big)^{2 - \alpha}.
\end{equation}
Note that, if $\overh{n} \leq \trar{Cn}$
and $J(Cn)$ occurs, then
\begin{equation}\label{unionf}
 \bigcup_{i=1}^{\tnb{n} - \delta n} \Big\{ j \in \N : 
 X(j) \in V \big(  T_{X(\trar{i})} \big) \Big\} \subseteq \big\{ 1,\ldots, \overh{n} \big\}.
\end{equation} 
Indeed, by the occurrence of $J(Cn)$, and $\trar{\tnb{n}} \leq \overh{n} \leq \trar{Cn}$,
there exists $r \in {\rm RG}$ satisfying 
$\trar{\tnb{n} - \delta n} < r < \trar{\tnb{n}}$, and, if 
$i
\in \N$ satisfies $\trar{i} < r$ for a certain $r \in {\rm RG}$, then
$$
\sup \Big\{ j \in \N: X(j) \in V \big(  T_{X(\trar{i})} \big) \Big\} < r.
$$
Thus, (\ref{unionf}).
Provided that $\tnb{n} \geq n \constt  (1-\delta)$, it follows from (\ref{unionf})
that
$\overh{n} \geq \sum_{i=1}^{n \constt (1-2\delta)}\tau_i$. We find that, for any $\delta
> 0$ and $n$ sufficiently high,
\begin{equation}\label{eqdnlb}
 \Big( \P_{f,\nu,\infty} \times \P_{T,\beta}^{\phi} \Big)
 \Big( \overh{n} \geq \sum_{i=1}^{n \constt (1-2\delta)} \tau_i \Big) \geq 1 - \delta,
\end{equation}
by (\ref{jnineq}), 
$$
\Big( \P_{f,\nu,\infty} \times  \P_{T,\beta}^{\phi} \Big) \Big( \overh{n} 
\leq \trar{Cn}
\Big) \geq 1 - \frac{\delta}{2},
$$
and Lemma \ref{lemf}.
The last inequality follows from 
the equivalence of $\overh{n} \geq \trar{Cn}$ and $\tnb{n} \geq Cn$, and
Lemma \ref{lemf}, with the choice $C > \constt$ being made.

Having obtained the required bound on the lower tail of $\overh{n}$ in (\ref{eqdnlb}), we turn to the bound on the upper tail.  Note that
\begin{equation}\label{incdn}
\big\{ 1,\ldots, \overh{n} \big\}
 \subseteq \bigcup_{i=1}^{\tnb{n}} \Big\{ j \in \N : X_j \in V \big(  T_{X(\trar{i})}
 \big) \Big\} \, \cup \, \Big\{ j \in \N: X_j \in V(B), j \leq \overh{n} \Big\}.
\end{equation}
We decompose 
\begin{equation}\label{eqdecom}
\big\{ j \in \N: X_j \in V(B)  \big\} = F_1
\cup F_2,
\end{equation}
the set $F_1$ consisting of those moments 
of visit to the backbone by the walk
 with the property that
the preceding instance of the walk visiting the backbone took place at its present location.

Note that each element of $F_1$ is a moment of return to the backbone
$B(T)$ by the walk $X$. For $j \in \N$, let $h_j \in V(B)$
denote the parent of $X(\trar{j})$. Define 
$$
\Delta_j = \Big\{ i \in \N: X(i) = h_j, X(i-1) = X\big(\trar{j}\big) \Big\}.
$$
Note that the union $\bigcup_{j=1}^{\infty} \Delta_j$ is disjoint and
partitions $F_1$. 

We claim that, for each $j \in \N$,
\begin{equation}\label{edelj}
\E_{\pf} \big( \big\vert \Delta_j \big\vert \big) \leq \frac{\bsup + \binf + 1}{\binf - 1}.
\end{equation}
Indeed, given $T$ and the values $\big\{ X(0),\ldots,X(\trar{j}) \big\}$, the
first return of $X$ to the backbone after time $\trar{j}$ forms the first
element of $\Delta_j$. Given that $\Delta_j$ has at least $k$ elements,
with $k \geq 1$, a $(k+1)$-st occurs only in the case that $X$ revisits
$X(\trar{j})$ after its visit to $h_j$ at the $k$-th element of $\Delta_j$. This
probability is $1 - \pesc$ as specified by Definition \ref{eqfdef},
for the backbone-tree pair for which $\entrance = X(\trar{j})$. The lower bound on $\pesc$
in Lemma \ref{lemebfbd} thus gives (\ref{edelj}).

Note further that
$$
 F_1 \cap \big\{ 1,\ldots, \overh{n} \big\} \subseteq \bigcup_{j=1}^{\tnb{n}} \Delta_j,
$$ 
where $\tnb{n}$ was defined in Lemma \ref{lemf}.
We find that, for any $\epsilon > 0$ and $\delta > 0$,
\begin{eqnarray}
& & \Big( \pf \Big) \Big(  \Big\vert F_1 \cap \big\{ 1,\ldots, \overh{n} \big\}
\Big\vert \geq \delta^{-1} \big( 1 + \bsup \big) \big( \constt   +
\epsilon \big) n \Big) \nonumber \\
& \leq &  \Big( \pf \Big) \Big(  \Big\vert  \bigcup_{j=1}^{(\constt +  \epsilon)n} \Delta_j
\Big\vert \geq \delta^{-1} \big( 1 + \bsup \big) \big( \constt   +
\epsilon \big) n \Big) \nonumber \\
 & & \, \, + \,\,  \Big( \pf \Big) \Big( \tnb{n} \geq \big( \constt   + \epsilon
\big) n \Big) \leq \delta + \epsilon_n, \label{fdfn}
\end{eqnarray}
where (\ref{edelj}) and Markov's inequality were used in the second
inequality, and where $\epsilon_n \to 0$ is due to Lemma \ref{lemf}.

Note that the sequence $\big\{ X(j):
j \in F_2 \big\}$ has the law of the walk $\P_{f,\nu,\infty} \times
\P_{B(T),\beta}^{\phi}$ on the backbone, so that  the sequence
$\big\{ d\big(\phi,X(j)\big):j \in F_2 \big\}$ stochastically dominates a
nearest-neighbour random
walk $Z: \N \to \Z$, $Z(0) = 0$,
with independent and identically distributed increments, with rightward transition probability
$\beta/(\beta + 1) > 1/2$. 
By this, and (\ref{fdfn}), we have that, for $C > 0$ large enough, and for
all $n$ sufficiently high, 
\begin{equation}\label{pfjdn}
 \Big( \P_{f,\nu,\infty} \times  \P_{T,\beta}^{\phi} \Big) 
 \Big(  \Big\vert \Big\{ j \in \N: X(j) \in V(B) , j \leq \overh{n}  \Big\} \Big\vert
 \leq C n \Big) \geq 1 - \delta.
\end{equation}
By Lemma \ref{lemf}, (\ref{incdn}) and (\ref{pfjdn}), we obtain the bound on the upper tail
of $\overh{n}$ asserted by the lemma. Alongside (\ref{eqdnlb}), this  completes the proof. \qed 

We now provide the short argument establishing the final assertion of Theorem \ref{thm}
for $\overh{n}$.  The actual statement, with $\hitting{n}$ in place of $\overh{n}$, will also require Lemma \ref{lemdnhn}. 
\begin{lemma}\label{lemdnconv}
If $\chiexp > 1$, then there exists $v \in (0,\infty)$ such that then 
$n^{-1} \overh{n} \to v$, $\pf$-a.s.
\end{lemma}
\noindent{\bf Proof.}
Recalling that $g_n = \sup \big\{ j \in \N: r_j < \overh{n} \big\}$, note that
$$
 r_1 + \sum_{i=1}^{g_n - 1} \big( r_{i+1} - r_i \big) < \overh{n} \leq r_1 + \sum_{i=1}^{g_n} \big( r_{i+1} - r_i \big).
$$  
The remark after Definition~\ref{defprop} implies that the summands here are independent and identically distributed. 
Recalling from (\ref{kapinv}) that $n^{-1} g_n \to \kappa^{-1}$, $\pf$-a.s., it is suffices to show that 
$\E_{\pf}\big( r_2 - r_1 \big) < \infty$. Note that $r_2 - r_1$ under $\pf$ has the distribution of $r_1$
under $\big( \pf \big) \big( \cdot \big\vert 0 \in {\rm RG} \big)$. By partitioning the time-interval $\big\{ 0,\ldots,r_1\big\}$ according to the walk's location in trap or backbone, note that 
\begin{equation}\label{eqroneone}
r_1 + 1 = \sum_{j=1}^{A_0} \tau_i \, + \, \Big\vert \Big\{ j \in \N: X(j) \in V(B), j \leq r_1 \Big\} \Big\vert.
\end{equation}
It follows from (\ref{eqnmean}) and (\ref{edelj}) that under $\big( \pf \big) \big( \cdot \big\vert 0 \in {\rm RG} \big)$ the mean of $\tau_i$ is at most $2 (\bsup + \binf -1) (\binf - 1)^{-1} \E_{h,\nu}(\omega(T))$ for each $i \in \N$.
However, Proposition \ref{theoremone}
and $\chiexp > 1$ imply that  $\E_{h,\nu}(\omega(T))$ is finite. By Lemma \ref{lemai}, we obtain then that $\sum_{j=1}^{A_0} \tau_i$ has finite mean under $\big( \pf \big) \big( \cdot \big\vert 0 \in {\rm RG} \big)$. 
To treat the second term on the right-hand-side of (\ref{eqroneone}), we use the decomposition (\ref{eqdecom}). Noting that 
$F_1 \cap \big\{ 0,\ldots, r_1 \big\} \subseteq \cup_{j=1}^{\vert A_0 \vert} \Delta_j$, we see that
$\big\vert F_1 \cap \big\{ 0,\ldots, r_1 \big\} \big\vert$ has finite mean under $\big( \pf \big) \big( \cdot \big\vert 0 \in {\rm RG} \big)$ by (\ref{edelj}) and Lemma \ref{lemai}. To establish  the lemma, it remains only to show the same assertion for $F_2$. To do this, note that, under $\pf$, $\big\{ X(j):j \in F_2 \big\}$ is the walk on the backbone given by the walk $Y$ with $S = B(T)$ specified in the remark preceding Lemma \ref{lemai}. Note further that $r_1$ is the first regeneration time of  
 $\big\{ X(j):j \in F_2 \big\}$ considered as a walk on the backbone. Hence,  $\big\vert F_1 \cap \big\{ 0,\ldots, r_1 \big\} \big\vert$, being at most one plus the first regeneration time of a $\pfba$-distributed walk, is shown to have finite mean by Lemma \ref{lemrgexp}. \qed 
\end{subsection}
\begin{subsection}{Proofs for step 2}
We seek to argue that the sum $\sum_{i=1}^n \tau_i$ under $\pf$ is well approximated by the sum of $n$ trap holding times in the idealized independent sequence $\idtrap$. Our strategy for proving this is to argue that the sum 
 $\sum_{i=1}^n \tau_i$ is dominated by sojourns in high-weight traps. These traps are far apart, so that the walk is likely to pass regeneration times as it moves from one to the next. This entails that the holding times in these deep traps are independent. Moreover, plenty of time passes between a regeneration time and the arrival at a new deep trap, with the effect that the backbone environment about the new trap entrance has time to converge towards $\queu$, the limiting law which was introduced in Proposition \ref{propqkq} and which forms the backbone marginal of our idealized trap environment $\qstarnoarg$. 

To implement this sketch, we begin by defining a procedure that 
will identify most of the big traps that the walk
encounters, while ensuring that the times spent in the selected traps form an independent
and identically distributed sequence.
\begin{definition}\label{defjiv}
For $v \in (0,\infty)$ and $k \in \N$,
we construct a subsequence of indices of the sequence $\big\{ \trar{j}: j \in \N
\big\}$ on the space 
$\pf$. Let $(T,X)$ denote a sample of $\pf$. 
Let $r_1^v$ denote the first regeneration time of $X$. 
Let 
$$
 j_1^v = \inf \Big\{   j \in \N: \trar{j} > r_1^v,  \, d \Big( X(\trar{j}) , X(r_1^v)   \Big) > k , \, \omenn{X(\trar{j})} \big(  T_{X(\trar{j})} \big) > v \Big\}.
$$
That is, $\trar{j_1}^v$ is the first time after $r_1^v$ at which $X$ visits a
vertex in $\bext = V(\baug) \setminus V(B)$ at distance greater than $k$ from $X(r_1^v)$ and whose descendent tree has weight exceeding
$v$.

Suppose that a sequence of times 
$$
 r_1^v < \trar{j_1^v} < \ldots < r_n^v < \trar{j_n^v}
$$
has been constructed, in such a way that
$r_i^v \in {\rm RG}$ is a regeneration time for each $i \leq n$.
We now define $r_{n+1}^v$ to be the first regeneration time after time
$\trar{j_n^v}$
such that $d \big( X(\trar{j_n^v}), X(r_{n+1}^v) \big) > k$.
We then set
$$
 j_{n+1}^v = \inf \Big\{   j \in \N: \trar{j} > r_{n+1}^v, \, d \Big( X(\trar{j}) , X(r_{n+1}^v)   \Big) > k , \,  \omenn{X(\trar{j})} \big(  T_{X(\trar{j})} \big) > v \Big\}.
$$
Note that the $k$-dependence is suppressed in denoting the two sequences.
\end{definition}
Inherent in the sequence of times just defined is the property that the $k$-neighbourhoods on the backbone of the 
trap entrances $\trar{j_i^v}$ are independent for distinct $i$. To prove an independence statement for the holding times $\tau_i$ under $\pf$, it is thus natural to approximate these times by cut-off versions, where we consider time only until the walk leaves the $k$-neighbourhood of the trap entrance on the backbone. The following definition introduces this notation, and its counterpart for the idealized traps under the measure $\idtrap$.

\begin{definition}\label{defru}
Under $\pf$, we set, for $k \in \N$ and $i \in \N$,
$$
\sigma_i^k = \inf \Big\{ j > \trar{i}: X(j) \in V(B), d \Big( X(j) , X(\trar{i}) \Big) =
k + 1 \Big\}
$$
and
$$
\tau_i^k = \Big\vert \Big\{  j \in \N: X(j) \in V \big(  T_{X(\trar{i})} \big), j
\leq \sigma_i^k  \Big\} \Big\vert.
$$
Under $\idtrap$, we set, for $k \in \N$ and $i \in \N$, 
$$
\tilde{\sigma}_i^{k} = \inf \Big\{ j \in \N : \tilde{X}_i(j) \in V\big( \tilde{B}^i \big),
 d \Big( \tilde{X}_i(j), \entrance \Big) = k + 1 \Big\}
$$
and
$$
\tilde{\tau}_i^{k} = \Big\vert \Big\{ j \in \big\{ 0,\ldots, \tilde{\sigma}_i^k \big\}: \tilde{X}_i(j) \in V \big( \tilde{T}^i_\entrance \big) \Big\} \Big\vert.
$$
\end{definition}

\begin{lemma}\label{lemjd}
Let $v > 0$ and $k \in \N$. The law $\pf$ has the following properties.
\begin{enumerate}
 \item 
The sequence of traps $\big\{ T \big( X(\trar{j^v_n}) \big): n \in \N \big\}$
 is independent and identically distributed, with common law
 $\P_{h,\nu} \big( \cdot \big\vert \omega(T) > v \big)$. 
  \item  
Recall the notions $\bout_k$ and $\qsubk$ from Definition \ref{defpsik} and Proposition \ref{propqkq}. 
The sequence $\big\{ \bout_k\big( X(\trar{j^v_n}) \big): n \in \N \big\}$
 of $k$-large backbone neighbourhoods of the successive trap entrances is independent and identically distributed.
  Let $\eta_v^k$ denote its common law.
  We have that ${\rm TV} \big( \eta_v^k , \qsubk \big) \to 0$ as $v \to \infty$.
 \item  Let $X'_n = X'_{n,v,k}:\N \to T$ be defined under the law $\pf$ so that, given $T$ and the initial segment 
$X:\{ 0,\ldots, \trar{j^v_n} \} \to T$, it is distributed according to $\P_{T,\beta}^{X(\trar{j^v_n})}$. Then, for any given $k \in \N$, and for any $n \in \N$, the construction of $X'_n$ may be effected so that $X'_n:\N \to T$ coincides with $X\big( \trar{j^v_n} + \cdot \big):\N \to T$
with $\pf$-probability tending to one as $v \to \infty$.
\end{enumerate}
\end{lemma}
\noindent{\bf Proof.}\\
\noindent{\bf (1).} This follows directly from Definition \ref{defjiv}.\\
\noindent{\bf (2).} Independence and identical distribution follows directly from the construction. 
Regarding the second assertion in this part, under the law $\pf$, let $Q_v = \inf \big\{ i \in \N: \omega_{X(\trar{i})}(T_{X(\trar{i}})) > v \big\}$. 
Let
 $\overline\eta_v^k$  
 denote the law of  $\bout_k\big( X(\trar{Q_v}) \big)$ under
$\big( \pf \big) \big( \cdot  \big\vert  0 \in {\rm RG} \big)$. 
It suffices to prove that  
\begin{equation}\label{tvone}
{\rm TV} \big( \overline\eta_v^k , \qsubk \big) \to 0
\end{equation}
 and that 
\begin{equation}\label{tvtwo}
 {\rm TV} \big( \eta_v^k , \overline\eta_v^k \big) \to 0
 \end{equation}
 as $v \to \infty$.

To establish (\ref{tvone}), note firstly that $Q_v$ under $\big( \pf \big) \big( \cdot  \big\vert  0 \in {\rm RG} \big)$ is distributed as a geometric random variable of parameter $p = \P_{h,\nu}\big( \omega(T) > v \big)$, and that it is independent of the walk induced on the augmented backbone (which is the walk $Y$ defined in the remark preceding the proof of Lemma \ref{lemai} when the choice $S = \baug(T)$ is made). Since $p \to 0$ as $v \to \infty$, we have that 
\begin{equation}\label{pfq}
   \Big( \pf \Big)  \Big( Q_v \in [0,C] \Big\vert   0 \in {\rm RG}  \Big) \to 0,
\end{equation}
for any given $C > 0$.

Adopting the notation of Corollary \ref{corqkq}, we see that
 $\overline\eta_v^k$ has the law of $\qsubk(n)$, where the index $n$ is independently randomized with the law of $Q_v$.
Hence, we obtain (\ref{tvone}) from Corollary \ref{corqkq}  and (\ref{pfq}).

To derive (\ref{tvtwo}), note that, under $\big( \pf \big) \big( \cdot \big\vert 0 \in {\rm RG} \big)$,  $Q_v = j_1^v$ unless $d \big( \phi, X(\trar{Q_v}) \big) \leq k$. However, if 
$Q_v = j_1^v$, then naturally $\bout_k\big(X(\trar{Q_v})\big) = \bout_k\big(X(\trar{j_1^v})\big)$, so that
 ${\rm TV} \big( \eta_v^k , \overline\eta_v^k \big) \leq 
 \big( \pf \big) \big(  d \big( \phi, \xsub{\trar{Q_v}} \big) \leq k  \big\vert  0 \in {\rm RG} \big)$. Hence, it suffices for (\ref{tvtwo}) to argue that 
 \begin{equation}\label{eqxuqv}
 \Big( \pf \Big) \Big(  d \big( \phi, \xsub{\trar{Q_v}} \big) \leq k  \Big\vert  0 \in {\rm RG} \Big) \to 0 
\end{equation}
 as $v \to \infty$. 
Noting that the sequence of terms $d(\phi,X(r))$ as $r$ ranges over $\rm RG$ is strictly increasing,
and that, for each $r \in {\rm RG}$, we have that $d(\phi,X(s)) \geq d(\phi,X(r))$ for all $s > r$,
we find from Lemma \ref{lemai} that, for each $\ell \in \N$,
\begin{equation}\label{eqpfrg}
 \Big( \pf \Big) \Big( d \big( \phi  , X(\trar{n}) \big) \geq \ell \, \, \textrm{for all $n \geq m$} \Big\vert 0 \in {\rm RG}  \Big) \to 1,
\end{equation}
as $m \to \infty$.
We obtain (\ref{eqxuqv}), and thus (\ref{tvtwo}), from (\ref{pfq}) and  (\ref{eqpfrg}).\\
\noindent{\bf (3).} In plainer words, we must show that, under $\pf$, the evolution of the walk $X$ after time $\trar{j^v_n}$
is well approximated by a biased walk from the same location in the same environment whose moves are made independently of the past of $X$. The reason that the approximation is not exact is that, given its history until time $\trar{j^v_n}$, the walk $X$ must subsequently remain in the descendent tree of its location at the last constructed regeneration time $r_n^v$, for otherwise $r_n^v$ would not be a regeneration time for $X$. The reason that the approximation is good is that the walk at time $\trar{j^v_n}$ is typically already at a safe distance from 
$X(r_n^v)$.

That is, the conditional distribution of $X\big( \trar{j^v_n} + \cdot \big): \N \to T$ given $T$, 
the sequence $r_1^v < \trar{j_1^v} < \ldots < r_n^v < \trar{j_n^v}$ and
$X(0),\ldots,X(\trar{j^v_n})$, is given by $\P_{T,\beta}^{X(\trar{j^v_n})}$ conditional on  $Y(i) \not= X(r_n^v)$ for each $i \geq 0$. (Here, we write $Y$ for the walk under $\P_{T,\beta}^{X(\trar{j^v_n})}$. Note that we might also write $X_n'$ for this walk.) Hence, the total variation distance between this conditional distribution and that of $X'_n$ is the conditional probability that 
 $X'_n(i) = X(r_n^v)$ for some $i \geq 0$. It suffices then for the third part of the lemma to show that
$$
 \Big( \pf \Big)  \Big(  X'_n(i) = X(r_n^v) \,\, \textrm{for some $i \in \N$} \Big) \to 0,
$$
as $v \to \infty$. By shifting the origin of time to be $r_n^v$, we see that it is equivalent to show that
\begin{equation}\label{eqpfrgpr}
 \Big( \pf \Big)  \Big(  X'_1(i) = \phi(T) \,\, \textrm{for some $i \in \N$} \Big\vert 0 \in {\rm RG} \Big) \to 0,
\end{equation}
as $v \to \infty$.

To this end, we note that, for each $\ell \in \N$, 
\begin{equation}\label{pftdrg}
 \Big( \pf \Big) \Big(  d \big( \phi, X(\trar{j_1^v}) \big) \geq \ell  
   \Big\vert 0
 \in {\rm RG} \Big) \to 1
\end{equation}
as $v \to \infty$. This follows from  (\ref{eqpfrg}), (\ref{pfq}), 
and the almost-sure inequality $Q_v \leq j_1^v$. 

For any $v \in \bext(T)$, the walk $X$ under $\P_{T,\beta}^v$ will reach $\phi(T)$
only by travelling a distance of $d(\phi,v) - 1$ along the backbone $B(T)$,
from the parent of $v$ to $\phi$. Thus, for any such $v$,
\begin{equation}\label{ptxph}
\P_{T,\beta}^v \Big(  X(i) = \phi(T) \, \, \textrm{for some $i > 0$} \Big)
 \leq \binf^{1 - d(\phi,v)}.
\end{equation}
It follows from (\ref{pftdrg}) and (\ref{ptxph}) that, for given $k \in \N$, 
for any $\epsilon > 0$
and for  $v$ sufficiently high, 
the sample of tree and walk under $\big( \pf \big)  \big( \cdot \big\vert 0 \in {\rm RG} \big)$ 
and the random variable $j_1^v$ are, 
with probability at least $1 - \epsilon$, 
such that 
$$
\P_{T,\beta}^{\trar{j_1^v}} \Big( X(i) \not= \phi(T) \, \, \textrm{for all $i >
  0$} \Big) \geq 1 - \epsilon,
$$
as required for (\ref{eqpfrgpr}). \qed

We are now ready to provide the coupling that underlies the one in the statement of Lemma \ref{lemp}.
The coupling identifies
the local environments and walks about a lengthy initial sequence of entrances to high-weight traps
under the measures $\pf$ and $\idtrap$.
\begin{lemma}\label{lemhp}
For the statement, we adopt the following notation. Let $\Theta$ be a coupling of the laws $\pf$ and $\idtrap$.
We say that $\Theta$ is good if, for each $i \in \N$, the traps $T_{X(\trar{i})}$ and $\tilde{T}^i_\entrance$ are coupled to be equal by $\Theta$. Let $k \in \N$ and $v > 0$. A sample point of a good coupling $\Theta$ is said to   to be $k$-accurate for the $i$-th $v$-large trap if the following properties hold:
\begin{enumerate}
\item the descendent trees $T_{\xsub{\trar{j_i^v}}}$ and
  $\tilde{T}^{j_i^v}_\entrance$ coincide,
\item the $k$-large backbone neighbourhood $\bout_k\big( X(\trar{j_i^v}) \big)$ in $T$ of $X(\trar{j_i^v})$
coincides with its counterpart $\bout_k(\entrance)$ in $\big( \tilde{B}^{j_i^v}, \tilde{T}^{j_i^v} \big)$,
\item and the walks $X: \big\{ \trar{j_i^v}, \ldots, \sigma_{j_i}^{k}  \big\}
  \to T$ and $\tilde{X}_{j_i^v}: \big\{ 0, \ldots, \tilde{\sigma}_{j_i^v}^{k}  \big\}
  \to \big( \tilde{B}^{j_i^v}, \tilde{T}^{j_i^v} \big)$ coincide (their domains
  being identified by translation). 
\end{enumerate} 
For $\epsilon > 0$ and $\ell \in \N$, there exists $v_1$ such that, for $v > v_1$, the following holds.
A good coupling $\Theta = \Theta_{v,k,\ell,\epsilon}$
of  $\pf$ and $\idtrap$ exists such that the $\Theta$-probability of an outcome that is $k$-accurate for the $i$-th $v$-large trap for all $i \in \{1,\ldots,\ell\}$ is at least $1 - \epsilon$.
\end{lemma}
\noindent{\bf Remark.} Regarding the third of the listed properties, note that the first walk is taking values in the union of the vertices of 
 $T_{\xsub{\trar{j_i^v}}}$ and
 $\bout_k\big( X(\trar{j_i^v}) \big)$, 
while the second one does so in the union of the vertices of   $\tilde{T}^{j_i^v}_\entrance$ and  $\bout_k(\entrance)$. The first two properties ensure then a shared environment for these two walk segments, and the third property asserts that the walks themselves coincide.\\ 
\noindent{\bf Proof.} 
Note firstly that the sequences $\big\{ T_{X(\trar{i})}: i \in \N \big\}$
and  $\big\{ \tilde{T}^i_\entrance: i \in \N \big\}$ each have the distribution of an independent sequence of $\P_{h,\nu}$-distributed terms. We naturally choose a coupling that identifies the respective traps, and this is a good coupling.

It is a simple exercise to show from the construction in Definition \ref{defjiv} that
the descendent trees $T_{\xsub{\trar{j_i^v}}}$,
neighbourhoods $\bout_k\big( X(\trar{j_i^v}) \big)$ 
and walks $X: \big\{ \trar{j_i^v}, \ldots, \sigma_{j_i}^{k}  \big\}
  \to V\big(T_{X(\trar{j_i^v})}\big) \cup V\big(\bout_k\big( X(\trar{j_i^v})\big) \big)$
are independent as $i$ varies over $\N$. This is true of their $\idtrap$ counterparts, by definition. 
For this reason, it is enough to prove the statement in the case that $\ell = 1$. 
Regarding the three listed properties, 
the trees in the first have already been coupled. 
The neighbourhoods and walks in the second and third have laws differing in total variation norm by a quantity tending to zero as $v \to \infty$ by Lemma \ref{lemjd} parts (2) and (3), and by the definition of $\idtrap$. Hence,
the statement of the lemma indeed holds for $\ell = 1$. \qed

We now state several lemmas which will be used to prove Lemma \ref{lemp}. 
The coupling provided by Lemma \ref{lemhp} will be invoked to show the existence of the coupling in Lemma \ref{lemp}. The plan is to tune the value of $v$ in Definition \ref{defjiv} so that, among the
first $n$ traps, most of the time is spent in traps of weight exceeding $v$; 
$v$ will also be chosen to be high enough that these high-weight traps are well-separated, so that the cutoff of the $k$-neighbourhood
about the trap entrance used in the coupling of Lemma \ref{lemhp} does not err significantly in describing total trap time.  The choice that we will make is $v = c n^{1/\chiexp}$, for some small $c = c(\epsilon)$ (where $\epsilon > 0$ appears in the statement of Lemma \ref{lemp}).
As we now record, this choice results in a number of traps among the first $n$ whose weight exceeds the cutoff which is tight in $n$.
\begin{definition}\label{defjandk}
For $v > 0$, let $\largev^v$ denote the set of $i \in \N$ for which
$\omega \big( T_{X(\trar{i})} \big) > v$. (We will call the traps so identified $v$-large traps.)
Let $J_k^v$ denote the set $\big\{ j_i^v: i \in \N \big\}$ appearing in
Definition \ref{defjiv}.
\end{definition}
\begin{lemma}\label{lempintmin}
For any $c > 0$ and $\epsilon > 0$, there
exists $K \in \N$ such that, for $n$ sufficiently high,
\begin{equation}\label{pfin}
 \Big( \pf \Big) \Big( \Big\vert \largev^{c n^{1/\chiexp}} \cap \Big\{ 1,\ldots, n \Big\} \Big\vert \leq K \Big) > 1 - \epsilon.
\end{equation}
\end{lemma}
\noindent{\bf Proof.} This is implied by Proposition \ref{theoremone}. \qed 

We anticipate that the $cn^{1/\chiexp}$-large traps are each responsible for a significant portion of
the time to reach distance $n$. For this reason, it is important to show that 
the procedure from Definition \ref{defjiv} is unlikely to skip any such trap.
\begin{lemma}\label{lempinttwo}
For any $\epsilon > 0$, $k \in \N$, $c > 0$, and for all $n$
sufficiently high,
\begin{equation}\label{ppinf}
 \Big( \pf  \Big) \Big( \largev^{c n^{1/\chiexp}} \cap \big\{ 1,\ldots,n \big\} = 
   J_k^{c n^{1/\chiexp}} \cap \big\{ 1,\ldots,n  \big\} \Big) > 1 - \epsilon.
\end{equation}
\end{lemma}
\noindent{\bf Proof.}
Let $\big\{ k_i^v: i \geq  1 \big\}$ enumerate $\largev^v$.
For $c,\delta > 0$ and $n \in \N$, set
\begin{eqnarray}
 H_1^{c,\delta}(n) & = & \Big\{ k_1^{c n^{1/\chiexp}} > \delta n, \, \textrm{and
   $k_{j+1}^{c n^{1/\chiexp}} - k_j^{c n^{1/\chiexp}} > \delta n$} \label{hone} \\
  & &  \qquad \qquad \qquad \textrm{for each $j
   \in \N$ such that $k_{j+1}^{c n^{1/\chiexp}} \leq n$} \Big\}, \nonumber
\end{eqnarray}
and
\begin{eqnarray}
H_2^\delta(n) & = & \Big\{ \textrm{if $1 \leq i,j \leq n$ satisfy $j \geq i +
  \delta n$, then there exists $r \in {\rm RG}$} \nonumber \\
  & & \qquad \textrm{for which $\trar{i} < r < \trar{j}$} \Big\} \cap \Big\{ \exists r \in
  {\rm RG}: r < \trar{\delta n} \Big\}. \label{htwo}
\end{eqnarray}
 We claim that, for any $c > 0$, $\delta > 0$ and $k \in \N$, and
for all $n \in \N$ sufficiently high,
\begin{equation}\label{hinc} 
 H_1^{c,\delta}(n) \cap H_2^{\frac{\delta}{2k+3}}(n)  
 \subseteq  \Big\{ \largev^{c n^{1/\chiexp}} \cap \{ 1,\ldots, n \} = 
   J_{c n^{1/\chiexp}}^k \cap \{ 1,\ldots, n \}  \Big\}. 
\end{equation}
Once (\ref{hinc}) is demonstrated, the proof will finish by showing  
the event on its left-hand side to be highly probable.

One of the inclusions needed for the event on the right-hand side of (\ref{hinc}) is trivial: indeed,  $\ell \in J_k^v$ satisfies $\omen \big( T_{X(\trar{\ell})} \big) > v$, so that $J_k^{c n^{1/\chiexp}}
\subseteq \largev^{c n^{1/\chiexp}}$. To prove the reverse inclusion,
we adopt the shorthand $k_\ell = k_\ell^{c n^{1/\chiexp}}$ and $j_\ell = j_\ell^{c n^{1/\chiexp}}$ for $\ell \in \N$. We assume throughout this derivation that 
the event on the left-hand-side of (\ref{hinc}) occurs.

We will obtain (\ref{hinc}) by verifying the following assertion: for any given $\ell \in \N$, 
if  $k_\ell \leq n$, then $k_\ell = j_\ell$.
We now establish this statement by induction on $\ell$. 
For this, we will make use of the statement that, for each $j \in \N$ such that  $k_j \leq n$, 
\begin{equation}\label{equrgk}
  \Big\vert  \Big\{ \trar{k_{j-1}}, \ldots, \trar{k_j} \Big\}  \cap {\rm RG} \Big\vert \geq 2k + 2.
\end{equation}
Accepting (\ref{equrgk}) for now, we perform the inductive argument, arguing that, for any $\ell > 1$, 
(\ref{equrgk}), $k_{\ell-1} = j_{\ell-1}$ and $k_\ell \leq n$ imply that $k_\ell = j_\ell$. 
(The initial step, indexed by $\ell = 1$, is a minor modification of the generic one.)
Let $r'_1,\ldots,r'_{2k+2}$ denote the first $2k+2$ elements of ${\rm RG}$
occurring after time $\trar{k_{\ell-1}}$. By the definition of regeneration times, we find that
$d \big( \phi, X(r'_1) \big) \geq d \big( \phi, X(\trar{k_{\ell-1}}) \big)$, and also that the sequence $d \big( \phi, X(r'_i) \big)$ is strictly increasing in $i$. Hence, $r^{c n^{1/\chiexp}}_\ell$, by its definition (in Definition \ref{defjiv}) and by $k_{\ell-1} = j_{\ell-1}$,
is equal to one among $\big\{ r'_i: 1 \leq i \leq k +2 \big\}$. This implies that 
$d \big( X(r^{c n^{1/\chiexp}}_\ell), X(\trar{k_\ell})   \big) \geq  d \big( X(r^{c n^{1/\chiexp}}_\ell), X(r'_{2k+2}) \big) + 1 \geq k + 1$.
From this, we find that $k_\ell$ meets the condition required to ensure that  $j_\ell = k_\ell$. 

To obtain (\ref{hinc}), it remains to verify (\ref{equrgk}). To this end, note that 
the event $H_1^{c,\delta}(n)$
ensures that $k_j - k_{j-1} \geq \delta n$.
Hence,  the interval  $\big\{ \trar{k_{j-1}},\ldots, \trar{k_j}
\big\}$ may be divided into $2k+2$ consecutive subintervals, each of which contains 
at least $\delta n/(2k+2) - 1$
values of $\big\{ \trar{i}:i \in \N \big\}$. By the occurrence  $H_2^{\frac{\delta}{2k+3}}(n)$,
each of these subintervals intersects ${\rm RG}$. We obtain (\ref{equrgk}) and thus   (\ref{hinc}).

Next we show that the events on the left-hand side of (\ref{hinc}) are indeed probable.

Under $\pf$, elements of $\N$ belong independently to $\largev^v$ with probability $\P_{h,\nu}\big(\omega(T) > v\big)$.
As such, by Proposition \ref{theoremone}, 
we find that, given $\epsilon > 0$ and $c > 0$,
there exists $\delta > 0$ and $n_0 \in \N$, such that, for $n \geq n_0$,
\begin{equation}\label{pfho}
 \Big( \pf \Big) \Big( H_1^{c,\delta}(n)  \Big) > 1 - \epsilon.
\end{equation}  
Note that $H_2^\delta(n)$ equals $J(n)$ in the proof of Lemma \ref{lemt}. Restating (\ref{jnineq}),
there exists $C > 0$ such that, for $n$ sufficiently high,
\begin{equation}\label{pfht}
\Big( \P_{f,\nu,\infty} \times \P_{T,\beta}^{\phi} \Big) 
 \Big(  H_2^{\frac{\delta}{2k+3}}(n)  \Big) \geq 1 - C  \big( \log n \big)^{2 - \alpha}.
\end{equation}
By  (\ref{hinc}), (\ref{pfho}) and (\ref{pfht}), 
we obtain Lemma \ref{lempinttwo}. \qed

Making a step towards Lemma \ref{lempint}, we now show that, if we discard returns made to traps after long (of distance exceeding a large finite $k$) interludes on the backbone, the laws of total trap time among $c n^{1/\chiexp}$-large traps are  almost the same for $\pf$ and for $\idtrap$.
\begin{lemma}\label{lempint} 
For any $k \in \N$, $c > 0$,
there exists a sequence of couplings $\big\{ \coup_n: n \in \N \big\}$ 
of the distributions $\idtrap$
and $\pf$ such that, 
for any $\epsilon > 0$ and for $n$ sufficiently high,
\begin{equation}\label{ppth}
  \coup_n \Big( \sum_{i=1}^n \tau_i^{(k)} 1\!\!1_{\omen \big( T_{X(\trar{i})} \big)
 \geq c n^{1/\chiexp}} = \sum_{i=1}^n \tilde{\tau}_i^{(k)} 1\!\!1_{\tilde{\omega}_i 
 \geq c n^{1/\chiexp}}  \Big) > 1 -  \epsilon,
\end{equation}
and for which $\omega \big( T_{X(\trar{i})}\big) = \tilde\omega_i$ for all $i \in \N$ under $\Theta_n$.
\end{lemma}
\noindent{\bf Proof.}
We will couple $k$-cutoff total trap time for each pair of actual and idealized traps using Lemma \ref{lemhp}; for each pair, there is a small failure probability, but, by Lemma \ref{lempintmin}, the number of 
$c n^{1/\chiexp}$-large traps is uniformly bounded in $n$, so that these small probabilities remain small cumulatively.

We use 
Lemma \ref{lemhp} to define the 
coupling $\coup_n$ for each $n \in \N$: fixing $\epsilon > 0$,
we set $\coup_n$ equal to $\Theta_{cn^{1/\chiexp},k,K,\epsilon}$ in  Lemma \ref{lemhp}.  
By so doing, we obtain that there exists $n_0 \in \N$ such that, for $n \geq n_0$,
and for each $\ell \in \{1,\ldots,K\}$,
\begin{equation}\label{tte}
 \coup_n \bigg( \sum_{m=1}^\ell \tau_{j_m^{c n^{1/\chiexp}}}^{(k)} = 
\sum_{m=1}^\ell  \tilde{\tau}_{j_m^{c n^{1/\chiexp}}}^{(k)}  \bigg) > 1 - \epsilon.
\end{equation}
Note that, because $\Theta_n$ is a good coupling, the event in (\ref{ppth}) may be written
$\big\{ \sum_{i=1}^n \tau_i^{(k)} 1 \!\! 1_{i \in \largev^{c n^{1/\chiexp}}} = \sum_{i=1}^n \tilde\tau_i^{(k)} 1\!\!1_{i \in \largev^{c n^{1/\chiexp}}} \big\}$. Hence, Lemma \ref{lempint} follows from 
(\ref{ppinf}),  (\ref{pfin}) and (\ref{tte}), with a relabelling of $\epsilon > 0$. \qed \\
\noindent{\bf Proof of Lemma \ref{lemp}.}
With Lemma \ref{lempint} available, to prove Lemma \ref{lemp}, it remains only to show that
small traps count for little time under the actual and idealized laws, and that the $k$-cutoff creates an error only with small probability for each of these laws. We state four assertions to this effect and then prove them.

Firstly, for any $\delta > 0$, there exists $c > 0$ and $n_0 \in \N$ such
that, for $n \geq n_0$,
\begin{equation}\label{ppo}
 \Big( \pf \Big) \Big( \sum_{i=1}^n \tau_i 1\!\!1_{\omen\big( T_{X(\trar{i})} \big)
 < c n^{1/\chiexp}} > \delta  \sum_{i=1}^n \omen\big( T_{X(\trar{i})} \big) \Big) < \delta.
\end{equation} 
Secondly, for any $\delta > 0$, there exists $c> 0$ and $n_0 \in \N$ such
that, for $n \geq n_0$,
\begin{equation}\label{pptw}
  \idtrap \Big( \sum_{i=1}^n \tilde{\tau}_i 1\!\!1_{\tilde{\omega}_i 
 < c n^{1/\chiexp}} > \delta  \sum_{i=1}^n \tilde{\omega}_i      \Big) < \delta.
\end{equation} 
Thirdly,  for any $\epsilon > 0$ and $c> 0$, there exists $k \in \N$ and  $n_0 \in \N$ such
that, for $n \geq n_0$,
\begin{equation}\label{ppf}
  \Big( \pf  \Big) \Big( \sum_{i=1}^n \tau_i^{(k)} 1\!\!1_{\omen\big( T_{X(\trar{i})} \big)
  > c n^{1/\chiexp}}  
= \sum_{i=1}^n \tau_i 1\!\!1_{\omen \big( T_{X(\trar{i})}\big)
  > c n^{1/\chiexp}}  \Big) > 1 - \epsilon.
\end{equation}
Finally, for any $\epsilon > 0$ 
 and $c> 0$, there exists $k \in \N$ and  $n_0 \in \N$ such
that, for $n \geq n_0$,
\begin{equation}\label{qtf}
  \idtrap \Big( \sum_{i=1}^n \tilde{\tau}_i^{(k)} 1\!\!1_{\tilde{\omega}_i
  > c n^{1/\chiexp}}  
= \sum_{i=1}^n \tilde{\tau}_i 1\!\!1_{\tilde{\omega}_i
  > c n^{1/\chiexp}}  \Big) > 1 - \epsilon.
\end{equation} 
Deferring the four derivations, we close out the proof of the lemma.
Note that (\ref{ppth}),
(\ref{ppf}) and (\ref{qtf})
imply that, for any $\epsilon > 0$ and $c > 0$, and for all $n$
sufficiently high,
$$
  \coup_n \Big( \sum_{i=1}^n \tau_i 1\!\!1_{\omen \big( T_{X(\trar{i})} \big)
 \geq c n^{1/\chiexp}} = \sum_{i=1}^n \tilde{\tau}_i 1\!\!1_{\tilde{\omega}_i 
 \geq c n^{1/\chiexp}}  \Big) > 1 -  3\epsilon,
$$
to which the use of (\ref{ppo}) and (\ref{pptw}) yields, for any $\epsilon
> 0$ and $\delta > 0$, and for all $n$ sufficiently high,
$$
  \coup_n \Big( \Big\vert \sum_{i=1}^n \tau_i  - \sum_{i=1}^n
  \tilde{\tau}_i 
  \Big\vert < \delta  \sum_{i=1}^n \omen \big( T_{X(\trar{i})} \big)  + \delta \sum_{i=1}^n
  \tilde{\omega}_i  \Big) > 1 -  3\epsilon - 2\delta.
$$
The coupling $\coup_n$ having been constructed so that
 $\tilde{\omega}_i =  \omen \big( T_{X(\trar{i})} \big)$, we obtain the
statement of the lemma (with $\epsilon$ replaced by $c$), by choosing $3 \epsilon + 2 \delta <
c$. 

Turning to the proofs of the four statements,  for (\ref{ppf}), note that, for some $c_1 \in (0,1)$ 
and for each $i \in \N$,
\begin{equation}\label{tauikt}
 \Big( \pf \Big) \Big( \tau_i^{(k)} = \tau_i \Big) \geq 1 - c_1^k.
\end{equation}
Indeed, recalling Definition \ref{deftauu}, given any backbone-tree pair $\big( B,\singent \big)$, 
$\singent =  (\head,\entrance)  \circ T_{\entrance}$,
we have that
$\tau_{T_{\entrance}}^k \not= \tau_{T_{\entrance}}$ under $\P_{(B,\singent),\beta}^\head$ precisely when 
$X$ returns to $\entrance$ after
time $\sigma_{(B,\singent)}^{\head}(k)$, an event whose probability is $\ekm(B,\singent) -
\pesc(B,\singent)$.
By Lemma \ref{lembt}, then, we derive (\ref{tauikt}).

Note that (\ref{tauikt}) and (\ref{pfin}) imply that
$$
 \Big( \pf \Big) \Big( \tau_i^{(k)} = \tau_i \, \textrm{for each $i \in
   \largev^{c n^{1/\chiexp}} \cap \big\{ 1,\ldots,n \big\}$} \Big) \geq 1 - \big( K c_1^k   + \epsilon \big),
$$
whence (\ref{ppf}). Note that (\ref{qtf}) follows similarly.

To derive (\ref{ppo}), note that, for any $c_1 > 0$, there exists $c > 0$
and $n_0 \in \N$ such that, for $n \geq n_0$,
\begin{equation}\label{pomc}
\Big( \pf \Big) \Big( \sum_{i=1}^n \omen\big( T_{X(\trar{i})} \big)
1\!\!1_{\omen\big( T_{X(\trar{i})} \big) < c n^{1/\chiexp}} \leq c_1  \sum_{i=1}^n
\omen\big( T_{X(\trar{i})} \big)  \Big) > 1 - c_1,
\end{equation}
by Proposition \ref{theoremone} and the assumption that $\chiexp < 1$. Note that,
by (\ref{eqnmean}) and Lemma \ref{lemebfbd}, there exists $C > 0$ such that,
for any $i \in \N$,
$$
\E_{\pf} \Big( \tau_i \Big\vert \omen\big( T_{X(\trar{i})} \big) \Big)
\leq C \omen\big( T_{X(\trar{i})} \big).
$$
Hence, by Markov's inequality, 
\begin{equation}\label{pomct}
\Big( \pf \Big) \Big( \sum_{i=1}^n \tau_i 
1\!\!1_{\omen\big( T_{X(\trar{i})} \big) < c n^{1/\chiexp}} >  C^2  \sum_{i=1}^n
\omen\big( T_{X(\trar{i})} \big) 
1\!\!1_{\omen\big( T_{X(\trar{i})} \big) < c n^{1/\chiexp}}  \Big) \leq C^{-1}.
\end{equation}
Choosing $C$, $c$ and $c_1$ in (\ref{pomc}) and (\ref{pomct}) such that 
$C^2 c_1 < \delta$ and  $\max \big\{ c,c_1, C^{-1} \big\} <  \delta/2$,
we obtain (\ref{ppo}). 

The statement (\ref{pptw}) has the same proof, the measure $\idtrap$
replacing $\pf$. This completes the proof.  \qed
\end{subsection}
\begin{subsection}{Completing the proof of Theorem \ref{thm}}
The time $\overh{n}$ to reach a vertex on $B$ at distance $n$, and the hitting
time $\hitting{n}$ appearing in Theorem \ref{thm}, are almost the same:
\begin{lemma}\label{lemdnhn}
For any $\epsilon > 0$, 
$\overh{n(1-\epsilon)}  \leq   \hitting{n}  \leq \overh{n}$ holds for all $n$ sufficiently high, 
$\pf$-a.s.
\end{lemma} 
\noindent{\bf Proof.}
We have $\hitting{n} \leq \overh{n}$ by definition.  The inequality $\overh{n(1-\epsilon)}  >   \hitting{n}$
entails that $T_{X(\trar{i})}$ has depth at least $\epsilon n$ for some $i
\leq \tnb{n(1 - \epsilon)}$. 
 The probability of this event is negligible, by
Lemma \ref{lemf} and the elementary bound 
\begin{equation}\label{eqkol}
\P_h \big( D(T) = n \big) \leq m_h^n
\end{equation}
where recall that $m_h \in (0,1)$ denotes the expected number of children of the
offspring law $h$. \qed \\
\noindent{\bf Proof of Theorem \ref{thm}.} 
By Lemma \ref{lemdnhn}, it suffices for (\ref{statone}) to prove the analogous statement for $\overh{n}$.
In this regard, note that
\begin{eqnarray}
 & & \Big( \pf \Big) \Big( \overh{n} \geq u \Big)
 \geq \Big( \pf \Big) \Big( \sum_{i=1}^{\constt n(1-\epsilon)} \tau_i \geq u \Big) - \epsilon \nonumber \\
 & \geq & \idtrap  \Big( \sum_{i=1}^{\constt n(1-\epsilon)} \tilde\tau_i \geq u + \epsilon \sum_{i=1}^{\constt n(1-\epsilon)}  
 \tilde\omega_i \Big) - 2\epsilon. \label{taillbd}
\end{eqnarray}
The first inequality here is due to Lemma \ref{lemt}, while the second makes use of the coupling $\Theta_{\constt n(1-\epsilon)}$ in Lemma \ref{lemp}, including the property that $\omega \big( T_{X(\trar{i})}\big) = \tilde\omega_i$.

Under $\idtrap$, the terms $\tilde\omega_i$ are independent and distributed as $\omega(T)$ under $\P_{h,\nu}$.
Hence, Proposition \ref{theoremone} and $\chiexp < 1$ implies that $\sum_{i=1}^n  
 \tilde\omega_i \leq C n^{1/\chiexp}$ with high probability. However, the terms $\tilde\tau_i$ are also independent under $\idtrap$, so that Theorem \ref{proptraplaw} implies that it is also probable that
 $\sum_{i=1}^n  
 \tilde\tau_i \geq c n^{1/\chiexp}$. Hence, for all $\epsilon > 0$, there is some $c > 0$ such that, for $n$ sufficiently high,
$$
\idtrap \Big( \sum_{i=1}^n \tilde\tau_i \geq c \sum_{i=1}^n  
 \tilde\omega_i \Big) > 1 - \epsilon.
$$
Using this in (\ref{taillbd}), we learn that, for any $\epsilon > 0$ and all $n$ sufficiently high,
\begin{equation}\label{qlowbd}
\Big( \pf \Big) \Big( \overh{n} \geq u \Big) \geq  \idtrap  \Big( \sum_{i=1}^{\constt n(1-\epsilon)} \tilde\tau_i \geq u (1 +  \epsilon)  \Big)  - \epsilon.
\end{equation}
The conclusion that for any $\epsilon > 0$ and all $n$ sufficiently high,
\begin{equation}\label{qupbd}
 \Big( \pf \Big) \Big( \overh{n} \geq u \Big) \leq  \idtrap  \Big( \sum_{i=1}^{\constt n(1 + \epsilon)} \tilde\tau_i \geq u (1 -  \epsilon)  \Big)  + \epsilon,
\end{equation}
may be reached by using the same ingredients, as well as the fact that for all $C > 0$, $\idtrap \big( \sum_{i=1}^n \tilde\tau_i \geq C n \big) \to 1$, as $n \to \infty$ (which follows from Theorem \ref{proptraplaw} and $\chiexp < 1$).
Now, by a direct computation that uses Theorem \ref{proptraplaw} and the independence of $\big\{ \tilde\tau_i: i \in \N \big\}$ under $\idtrap$, 
$$
\E_{\idtrap} \exp \Big\{ - \frac{\lambda}{n^{1/\chiexp}} \sum_{i=1}^{\constt n} \tilde\tau_i \Big\}
 \to \exp \Big\{ - \xi \lambda^{\chiexp} \Big\}
$$
for any $\lambda > 0$, from which it follows that, under $\idtrap$,
$$
 \frac{1}{n^{1/\chiexp}} \sum_{i=1}^{\constt n} \tilde\tau_i  \rightarrow \xi^{1/\chiexp} S_{\chiexp} \quad \textrm{in distribution as $n \to \infty$,}
$$
where $S_\chiexp$ is the stable law of index $\chiexp$.

We may then use (\ref{qlowbd}) and (\ref{qupbd}) to reach the conclusion that (\ref{statone}) holds for $\overh{n}$ under $\pf$. As we noted earlier, this proves (\ref{statone}) itself.

To derive (\ref{stattwo}), define $\overline{M}:\N \to \N$ under $\pf$ according to 
$\overline{M}_n = \sup \big\{ d \big(\phi, X_i \big) : 0 \leq i \leq n \big\}$ for each $n \in \N$. Note that (\ref{statone}) may be recast as the statement that
$n^{-\chiexp} \overline{M}_n$ converges to $\xi^{-1}  S_\chiexp^{-\chiexp}$ in distribution. 
Write $\overline{X}_n =  d \big(\phi, X_i \big)$ for each $n \in \N$. To obtain (\ref{stattwo}), it remains only to argue that, for each $\epsilon > 0$,
\begin{equation}\label{eqfinal}
\big( \pf \big) \Big( \big\vert \overline{X}_n  - \overline{M}_n \big\vert > \epsilon  \overline{M}_n \Big) \to 0
\quad
\textrm{as 
$n \to \infty$.}
\end{equation}
As the ensuing sketch indicates, the proof of this statement is straightforward.  
We note that, by Lemma \ref{lemf} and (\ref{eqkol}), 
there is a negligible probability that the walk encounters a trap of depth at least 
$\epsilon \overline{M}_n$ by time $n$; 
and, 
by virtue of a brief argument whose main tool is Lemma \ref{lemrgexp}, there is also a negligible probability that the walk backtracks along the backbone by a distance of $\epsilon \overline{M}_n$ before time $n$. Hence, we obtain (\ref{eqfinal}) and complete the proof of (\ref{stattwo}). 
 
The final statement in the theorem, regarding the case $\chiexp > 1$, follows from Lemmas \ref{lemdnconv} and \ref{lemdnhn}. \qed
\end{subsection}
\begin{subsection}{The environment about the trap entrance at late time}\label{secenvproof}
We overview a proof of the assertion made in Section \ref{secenv}, that the measure $\hat\rho$ is the environment viewed from the entrance of the trap to which the walker under $\pf$ belongs at late time. 
\begin{definition}
Recall the idealized sequence $\idtrap$ of $\qstarnoarg$-distributed backbone-tree pairs 
$\big( \tilde{B}^i, \tilde{T}^i \big)$ from Definition \ref{defidseq}.
Under the law $\idtrap$, we define the idealized walker $X_{\rm id}$, taking values in the union of the listed backbone-tree pairs. This walker begins at time zero at the entrance of the first such pair
$\big( \tilde{B}^1, \tilde{T}^1 \big)$, and then follows a trajectory in that graph according to the walk law
 $\P_{(\tilde{B}^1,\tilde{T}^1),\beta}^\entrance$ until its final visit to the entrance. At the next integer time, the walk is transported to the entrance of 
$\big( \tilde{B}^2, \tilde{T}^2 \big)$, and begins moving according to the law 
 $\P_{(\tilde{B}^2,\tilde{T}^2),\beta}^\entrance$. The walker successively visits each 
$\big( \tilde{B}^i, \tilde{T}^i \big)$ in this way.
\end{definition}
In this section, we coupled the laws $\pf$ and $\idtrap$ together, to derive such results as Lemma \ref{lemp}.
With minor modifications to this construction, we may couple the two laws in such a way that, with probability tending to one as time tends to infinity, there is agreement in the two marginals
in all of the following respects: the trap in which the walker lies at present, the neighbourhood $\bout_k$ about the trap entrance (for fixed $k$), and the position of the walker inside the trap.  
This coupling reduces the question to understanding the limiting distribution of the environment about the trap entrance 
for the idealized walker $X_{\rm id}$(n) as $n \to \infty$.
The total time spent in the trap to which $X_{\rm id}(n)$ belongs being likely to be high, it is enough to find the conditional distribution about the trap entrance of an idealized backbone-tree pair (having the law $\qstarnoarg$)  conditional on the total time spent in the trap being high. Denoting this total time by $\tau$, a slight reworking of Lemma \ref{lemsmallerr} provides the approximate representation $\tau \approx E 1\!\!1_\deepfall$. Lemma \ref{lemrhox} states that  the backbone-tree pair law $\rho$ arises by conditioning $\qstarnoarg$ on $E = x$ and then taking $x \to \infty$. This limit taken, we have further to condition on the occurrence of $\deepfall$ 
in order to obtain the limiting distribution of $\qstarnoarg$ conditionally on $\tau = x$ with $x \to \infty$.
This is the reconditioning of $\rho$ that appears in the definition of $\hat\rho$. \qed 
\end{subsection}
\end{section}
\bibliographystyle{plain}
\bibliography{trwbib}
\end{document}